\documentclass[12pt,dvips,psamsfonts,leqno,oneside,letterpaper]{amsart}
\usepackage[letterpaper,left=1.75truein, top=1truein]{geometry}
\usepackage{amssymb,amsmath,amscd}
\usepackage[dvips]{graphicx}
\usepackage[colorlinks,linkcolor=blue,citecolor=blue]{hyperref}

\newtheorem{thm}{Theorem}[section]
\newtheorem{cor}[thm]{Corollary}
\newtheorem{lemma}[thm]{Lemma}
\newtheorem{prop}[thm]{Proposition}
\newtheorem{conj}[thm]{Conjecture}

\newtheorem{obs}[thm]{Observation}

\theoremstyle{definition}
\newtheorem{defn}[thm]{Definition}
\newtheorem{defs}[thm]{Definitions}
\newtheorem{con}[thm]{Conventions}

\newtheoremstyle{cases}
  {12pt plus 6 pt}
  {2pt}
  {\bfseries}   
  {}
  {\bfseries}
  {.}
  {.5em}
  {}
\theoremstyle{cases}
\newtheorem{ccase}{Case}
\newtheorem{case}{Case}
\newtheorem{subcase}{Subcase}
\newtheorem{subsubcase}{Subsubcase}
\numberwithin{subcase}{case}
\numberwithin{subsubcase}{subcase}
\numberwithin{equation}{thm}

\newcommand\G{{\Gamma}}

\newcommand\e{{\epsilon}}

\newcommand\n{{\nu}}

\newcommand\s{{\sigma}}
\newcommand\p{{\partial}}
\newcommand\ra{{\rightarrow}}
\newcommand\lra{{\longrightarrow}}

\newcommand\D{{\Delta}}

\newcommand\z{{\mathbb Z}}

\begin{document}
\title{Characteristic Subsurfaces, Character Varieties and Dehn Fillings}
\author{Steve Boyer, Marc Culler, Peter B.~Shalen and Xingru Zhang}

\address{
D\'epartement de math\'ematiques\\
Universit\'e du Qu\'ebec \`a Montr\'eal
P. O. Box 8888, Postal Station Centre-ville
Montr\'eal, Qc, H3C 3P8, Canada}
\email{boyer@math.uqam.ca}

\address{
Department of Mathematics, Statistics, and Computer Science (M/C 249)\\
University of Illinois at Chicago\\
851 S. Morgan St.\\
Chicago, IL 60607-7045}
\email{culler@math.uic.edu}

\address{
Department of Mathematics, Statistics, and Computer Science (M/C 249)\\
University of Illinois at Chicago\\
851 S. Morgan St.\\
Chicago, IL 60607-7045, USA}
\email{shalen@math.uic.edu}

\address{
Department of Mathematics\\
SUNY at Buffalo\\
Buffalo, NY, 14260-2900, USA}
\email{xinzhang@math.buffalo.edu}

\thanks{Steve Boyer is partially  supported by NSERC grant OGP0009446
and FCAR grant ER-68657.}
\thanks{Marc Culler is partially supported by NSF grants DMS-0204142
and DMS-0504975.}
\thanks{Peter Shalen is partially supported by NSF grants DMS-0204142
and DMS-0504975.}
\thanks{Xingru Zhang is partially supported by NSF grant DMS-0204428.}

\begin{abstract}
  Let $M$ be a one-cusped hyperbolic manifold.  A slope on the
  boundary of the compact core of $M$ is called exceptional if the
  corresponding Dehn filling produces a non-hyperbolic manifold.  We
  give new upper bounds for the distance between two exceptional
  slopes $\alpha$ and $\beta$ in several situations.  These include
  cases where $M(\beta)$ is reducible and where $M(\alpha)$ has
  finite $\pi_1$; or $M(\alpha)$ is very small; or $M(\alpha)$ admits a
  $\pi_1$-injective immersed torus.
\end{abstract}

\maketitle

\section{Introduction}

Throughout this paper $M$ will denote a compact, connected,
orientable, irreducible, atoroidal $3$-manifold whose boundary is a
torus.  Thus $M$ is homeomorphic to the compact core of a
finite-volume hyperbolic $3$-manifold with one cusp.  A slope $\alpha$
on $\partial M$ (defined in Section \ref{notation}) is said to be {\it
  exceptional} if the Dehn filling $M(\alpha)$ does not admit a
hyperbolic structure.  By the {\it distance} between two slopes
$\alpha$ and $\beta$ we will mean their geometric intersection number
$\Delta(\alpha,\beta)$.

Cameron Gordon has conjectured in [Go] that the distance between any
two exceptional slopes for $M$ is at most $8$, and also that there are
exactly four specific manifolds $M$ which have a pair of exceptional
slopes with distance greater than $5$.  The results in this paper give
upper bounds for the distance between two exceptional slopes in several
special cases.  We assume for most of these results that $M(\beta)$ is
reducible, and that $M(\alpha)$ is a non-hyperbolic manifold of one of
several types.  Here, and throughout the paper, we will write $L_p$ to
denote a lens space whose fundamental group has order $p \geq 2$.

Our first result applies in the case that $M(\alpha)$ has finite
fundamental group.

\begin{thm} \label{redfin}
  If $M(\beta)$ is reducible and if $\pi_1(M(\alpha))$ is finite then
  $\Delta(\alpha,\beta)\leq 2$.  Moreover, if $\Delta(\alpha, \beta) =
  2$, then $H_1(M) \cong \mathbb Z \oplus \mathbb Z/2$, $M(\beta) =
  L_2 \# L_3$ and $\pi_1(M(\alpha)) \cong O_{24}^* \times \mathbb
  Z/j$, where $O_{24}^*$ denotes the binary octahedral group.
\end{thm}
 
Although we expect that the case $\Delta(\alpha,\beta) = 2$ does not
arise, this theorem is a considerable improvement on the previously
known bounds [BZ2].

Recall that a closed $3$-manifold $N$ is said to be {\it very small}
if $\pi_1(N)$ has no non-Abelian free subgroup.  The next result deals
with the situation where $M(\beta)$ is reducible and $M(\alpha)$ is
very small. The proof is based on an analysis of the $PSL_2(\mathbb C)$
character variety of a free product of cyclic groups.  (See Section
\ref{notation} for the definition of a strict boundary slope.)

\begin{thm} \label{redvsmall2}
  Suppose that $M(\beta)$ is a reducible manifold and $\beta$ is a
  strict boundary slope.  If $M(\alpha)$ is very small
  then $\Delta(\alpha, \beta) \leq 3$.
\end{thm}

A closed orientable $3$-manifold $N$ is said to admit a {\it geometric
  decomposition} if the pieces of its prime and torus decompositions
either admit geometric structures or are $I$-bundles over the torus.
According to Thurston's Geometrization Conjecture, which has been
claimed by Perelman, any closed orientable $3$-manifold admits a
geometric decomposition.  If we add the information that $M(\alpha)$
admits a geometric decomposition we obtain the following stronger
result.

\begin{thm} \label{redvsmall}
  Suppose that $M(\beta)$ is a reducible manifold and $M(\alpha)$ is a
  very small manifold that admits a geometric decomposition, then
  $\Delta(\alpha, \beta) \leq 2$.
\end{thm}

The next result applies in the case where $M(\alpha)$ contains an
immersed $\pi_1$-injective torus.  Note that in this case, 
$M(\alpha)$ is either reducible, toroidal, or a Seifert fibred space with base
orbifold of the form $S^2(r,s,t)$ (see Torus Theorem in [Sc] and Corollary 8.3 of [Ga3]). 
The bound $\Delta(\alpha, \beta) \leq 3$ holds in the first
two cases by [GLu], [Oh], [Wu2]. Thus the new information contained in
this theorem concerns the case where $M(\alpha)$ is Seifert fibred and geometrically atoroidal. 

\begin{thm}  \label{redtor}
  Suppose that $\beta$ is a strict boundary slope for $M$. If
  $M(\beta)$ is a reducible manifold and if $M(\alpha)$ admits a
  $\pi_1$-injective immersed torus, then $\Delta(\alpha,\beta)\leq 4$.
  Moreover, if $\Delta(\alpha, \beta) = 4$, then $M(\alpha)$ is a
  Seifert-fibred manifold with base orbifold $S^2(r,s,t)$, where
  $(r,s,t)$ is a hyperbolic triple and at least one of $r$, $s$ or $t$
  is divisible by $4$.
\end{thm}

The inequalities we obtain in the last two results are significantly
sharper than those obtained under comparable hypotheses in [BCSZ]. For
Theorem \ref{redtor}, this is due to the fact that in [BCSZ] it is
only assumed that $\beta$ is the boundary slope of an essential,
planar surface in $M$.  Here we are using additional information about
the topological structure of the connected sum decomposition of
$M(\beta)$.

Since a Seifert fibred manifold is either very small or contains a
$\pi_1$-injective immersed torus, the results above immediately
yield the following corollary.

\begin{cor} \label{redseif}
  If $M(\beta)$ is a reducible manifold, $\beta$ is a strict boundary
  slope, and $M(\alpha)$ is Seifert fibred, then
  $\Delta(\alpha,\beta)\leq 4$. If $\Delta(\alpha, \beta) = 4$, then
  the base orbifold ${\mathcal B}$ of $M(\alpha)$ is $S^2(r,s,t)$
  where $(r,s,t)$ is a hyperbolic triple and $4$ divides at least one
  of $r,s,t$.
\end{cor}

We also obtain the following result in the case where $M(\beta)$ is
only assumed to be non-Haken rather than reducible.

\begin{thm} \label{planar}
  If $\beta$ is a strict boundary slope and $M(\beta)$ is not a Haken
  manifold, then
\begin{enumerate}
\item[(1)] $\Delta(\alpha,\beta)\leq 2$ if $M(\alpha)$ has finite
  fundamental group;
\item[(2)] $\Delta(\alpha,\beta)\leq 3$ if $M(\alpha)$ is very small;
\item[(3)] $\Delta(\alpha,\beta)\leq 4$ if $M(\alpha)$ admits a
  $\pi_1$-injective immersed torus.
\end{enumerate}
\end{thm}

We will show that our results imply the following restricted
version of Gordon's conjecture.
 
\begin{thm} \label{redexc} 
  If $M(\beta)$ is a reducible manifold and $\beta$ is a strict
  boundary slope, then $M(\alpha)$ is a hyperbolic manifold for any
  slope $\alpha$ such that $\Delta(\alpha, \beta) > 5$. If we assume
  that the geometrization conjecture holds, then $M(\alpha)$ is a
  hyperbolic manifold for any slope $\alpha$ such that $\Delta(\alpha,
  \beta) > 4$.
\end{thm}

We remark that we expect the following to hold in this subcase of
Gordon's Conjecture.

\begin{conj}
  If $M(\beta)$ is a reducible manifold, then $M(\alpha)$ is a
  hyperbolic manifold for any slope $\alpha$ such that $\Delta(\alpha,
  \beta) > 3$.
\end{conj}

The bound in the conjecture cannot be lowered. For instance, if $M$ is the hyperbolic manifold obtained by doing a Dehn filling of slope $6$ on one boundary component of the Whitehead link exterior, then $M(1) \cong L_2 \# L_3$ is reducible while $M(4)$ is toroidal.  

The paper is organized as follows.  Basic definition and notational
conventions are given in Section \ref{notation}.  We review the notion
of a singular slope in Section \ref{singularslopes} and prove
Proposition \ref{sing summary}, which characterizes the situations in
which a boundary slope can fail to be a singular slope.  At the end of
Section \ref{singularslopes} we prove Theorems \ref{planar} and
\ref{redexc}, assuming Theorems \ref{redfin}, \ref{redvsmall} and
\ref{redtor}.  Section \ref{special singular} contains the proof of a
technical result (Proposition \ref{strictsinggen}) about singular
slopes in $L(p,1)\#L(q,1)$ which is stated and applied earlier, in Section
\ref{singularslopes}.  In Section \ref{preliminary} we reduce the
proofs of Theorems \ref{redfin} -- \ref{redtor} to more specific
propositions, which are proved in 
Sections \ref{PTLF},
\ref{weakPTLVS}, \ref{PTLVS} and \ref{PTLS} respectively.
Section \ref{background} is a review of $PSL_2(\mathbb
C)$-character variety theory and Section \ref{very small reps} contains
results about the representation varieties of fundamental groups of very
small $3$-manifolds.  Section \ref{CharacterSurfaces} is based on the
Characteristic Submanifold methods used in [BCSZ], and extends some of
those results under the additional topological assumptions that are
available in the setting of this paper.  These results are applied
in Section \ref{PTLS}. 

\section{Notation and Definitions}\label{notation}

We will use the notation $|X|$ to denote the number of components of a
topological space $X$.  The first betti number of a space $X$ will be
denoted $b_1(X)$.  

By a {\it lens space} we mean a closed orientable $3$-manifold with a
genus $1$ Heegaard splitting.  A lens space will be called non-trivial
if it is not homeomorphic to $S^2\times S^1$ or $S^3$.

By an {\it essential surface} in a compact, orientable $3$-manifold,
we mean a properly embedded, incompressible, orientable surface such
that no component of the surface is boundary-parallel and no
$2$-sphere component of the surface bounds a $3$-ball.

A {\it slope} $\alpha$ on $\partial M$ is a pair $\{\pm a\}$ where $a$ is a
primitive class in $H_1(\partial M)$.  The manifold $M(\alpha)$ is
the Dehn filling of $M$ obtained by attaching a solid torus to the
boundary of $M$ so that the meridian is glued to an unoriented curve
representing the classes in $\alpha$.

\begin{defn} A slope
  $\beta$ on $\p M$ is called a {\it boundary slope} if there is an
  essential surface $F$ in $M$ such that $\partial F$ is a non-empty
  set of parallel, simple closed curves in $\partial M$ of slope
  $\beta$. In this case we say that $F$ has slope $\beta$.  If
  $M(\beta)$ is reducible then of course $\beta$ is a boundary slope.
\end{defn}

Next consider a connected surface $F$ properly embedded in a $3$-manifold
$W$ with bicollar $N(F) = F \times [-1,1]$ in $W$. Denote by $W_{F}$
the manifold $W \setminus F \times (-\frac12, \frac12)$ and set $F_+ =
F \times \{\frac12\}, F_- = F \times \{-\frac12\} \subset \partial
W_F$.  We say that $W$ {\it fibres over $S^1$ with fibre $F$} if $W_F$
is connected and $(W_F, F_+ \cup F_-)$ is a product $(I, \partial
I)$-bundle pair. We say that $W$ {\it semi-fibres over $I$ with
  semi-fibre $F$} if $W_F$ is not connected and $(W_F, F_+, \cup F_-)$
is a twisted $(I, \partial I)$-bundle pair. 

\begin{defn}
  A slope $\beta$ on $\partial M$ is called a {\it strict boundary
    slope} if there is an essential surface $F$ in $M$ of slope
  $\beta$ which is neither a fibre or a semi-fibre.
\end{defn}

\begin{defn}
  Given a closed, essential surface $S$ in $M$, we let ${\mathcal
    C}(S)$ denote the set of slopes $\delta$ on $\partial M$ such that $S$
  compresses in $M(\delta)$.  A slope $\eta$ on $\partial M$ is called a
  {\it singular slope} for $S$ if $\eta \in {\mathcal C}(S)$ and
  $\Delta(\delta, \eta) \leq 1$ for each $\delta \in {\mathcal C}(S)$.
\end{defn}

\section{Reducible Dehn Fillings and Singular Slopes}\label{singularslopes}

A fundamental result of Wu [Wu1] states that if ${\mathcal C}(S) \ne
\emptyset$, then there is at least one singular slope for $S$.

The following result, which links singular slopes to exceptional
surgeries, is due to Boyer, Gordon and Zhang.

\begin{prop}\label{bgz}
  {\rm([BGZ, Theorem 1.5])} If $\eta$ is a singular slope for some
  closed essential surface $S$ in $M$ then for an arbitrary slope
  $\alpha$ we have
$$\Delta(\alpha, \eta) \leq
\left\{
\begin{array}{ll}
1 & \mbox{if $M(\alpha)$ is either small or reducible} \\
1 & \mbox{if $M(\alpha)$ is Seifert fibred and $S$ does not separate}\\
2 & \mbox{if $M(\alpha)$ is toroidal and ${\mathcal C}(S)$ is infinite}\\
3 & \mbox{if $M(\alpha)$ is toroidal and ${\mathcal C}(S)$ is finite}.\\
\end{array} \right.$$
Consequently if $M(\alpha)$ is not hyperbolic, then
$\Delta(\alpha, \eta) \leq 3$.
\end{prop}

If $b_1(M)\ge 2$ and $M(\beta)$ is reducible, then work of Gabai [Ga1,
Corollary] implies that $\beta$ is a singular slope of some closed,
essential surface.  This is also true generically when $b_1(M) = 1$,
as the following result indicates.

\begin{thm} {\rm (Theorem 2.0.3 and Addendum 2.0.4  of [CGLS])} \label{2.0.3}
  Suppose that $b_1(M) = 1$ and that $\eta$ is a boundary slope on
  $\partial M$. Then one of the following possibilities holds.
\begin{enumerate}
\item[(1)] $M(\eta)$ is a Haken manifold.
\item[(2)] $M(\eta)$ is a connected sum of two non-trivial lens
  spaces.
\item[(3)] $\eta$ is a singular slope for some closed essential
  surface in $M$.
\item[(4)] $M(\eta) \cong S^1 \times S^2$ and $\eta$ is not a strict
  boundary slope.
\end{enumerate}
\end{thm}
\qed

Thus when $M(\beta)$ is reducible, either $\beta$ is a singular slope
for some closed, essential surface in $M$; or $M(\beta)$ is $S^1 \times
S^2$ and $\beta$ is not a strict boundary slope; or $M(\beta)$ is a
connected sum of two lens spaces. In particular the inequalities of
Proposition \ref{bgz} hold unless, perhaps, $M(\beta)$ is
a very special sort of reducible manifold.

In order to prove our main results we must narrow the profile of a
reducible filling slope which is not a singular slope.

The following result will be proved in the next section of the paper.
\begin{prop} \label{strictsinggen}
  Suppose that $M(\beta)=L(p, 1)\#L(q,1)$ and there are at least two
  isotopy classes of essential surfaces in $M$ of slope $\beta$. Then
  $\beta$ is a singular slope for some closed essential surface in $M$.
\end{prop}

\begin{cor} \label{strictsing}
  Suppose that $M(\beta)=P^3\#P^3$ and $\beta$ is a strict boundary slope.
  Then $\beta$ is a singular slope for some closed essential surface in
  $M$.
\end{cor}

The proposition below, which follows immediately from Theorem
\ref{2.0.3} and Corollary \ref{strictsing}, summarizes the situation.

\begin{prop}\label{sing summary}
  Suppose that $\beta$ is a boundary slope for $M$.  Then one of the
  following three possibilities occurs:
\begin{enumerate}
\item[(1)] $\beta$ is a singular slope for $M$; or
\item[(2)] $M(\beta)$ is homeomorphic to $L_p\#L_q$, where $q>2$; or
\item[(3)] $M(\beta)$ is homeomorphic to $S^2\times S^1$ or $P^3 \#
  P^3$, and $\beta$ is not a strict boundary slope.
\end{enumerate}
\end{prop}
\qed

We end this section by giving the proofs of Theorems \ref{planar} and
\ref{redexc}, assuming Theorems \ref{redfin}, \ref{redvsmall} and
\ref{redtor}.

\begin{proof}[Proof of Theorem \ref{planar}]
  Since we have assumed that $\beta$ is a strict boundary slope, if
  $M(\beta)$ is reducible, then Theorems \ref{redfin}, \ref{redvsmall}
  and \ref{redtor} imply that the corollary holds. On the other hand,
  if $M(\beta)$ is irreducible, then $b_1(M) = 1$ as $M(\beta)$ is
  non-Haken.  Since $\beta$ is a boundary slope, Theorem \ref{2.0.3}
  implies that $\beta$ is a singular slope for a closed essential
  surface in $M$.  Proposition \ref{bgz} now shows that the conclusion
  holds.
\end{proof}

\begin{proof}[Proof of Theorem \ref{redexc}]
  First suppose that $M(\beta)$ is either $S^1 \times S^2$ or $P^3 \#
  P^3$. Since $\beta$ is a strict boundary slope, it follows from
  Proposition \ref{sing summary} that it must be a singular slope for
  $M$. Thus Proposition \ref{bgz} shows that the desired conclusion
  holds.
  
  Next suppose that $M(\beta) \ne S^1 \times S^2, P^3 \# P^3$. Theorem
  0.6 of [BZ3] implies that if $\Delta(\alpha, \beta) > 5$, then
  $M(\alpha)$ is virtually Haken.  In particular, $M(\alpha)$ admits a
  geometric decomposition ([CJ], [Ga2], [Ga3], [GMT]). According to
  [GLu] and either [Wu2] or [Oh], $M(\alpha)$ is irreducible and
  geometrically atoroidal as long as $\Delta(\alpha, \beta) > 3$.
  Further, Theorem \ref{redtor} shows that $M(\alpha)$ is not Seifert
  fibred as long as $\Delta(\alpha, \beta) > 4$. Thus $M(\alpha)$ is
  hyperbolic if $\Delta(\alpha, \beta) > 5$. This proves the first
  claim of the theorem. The second follows similarly since $M(\alpha)$
  admits a geometric decomposition for any slope $\alpha$ if the
  geometrization conjecture holds.
\end{proof}

\section{Singular slopes for $L(p,1)\#L(q,1)$}\label{special singular}

This section contains the proof of Proposition \ref{strictsinggen}.

Let ${\mathcal S}(M)$ denote the set of essential surfaces in $M$.
For each slope $\beta$ on $\partial M$ we set
$${\mathcal S}_\eta(M) = \{ F \in {\mathcal S}(M) : \partial F
\not=\emptyset \mbox{ and $\beta$ is the boundary slope of $F$}\}.$$
For each surface $F\in {\mathcal S}_\beta(M)$ we use $\hat F$ to
denote the closed surface in $M(\beta)$ obtained by attaching meridian
disks to $F$.

We begin with two propositions that give conditions on ${\mathcal
  S}_\beta(M)$ which guarantee that $\beta$ is a singular slope for
some closed essential surface in $M$.  The first is a consequence of
the proof of Theorem \ref{2.0.3} (cf. chapter 2 of [CGLS]).

\begin{prop}  \label{cgls-sing}    {\rm ([CGLS])}
  Suppose that $M(\beta) \cong L_p \# L_q$ and that $F \in {\mathcal
    S}_\beta(M)$ satisfies $|\partial F| \leq |\partial F'|$ for each
  $F' \in {\mathcal S}_\beta(M)$.  If $\hat F$ is not an essential
  $2$-sphere in $M(\beta)$ then $\beta$ is a singular slope for a
  closed, essential surface in $M$.
\end{prop}
\qed

\begin{prop}  \label{bcsz-sing}
  Suppose that $M(\beta) \cong L_p \# L_q$ and let $F \in {\mathcal
    S}_\beta(M)$. If there exists a closed, essential surface $S$ in
  $M$ which is disjoint from $F$, then $\beta$ is a singular slope for
  $S$.
\end{prop}

\begin{proof}
  Since $S$ is closed, essential, and disjoint from $F$, $F$ is not a
  semi-fibre in $M$. On the other hand, $S$ compresses in $M(\beta)
  \cong L_p \# L_q$, so $\beta\in {\mathcal C}(S)$.  Corollary 6.2.3
  of [BCSZ] then shows that $S$ is incompressible in $M(\gamma)$ for
  each slope $\gamma$ on $\partial M$ such that $\Delta(\gamma, \beta)
  \gg 0$. Wu's theorem [Wu1] states that either $\Delta(\gamma,
  \gamma') \leq 1$ for each $\gamma, \gamma' \in {\mathcal C}(S)$, or
  there is a slope $\gamma_0 \in {\mathcal C}(S)$ such that ${\mathcal
    C}(S) = \{\gamma : \Delta(\gamma, \gamma_0) \leq 1\}$.  In the
  first case it is immediate that $\beta$ is a singular slope for $S$.
  In the second case, observe that we must have $\gamma_0 = \beta$
  since otherwise there would exist compressible slopes $\gamma$ with
  $\Delta(\gamma,\beta)$ arbitrarily large.  Thus $\beta$ is a
  singular slope for $S$ in either case.
\end{proof}

We now proceed with the proof of Proposition \ref{strictsinggen},
which depends on the three lemmas below.  We introduce some notational
conventions that will be used in the lemmas.

\begin{con}\label{strictsingsetup}
  Suppose that $M(\beta) \cong L_p \# L_q$ and that $\beta$ is not a
  singular slope for a closed essential surface.  It is evident that
  $b_1(M) = 1$ and, since $\beta$ is not the slope of the rational
  longitude of $M$, that each surface $F \in {\mathcal S}_\beta(M)$ is
  separating.  Fix a surface $P \in {\mathcal S}_\beta(M)$ such that
  $$|\partial P| \leq |\partial F| \mbox{ for each } F \in {\mathcal
    S}_\beta(M).$$
  Since $P$ is connected and separating we have that
  $|\partial P|$ is even, and we set $|\partial P| = n$.  It follows
  from Proposition \ref{cgls-sing} that $\hat P$ is an essential $2$-sphere
  which bounds two punctured lens spaces $\hat X$ and $\hat X'$ in
  $M(\beta)$.  We shall make the convention that $\hat X$ is a
  punctured $L_p$ and $\hat X'$ is a punctured $L_q$.  We let $X$ and
  $X'$ denote the submanifolds bounded by $P$ in $M$, where $X\subset
  \hat X$ and $X'\subset \hat X'$.
\end{con}

In the situation of \ref{strictsingsetup}, we shall say that $(X,P)$
is {\it unknotted} if there is a solid torus $V \subset X$ and an
$n$-punctured disk $D_n$ with outer boundary $\partial_o D_n$ such
that
$$X = V \cup_A (D_n \times I)$$ 
where $A = (\partial_o D_n) \times I$ is identified with an essential
annulus in $\partial V$.

Note that if $(X,P)$ is unknotted and $p = 2$, then $(V, A)$ is a
twisted $I$-bundle pair over a M\"obius band and the induced
$I$-fibring of $A$ coincides with that from $D_n \times I$. Thus $(X,
P)$ is a twisted $I$-bundle.

\begin{lemma} \label{isotopic}
  Assume that $M(\beta) \cong L_p \# L_q$ and that $\beta$ is not a
  singular slope for a closed essential surface.  Let $P\in{\mathcal
    S}_\beta(M)$ be chosen to have the minimal number of boundary
  components.  Suppose that $(X,P)$ is unknotted. If $F \in {\mathcal
    S}_\beta(M)$ is contained in $X$ then $F$ is isotopic to $P$.
\end{lemma}

\begin{proof}
  Write $X = V \cup_A (D_n \times I)$ as above and isotope $F$ so as
  to minimize $|A \cap F|$. Then $F$ intersects $V$ and $D_n \times I$
  in incompressible surfaces. If $A \cap F = \emptyset$, then $F
  \subset D_n \times I$, and therefore Proposition 3.1 of [Wa] implies
  that $F$ is parallel into $D_n \times \{0\} \subset P$. But then
  $|\partial F| \leq n = \frac{1}{2}|\partial P|$, which contradicts
  our choice of $P$. Thus $F \cap A$ consists of a non-empty family of
  core curves of $A$. Another application of [Proposition 3.1, Wa]
  implies that up to isotopy, each component of $F \cap (D_n \times
  I)$ is of the form $D_n \times \{t\}$ for some $t \in (0,1)$. Since
  $|A \cap F|$ has been minimized, it also follows that each component
  of $F \cap V$ is parallel into $\overline{\partial V \setminus A}$.
  It is now simple to see that $F$ is of the form $D_n \times \{t_1\}
  \cup B \cup D_n \times \{t_2\}$ where $0 , t_1 < t_2 < 1$ and $B
  \subset V$ is an annulus as described in the previous sentence. It
  follows that $F$ is isotopic to $P$.  \end{proof}

\begin{lemma}  \label{unknotted} 
  Suppose that $M(\beta) \cong L_p \# L_q$ and that $\beta$ is not a
  singular slope for a closed essential surface in $M$.  Let
  $P\in{\mathcal S}_\beta(M)$ be chosen to have the minimal number of
  boundary components.
  \begin{enumerate}
  \item[(1)] If $M(\beta) \cong L_p \# L_q$ where $L_p \cong \pm L(p,
    1)$, then $(X,P)$ is unknotted.
  \item[(2)] If $M(\beta) \cong L_p \# L_q$ where $L_p \cong \pm L(p,
    1)$ and $L_q \cong \pm L(q, 1)$, then each planar surface in
    ${\mathcal S}_\beta(M)$ is isotopic to $P$.
\end{enumerate}
\end{lemma} 

\begin{proof}
  (1) Suppose that $M(\beta) \cong L_p \# L_q$ where $L_p \cong \pm
  L(p, 1)$.  We will follow the conventions in \ref{strictsingsetup};
  in particular $\hat X$ is the punctured $L_p$ and $|\p P| = 2n$.
  The desired conclusion follows from a combination of [CGLS] and
  [Wu3].  In order to make the application of these two papers clear,
  we must first set up some notation and recall some definitions.
  
  Since $M$ is hyperbolic, $n\geq 2$. The boundary of $P$ cuts the
  boundary of $M$ into $2n$ annuli $A_1, A_1', A_2, A_2',..., A_n,
  A_n'$, occurring successively around $\p M$, such that $\p
  X=P\cup(\cup_{i=1}^n A_i)$ and $\p X'=P\cup(\cup_{i=1}^nA_i')$. Let
  $V$ be the attached solid torus used in forming $M(\beta)$. Then $V$
  may be considered as a union of $2n$ $2$-handles
  $H_1,H_1',H_2,H_2',...,H_n,H_n'$ with attaching regions
  $A_1,A_1',A_2,A_2',...,A_n,A_n'$ respectively.  Let $\hat X$ be the
  manifold obtained from $X$ by adding the 2-handles $H_1,...,H_n$
  along $A_1,...,A_n$ respectively and similarly let $\hat X'$ be the
  manifold obtained from $X'$ by adding the 2-handles $H_1',...,H_n'$
  along $A_1',...,A_n'$. Then $M(\beta)= M\cup V= \hat X\cup_{\hat
    P}\hat X'$, where $\hat P$ is the 2-sphere obtained from $P$ by
  capping off $\p P$ with meridian disks of $V$.  Let $K$ be the core
  curve of the solid torus $V$.  Then $K$ is the union of $2n$ arcs
  $\alpha_1,\alpha_1',\alpha_2,\alpha_2',...,\alpha_n,\alpha_n'$ such
  that $\alpha_1,\alpha_2, ..., \alpha_n$ are properly embedded in
  $\hat X$ with regular neighborhoods $H_1, H_2,...,H_n$ and
  $\alpha_1',\alpha_2', ..., \alpha_n'$ are properly embedded arcs in
  $\hat X'$ with regular neighborhoods $H_1', H_2',...,H_n'$.
  
  Consider the $n$-string tangle $(\hat X; \alpha_1,...,\alpha_n)$ in
  $\hat X$ with strings $\alpha_1,...,\alpha_n$.  Let $P_i=P\cup A_i$
  and call it the $A_i$-tubing surface of $P$.  The surface $P$ is
  said to be {\it $A_i$-tubing compressible} if $P_i$ is compressible
  in $X$, and is said to be {\it completely $A_i$-tubing compressible}
  if $P_i$ can be compressed in $X$ until it becomes a set of annuli
  parallel to $\cup_{j\ne i}A_j$.  The tangle $(\hat X,
  \alpha_1,...,\alpha_n)$ is called {\it completely tubing
    compressible} if it is completely $A_i$-tubing compressible for
  each of $i=1,...,n$.  Since $M$ does not contain an essential torus,
  the argument of [CGLS, 2.1.2] proves that $(\hat X;
  \alpha_1,...,\alpha_n)$ is completely tubing compressible.  Thus for
  each of $i=1,...,n$, there exist disjoint properly embedded disks
  $E_i^{j}$ in $X$, $j\ne i$, such that $\p E^{j}_i$ meets $A_j$ in a
  single essential arc of $A_j$ and is disjoint from $A_k$ if $k\ne
  i,j$ (see [CGLS, 2.1.2] for details).  This in turn implies that if
  $\Omega$ is a proper subset of $\{H_1,...,H_n\}$ then the manifold
  obtained by attaching 2-handles from $\Omega$ to $X$ is a
  handlebody.  In particular for each of $i=1,...,n$,
  $X\cup(\cup_{j\ne i}H_j)$ is a solid torus.  Thus each $\alpha_i$ is
  a core arc of $\hat X$, i.e. its exterior in $\hat X$ is a solid
  torus.
  
  Recall from [Wu3] that a band in a compact 3-manifold $W$ whose
  boundary is a 2-sphere is an embedded disk $D$ in $W$ such that $\p
  D\cap \p W$ consists of two arcs on $\p D$.  A collection of
  properly embedded arcs in $W$ is said to be parallel in $W$ if there
  is a band $D$ in $W$ which contains all these arcs.  It is proved in
  [Wu3] that if $W$ is homeomorphic to a once punctured lens space
  $L(p,1)$ and $(W; \alpha_1,...,\alpha_n)$ is a completely tubing
  compressible tangle, then the arcs $\alpha_1,...,\alpha_n$ are
  parallel in $W$.  Though this result is not explicitly stated in
  [Wu3], its proof is explicitly dealt with in the proof of Theorem 1
  of that paper. Hence in our current situation,
  $\alpha_1,...,\alpha_n$ are parallel arcs in $\hat X$.  Let $D$ be a
  band in $\hat X$ which contains all the arcs and $H$ a regular
  neighborhood of $D$ in $\hat X$.  We may assume that $H$ contains
  every $H_i$.  Since each $\alpha_i$ is a core arc of $H$, $V=\hat
  X\setminus int(H)$ is a solid torus.  More precisely $H$ can be
  considered as a 2-handle and $\hat X$, a once punctured $L_p$, is
  obtained by attaching $H$ to the solid torus $V$ along an annulus
  $A$ in $\p V$. Thus (1) holds.
  
  (2) Now suppose that $M(\beta) \cong L_p \# L_q$ where $L_p \cong
  \pm L(p, 1)$ and $L_q \cong \pm L(q, 1)$. Part (1) of this lemma
  implies that both $(X, P)$ and $(X', P)$ are unknotted. Fix a planar
  surface $F \in {\mathcal S}_\beta(M)$ whose boundary is disjoint
  from $\partial P$, and which has been isotoped to be transverse to
  $P$ and so that $|F \cap P|$ has been minimized. Let ${\mathcal F}$
  be the set of surfaces in ${\mathcal S}_\beta(M)$ isotopic to $F$
  and which satisfy the conditions of this paragraph.
  
  If $F \cap P = \emptyset$, then Lemma \ref{isotopic} implies the
  desired result. Assume then that $F \cap P \ne \emptyset$ and
  consider a component $C$ of $F \cap P$ which is innermost in the
  $2$-sphere $\hat F$. Let $F_0$ be a subset of $F $ whose boundary is
  the union of $C$ and $k$, say, components of $\partial F$. We assume
  that $F $ and $F_0$ are chosen from all the surfaces in ${\mathcal
    F}$ so that $k$ is minimized.  Note that $k > 0$ by the minmality
  of $|F \cap P|$.
  
  Without loss of generality we take $F_0 \subset X = V \cup_A (D_n
  \times I)$ where $A \subset \partial V$ wraps $p$ times around $V$,
  and after an isotopy of $F$ which preserves $P$, we may arrange for
  $F_0$ to be transverse to $A$ and $|F_0 \cap A|$ to be minimal. The
  components of $F_0 \cap A$ are either core circles of $A$ or arcs
  properly embedded in $A$.
  
  First assume that $C \cap A = \emptyset$. Then $F_0 \cap A$ consists
  of core circles of $A$ and an argument like that used in the proof
  of Lemma \ref{isotopic} implies that $F_0$ is parallel into $P$,
  contrary to the minimality of $|F \cap P|$. Thus $C \cap A \ne
  \emptyset$. It follows that $F_0 \cap A$ contains arc components.
  Choose such an arc $\alpha$ which is outermost in the disk $\hat
  F_0$ and let $D_0$ be a planar subsurface of $F_0$ it subtends and
  whose interior is disjoint from $A$. Set $\alpha' =
  \overline{\partial D_0 \setminus \alpha}$.
  
  If $D_0 \subset V$, then $D_0$ is a disk. If $\alpha'$ is an
  essential arc in the annulus $E = \overline{\partial V \setminus
    A}$, then it connects $D_n \times \{0\}$ to $D_n \times \{1\}$.
  Hence $D_0$ is a meridian disk of $V$ and $\partial D_0$ is a dual
  curve on $\partial V$ to the core of $A$. But this is impossible as
  $A$ wraps $p > 1$ times around $V$. Thus $\alpha'$ is an inessential
  arc in $E$. It follows that $\alpha$ is inessential in $A$ and it is
  easy to see that $\alpha$ can be eliminated from $F_0 \cap A$ by an
  isotopy of $X$, contrary to the minimality of $|F_0 \cap A|$.
  
  Suppose next that $D_0 \subset D_n \times I$ so that $\alpha'
  \subset D_n \times \partial I$, say $\alpha' \subset D_n \times
  \{0\}$. Note then that $\alpha$ is inessential in $A$. An argument
  like that used in the previous paragraph shows that $D_0$ cannot be
  a disk. Thus $D_0 \cap \partial M \ne \emptyset$. By Proposition 3.1
  of [Wa], $D_0$ is parallel into $D_n \times \{0\}$, and it is now
  easy to see that $F$ can be isotoped in $M$ to reduce $k$, contrary
  to our choices. This contradiction completes the proof.  \end{proof}

\begin{proof}[Proof of Proposition \ref{strictsinggen}]
  Let ${\mathcal S}_\beta^0(M)\subset {\mathcal S}_\beta(M)$ consist
  of the surfaces in ${\mathcal S}_\beta(M)$ which are isotopic to
  $P$, and set ${\mathcal S}_\beta^1(M) = {\mathcal S}_\beta(M)
  \setminus {\mathcal S}_\beta^0(M)$. By hypothesis, ${\mathcal
    S}_\beta^1(M) \ne \emptyset$. Choose $F \in {\mathcal
    S}_\beta^1(M)$ so that $|\partial F| \leq |\partial F'|$ for all
  $F' \in {\mathcal S}_\beta^1(M)$ and let $Y, Y'$ be the components
  of $M$ split along $F$. Part (2) of Lemma \ref{unknotted} shows that
  $F$ is not planar.
  
  Let $B$ be a component of $Y \cap \partial M$ and consider $F_0 = F
  \cup B$. Let $C_1, C_2, \ldots , C_k$ be the components of the inner
  boundary $F_0^-$ of the maximal compression body of $F_0$ in $Y$. If
  any of the $C_i$ are closed, Proposition \ref{bcsz-sing} shows that
  $\beta$ is a singular slope for a closed essential surface in $M$.
  Suppose then that no $C_i$ is closed. If some $C_i$ is essential,
  the fact that $|\partial C_i| < |\partial F|$ implies that $C_i \in
  {\mathcal S}_\beta^0(M)$ and therefore is isotopic to $P$. Since $F$
  is disjoint from $C_i$, we can isotope $F$ into the complement of
  $P$.  But this is impossible as Lemma \ref{isotopic} would then
  imply that $F \in {\mathcal S}_\beta^0(M)$. Thus each $C_i$ is a
  $\partial$-parallel annulus. Similar arguments show that either
  $\beta$ is a singular slope for a closed essential surface in $M$ or
  for each component $B'$ of $\partial M \cap Y'$, the inner boundary
  of the maximal compression body of $F \cup B'$ in $Y'$ is a family
  of $\partial$-parallel annuli. Hence if $\beta$ is not a singular
  slope for a closed essential surface in $M$, the arguments of \S 2.2
  of [CGLS] imply that $\hat F$ is essential in $M(\beta) \cong L_p \#
  L_q$. This cannot occur since the genus of $F$ is positive. Thus
  $\beta$ is a singular slope for a closed essential surface in $M$.
\end{proof}

\section{Preliminary reductions}\label{preliminary}

In this section we state four propositions which, together with known
results, respectively imply our main theorems: \ref{redfin} --
\ref{redtor}.  Recall that $M$ always denotes a compact, connected,
orientable, irreducible, atoroidal $3$-manifold whose boundary is a
torus.

If $M(\beta)$ is a reducible manifold then it follows from [GLu] that
$\Delta(\alpha,\beta)\leq 1$ for any slope $\alpha$ such that
$M(\alpha)$ is reducible.  If $b_1(M) \geq 2$ then it follows from
[BGZ, Proposition 5.1] that $\Delta(\alpha,\beta)\leq 1$ for any slope
$\alpha$ such that $M(\alpha)$ is not hyperbolic.  The conclusions of
all three of the main theorems hold when $\Delta(\alpha,\beta)\leq 1$.
Thus in the proofs of these theorems we may assume without loss of
generality that $M(\alpha)$ is irreducible and $b_1(M) = 1$.

Next we recall that, since $\beta$ is a boundary slope, it follows
from Proposition \ref{sing summary} that one of the following three
possibilities occurs:
\begin{enumerate}
\item[(1)] $\beta$ is a singular slope for a closed essential surface
  in $M$; or
\item[(2)]  $M(\beta)$ is homeomorphic to $L_p\#L_q$, where $q>2$; or
\item[(3)] $M(\beta)=S^2\times S^1$ or $P^3 \# P^3$ and $\beta$ is not
  a strict boundary slope.
\end{enumerate}
Since the conclusion of Proposition \ref{bgz} implies that of each of the
three theorems, we may also assume that neither $\alpha$ nor $\beta$ is
a singular slope for any closed essential surface in $M$.

Therefore Theorems \ref{redfin} -- \ref{redtor}
follow respectively from the following four propositions.

\begin{prop} \label{TLF}
  Suppose that $b_1(M) = 1$ and neither $\alpha$ nor $\beta$ is a
  singular slope for a closed, essential surface in $M$. Assume as
  well that $M(\beta)$ is either a connected sum of two non-trivial
  lens spaces or $S^1 \times S^2$.  If $M(\alpha)$ has finite
  fundamental group, then $\Delta(\alpha,\beta)\leq 2$. Furthermore,
  if $\Delta(\alpha,\beta) = 2$ then $H_1(M) \cong \mathbb Z \oplus
  \mathbb Z/2$, $M(\beta) \cong L_2 \# L_3$ and $\pi_1(M(\alpha))
  \cong O_{24}^* \times \mathbb Z/j$ where $O_{24}^*$ is the binary
  octahedral group.
\end{prop}

\begin{prop} \label{weakTLVS}
  Suppose that $b_1(M) = 1$ and neither $\alpha$ nor $\beta$ is a
  singular slope for a closed, essential surface in $M$. Assume as
  well that $M(\alpha)$ is irreducible and $M(\beta)$ is either a
  connected sum of two non-trivial lens spaces or $S^1 \times S^2$.  If
  $M(\alpha)$ is a very small manifold and $\beta$ is a strict
  boundary slope, then $\Delta(\alpha,\beta)\leq 3$.
\end{prop}

\begin{prop} \label{TLVS} Suppose that $b_1(M) = 1$ and neither
  $\alpha$ nor $\beta$ is a singular slope for a closed, essential
  surface in $M$. Assume as well that $M(\alpha)$ is irreducible and
  $M(\beta)$ is either a connected sum of two non-trivial lens spaces
  or $S^1 \times S^2$. If $M(\alpha)$ is a very small manifold which
  admits a geometric decomposition, then $\Delta(\alpha, \beta) \leq
  2$.
\end{prop}

\begin{prop} \label{TLS}
  Suppose that $b_1(M) = 1$ and neither $\alpha$ nor $\beta$ is a
  singular slope for a closed, essential surface in $M$. Assume as
  well that $M(\beta)$ is a connected sum of two non-trivial lens
  spaces. If $\beta$ is a strict boundary slope and $M(\alpha)$ a
  $\pi_1$-injective immersion of a torus, then
  $\Delta(\alpha,\beta)\leq 4$. Moreover, if $\Delta(\alpha, \beta) =
  4$, then $M(\alpha)$ is a Seifert fibred space with base orbifold
  ${\mathcal B}$ of $M(\alpha)$ of the form $S^2(r,s,t)$ where
  $(r,s,t)$ is a hyperbolic triple and $4$ divides at least one of
  $r,s,t$.
\end{prop}

These four propositions will be proved in Sections \ref{PTLF},
\ref{weakPTLVS}, \ref{PTLVS} and \ref{PTLS} respectively.

\section{Background Results on $PSL_2(\mathbb C)$-Character Varieties}
\label{background}

In this section we gather together some background material on
$PSL_2(\mathbb C)$-character varieties that will be used in the
proofs of our main results. See [CS], [CGLS], [BZ1], [BZ2], and [BZ5]
for more details. As above, $M$ will denote a compact, connected,
orientable, hyperbolic $3$-manifold with boundary a torus.

\begin{defs}\label{defs}
  Let $\pi$ be a finitely generated group.  We shall denote by
  $R_{PSL_2}(\pi)$ and $X_{PSL_2}(\pi)$ respectively the
  $PSL_2(\mathbb C)$-representation variety and the $PSL_2(\mathbb
  C)$-character variety of $\pi$ .  (Note that these are affine
  algebraic sets, but are not necessarily irreducible.) The map $t:
  R_{PSL_2}(\pi) \to X_{PSL_2}(\pi)$ which sends a representation
  $\rho$ to its character $\chi_\rho$ is a regular map.  When
  $\pi$ is the fundamental group of a path-connected space $Y$ we will
  frequently denote $R_{PSL_2}(\pi)$ by $R_{PSL_2}(Y)$ and
  $X_{PSL_2}(\pi)$ by $X_{PSL_2}(Y)$.
  
  There is a unique conjugacy class of homomorphisms $\eta:
  H_1(\partial M) \to \pi_1(M)$, obtained by composing the inverse of
  the Hurewicz isomorphism $\pi_1(\partial M) \to H_1(\partial M)$
  with some homomorphism $\pi_1(\partial M) \to \pi_1(M)$ induced by
  inclusion.  To simplify notation, we shall often suppress $\eta$ in
  statements that are invariant under conjugation in $PSL_2(\mathbb
  C)$.  For instance, given $\rho \in R_{PSL_2}(\pi)$ and $\alpha \in
  H_1(\partial M)$ we may write $\rho(\alpha) = \pm I$ to indicate
  that $\eta(\alpha)$ is contained in the kernel of $\rho$ for every
  choice of $\eta$.
  
  By a {\it curve} in an an affine algebraic set we will mean an
  irreducible algebraic subset of dimension 1.  Suppose that $X_0$ is
  a curve in $X_{PSL_2}(M)$ and let $\widetilde X_0$ denote the smooth
  projective model of $X_0$.  There is a canonically defined
  quasi-projective curve $X_0^\nu \subset\widetilde X_0$ which
  consists of all points of $\widetilde X_0$ that correspond to points
  of $X_0$.  In particular there is a regular, surjective, birational
  isomorphism $\nu: X_0^\nu \to X_0$.  The points of $X_0^\nu$ are
  called {\it ordinary} points and the points in the finite set
  $\widetilde X_0 - X_0^\nu$ are called {\it ideal} points.
  
  It follows from [Lemma 4.1, BZ2] that for every curve $X_0$ in
  $X_{PSL_2}(M)$ there exists an algebraic component $R(X_0)$ of
  $R_{PSL_2}(M)$ such that $t(R(X_0)) = X_0$.
\end{defs}

To each homology class $a \in H_1(\partial M)$ we can associate a
regular function $f_a: X_0 \to \mathbb C$ given by $f_a(\chi) =
\chi(a)^2 -4$.  Each $f_a$ lifts to a rational function, also denoted
by $f_a$, on $\widetilde X_0$.  It is shown in [CGLS] (see also [BZ2])
that the degrees of these functions on $\widetilde X_0$ vary in a
coherent fashion. Indeed, there is a seminorm $\|\cdot\|_{X_0}:
H_1(\partial M; \mathbb R) \to [0, \infty)$, called the {\it
  Culler-Shalen seminorm} of $X_0$, determined by the condition that
for each $a \in H_1(\partial M)$, $\|a\|_{X_0}$ is the degree of $f_a$
on $\widetilde X_0$. As in [CGLS] we use $Z_x(f)$ to denote the order
of zero of a rational function $f$ on $\widetilde X_0$ at a point
$x\in \widetilde X_0$, and use $\Pi_x (f)$ to denote the order of pole
of $f$ at a point $x\in \widetilde X_0$.  Then
\begin{equation}\label{normformula}
\|a\|_{X_0} = \sum_{x \in \widetilde X_0} Z_x(f_a) = \sum_{x \in
    \widetilde X_0} \Pi_x(f_a).
\end{equation}
If $\|\cdot\|_{X_0} \ne 0$, we define $$s_{X_0} = \mbox{
  min}\{\|a\|_{X_0} \; | \; a \in H_1(\partial M), \|a\|_{X_0} \ne 0\}
\in \mathbb Z \setminus \{0\}.$$
  
We note that $f_a = f_{-a}$.  As a notational convenience, if $\alpha
= \{\pm a\}$ is a slope on $\partial M$ then we shall set $f_\alpha
\dot= f_a = f_{-a}$, and define $\|\alpha\|_{X_0} \dot= \|a\|_{X_0} =
\|-a\|_{X_0}$.

It is possible that $\|\cdot\|_{X_0} \ne 0$ but $\|\beta\|_{X_0} = 0$
for some slope $\beta$ on $\partial M$.  In this case the slope
$\beta$ is the unique slope on $\partial M$ of norm $0$, and we shall
call $X_0$ a {\it $\beta$-curve}. If $X_0$ is a $\beta$-curve then for
any slope $\alpha$ on $\partial M$ we have
\begin{equation}\label{minformula}
\|\alpha\|_{X_0} = \Delta(\alpha, \beta) s_{X_0} .
\end{equation}
Hence if $\beta^*$ is a {\it dual slope} for $\beta$,
that is, a slope such that $\Delta(\beta, \beta^*) = 1$, then 
$$s_{X_0} = \|\beta^*\|_{X_0}.$$

If $\beta$ is any slope on $\partial M$ then we may regard
the character variety $X_{PSL_2}(M(\beta))$ as an algebraic subset of
$X_{PSL_2}(M)$. To see this, note that $R_{PSL_2}(M(\beta))$ can be identified with the Zariski closed, conjugation invariant subset $R_\beta(M) :=\{\rho \in R_{PSL_2}(M) : \rho(\beta) = \pm I\}$ of $R_{PSL_2}(M)$. Theorem 3.3.5(iv) of [Ne] shows that the image of $R_\beta(M)$ in $X_{PSL_2}(M)$ is Zariski closed and can be identified with $X_{PSL_2}(M(\beta))$. We
note that if $X_0$ is a curve in $X_{PSL_2}(M(\beta)) \subset X_{PSL_2}(M)$
such that $\|\cdot\|_{X_0} \ne 0$, then $X_0$ is a $\beta$-curve.

The following proposition is proved in [B].

\begin{prop} \label{betacurve} {\rm ([Proposition 6.2, B])}
  Let $X_0 \subset X_{PSL_2}(M(\beta))$ be a $\beta$-curve for a slope
  $\beta$ on $\partial M$.  Let $\beta^*$ be a dual slope for $\beta$
  and let $\alpha\not=\beta$ be a slope on $\partial M$.  Then
\begin{enumerate}
\item[(1)] For any point $x \in X_0^\n$ and any representation $\rho$
  such that $\chi_\rho = \nu(x)$ we have
\begin{enumerate}
\item[(a)] If $Z_x(f_\alpha) > 0$ then $\rho(\pi_1(\partial M))$ is either
  parabolic, or a finite cyclic group whose order divides
  $\Delta(\alpha, \beta)$; and
\item[(b)] $Z_x(f_\alpha) \geq Z_x(f_{\beta^*})$, with equality if and
  only if $\rho(\pi_1(\partial M))$ is parabolic or trivial.
\end{enumerate}
\item[(2)] If $f_{\beta^*}$ has a pole at each ideal point of
  $\widetilde X_0$ then for every divisor $d > 1$ of $\Delta(\alpha,
  \beta)$ there exists $x \in X_0^\nu$ such that $Z_x(f_\alpha) >
  Z_x(f_{\beta^*})$, and $\rho(\pi_1(\partial M))$ is a cyclic group of order $d$
for every representation $\rho$ such that $\chi_\rho = \nu(x)$ 
\end{enumerate}
\end{prop}
\qed

We call a subvariety $X_0$ of $X_{PSL_2}(M)$ {\it non-trivial} if it
contains the character of an irreducible representation.

For some applications we need a stronger condition on $X_0$ than
non-triviality. A character $\chi_\rho \in X_0$ is called {\it
  virtually reducible} if there is a finite index subgroup $\tilde
\pi$ of $\pi_1(M)$ such that $\rho|\tilde \pi$ is reducible. We will
say that $X_0$ is {\it virtually trivial} if every point of $X_0$ is a
virtually reducible character.  The proof of Proposition 4.2 of [BZ5]
shows that if a subvariety $X_0$ of $X_{PSL_2}(M)$ is non-trivial, but
contains infinitely-many virtually reducible characters, then $X_0$ is
virtually trivial and $X_0$ is a curve of characters of
representations $\pi_1(M) \to {\mathcal N} \subset PSL_2(\mathbb C)$
where
$${\mathcal N} = \{\pm \left( \begin{array}{cc} z & 0 \\ 0 & z^{-1} \end{array}
\right), \pm \left( \begin{array}{cc} 0 & w \\ -w^{-1} & 0 \end{array}
\right) \; | \; z, w \in {\mathbb C}^*\} \subset PSL_2(\mathbb C).$$

\subsection*{Ideal points, essential surfaces, and singular slopes} 

One of the key relations between $3$-manifold topology and
$PSL_2(\mathbb C)$-character varieties is the construction
described in [CS] which associates essential surfaces in a
$3$-manifold $M$ to ideal points of curves in $X_{PSL_2}(M)$.

\begin{prop} {\rm{([CS], [CGLS, \S 1.3], [BZ2])}} \label{idealsurf} 
Let $X_0$ be a non-trivial curve in $X_{PSL_2}(M)$ and $x$ an ideal
point of $X_0$.  One of the following mutually exclusive
alternatives holds: Either
\begin{enumerate}
\item[(1)] there is a unique slope $\alpha$ on $\partial M$ such
that $f_\alpha$ is finite-valued at $x$; or
\item[(2)] $f_\alpha$ is finite-valued for every slope $\alpha$ on
$\partial M$.
\end{enumerate}
In case $(1)$ the slope $\alpha$ is a boundary slope.  Moreover, if $X_0$ is
not virtually trivial then $\alpha$ must be a strict boundary slope.  In
case $(2)$ $M$ contains a closed, essential surface.
\end{prop}
\qed

If, as in case (1) of the proposition, there is a unique slope
$\alpha$ on $\partial M$ such that $f_\alpha(x) \in \mathbb C$, we say
that the boundary slope $\alpha$ is {\it associated to $x$}.

\begin{prop} {\rm (Propositions 4.10 and 4.12 of [BZ2])} \label{4.10-4.12}
  Suppose that $x$ is an ideal point of a non-trivial curve $X_0$ in
  $X_{PSL_2}(M)$ and that $\beta$ is a slope on $\partial M$ such that
  every closed, essential surface in $M$ associated to $x$ is
  compressible in $M(\beta)$.  Suppose further that $f_\delta$ is
  finite-valued at $x$ for every slope $\delta$ on $\partial M$. If
  either
\begin{itemize}
\item $X_0 \subseteq X_{PSL_2}(M(\beta))$, or
\item $Z_x(f_\beta) > Z_x(f_\delta)$ for some slope $\delta$ on
  $\partial M$
\end{itemize}
then $\beta$ is a singular slope for some closed essential surface in
$M$.
\end{prop}
 \qed

\subsection*{The $PSL_2$ character variety of $L_p \# L_q$}

It was shown in Example 3.2 of [BZ2] that $X_{PSL_2}(\mathbb Z/p *
\mathbb Z/q)$ is a disjoint union of a finite number of isolated
points and $[\frac{p}{2}][\frac{q}{2}]$ non-trivial curves, each
isomorphic to a complex line. If we fix generators $x$ and $y$ of the
two cyclic free factors of $\mathbb Z/p * \mathbb Z/q$, then each
curve consists of characters of representations which send $x$ and $y$
to elliptic elements of orders dividing $p$ and $q$ respectively.
Such a curve is parametrized by the complex distance between the axes
of these two elliptic elements.

Explicit parametrizations of the curves in $X_{PSL_2}(\mathbb Z/p *
\mathbb Z/q)$ can be given as follows.
For integers $j, k$ with $1 \leq j \leq
[\frac p2]$ and $1 \leq k \leq [\frac q2]$, set $$\lambda = e^{\pi
  ij/p}, \mu = e^{\pi ik/q},\tau = \mu + \mu^{-1}.$$
For $z \in
\mathbb C$ define $\rho_z \in R_{PSL_2}(\mathbb Z/p * \mathbb Z/q)$ by
$$\rho_z(x)=\pm \left(\begin{array}{cc} \lambda & 0\\0&
\lambda^{-1}
\end{array}\right), \bar \rho_z(y)=\pm \left(\begin{array}{cc}
z & 1\\z(\tau -z) -1 & \tau -z
\end{array}\right) .$$
The characters of the representations $\rho_z$ parameterize a
curve $X(j,k) \subset X_{PSL_2}(\mathbb Z/p * \mathbb Z/q)$. Moreover 
the correspondence
$\mathbb C \to X(j,k)$, $z \mapsto \chi_{\rho_z}$, is bijective if $j < [\frac p2]$ and  $k < [\frac
q2]$ and a 2-1 branched cover otherwise.

We shall denote by $D_k$ the dihedral group of order $2k$.
Recall that a finite subgroup of $PSL_2(\mathbb C)$ is either cyclic
or dihedral, or else it is isomorphic to the tetrahedral goup
$T_{12}$, the octahedral group $O_{24}$, or the icosahedral group
$I_{60}$.

The following elementary, but tedious, lemma characterizes
the points in the curve $X(j,k)$ which correspond to the character of a
representation with finite image.  We leave its verification to the
reader.

\begin{lemma}  \label{charcount}
  Fix integers $2 \leq p \leq q$.  Let $X_{p,q}$ be the union of all
  curves $X(j,k) \subset X_{PSL_2}(\mathbb Z/p * \mathbb Z/q)$ such
  that $j$ and $k$ are relatively prime to $p$ and $q$ respectively.
  Then
\begin{enumerate}
\item[(1)] An irreducible component $X(j,k)\subset X_{p,q}$ contains
  exactly two reducible characters if $p > 2$, and one if $p = 2$.
\item[(2)] An irreducible component $X(j,k)\subset X_{p,q}$ contains
  the character of an irreducible representation $\rho$ whose image
  lies in $\mathcal N$ if and only if $p=2$.  Moreover, if $p = 2$ and
  $q > 2$ then there is exactly one such character $\chi_\rho$ and the
  image of $\rho$ is $D_q$.
\item[(3)] $X_{p,q}$ contains the character of a representation whose
  image is $T_{12}$ if and only if $(p,q) \in\{(2,3), (3,3)\}$. If
  $(p,q)=(2,3)$ there is a unique such character and if $(p,q)=(3,3)$
  then there are two.
\item[(4)] $X_{p,q}$ contains the character of a representation whose
  image is $O_{24}$ if and only if $(p,q) \in \{(2,3), (2,4), (3,4),
  (4,4)\}$.  If $(p, q) = (3,4)$ there are two such characters, and in
  the remaining cases there is only one.
\item[(5)] $X_{p,q}$ contains the character of a representation whose
  image is $I_{60}$ if and only if $(p,q) \in \{(2,3), (2,5), (3,3),
  (3,5), (5,5)\}$. There are eight such characters if $(p, q)=(3,5)$
  or $(p,q)=(5,5)$, four if $(p,q)=(2,5)$, and two if $(p,q)=(2,3)$ or
  $(p,q)=(3,3)$.
\end{enumerate}
\end{lemma}
\qed

The next result follows from Proposition \ref{4.10-4.12} and work of
Culler, Shalen and Dunfield. Recall that if $X_0 \subset X_{PSL_2}(M)$
is a $\beta$-curve and $\beta^*$ is a dual class to $\beta$, then
$s_{X_0} = \|\beta^*\|_{X_0} = \sum_{x \in \widetilde X_0}
\Pi_x(f_{\beta^*})$.

\begin{prop} \label{lb}
  Suppose that $M(\beta) \cong L_p \# L_q$ and let $x$ be an ideal
  point of the curve $X(j,k)\subset X_{PSL_2}(M(\beta)) \subset
  X_{PSL_2}(M)$.  Then either
\begin{enumerate}
\item[(1)] $\beta$ is a singular slope for a closed essential surface
  in $M$, or
\item[(2)] $\|\cdot\|_{X(j,k)} \ne 0$ and
$$s_{X(j,k)} \geq  \Pi_x(f_{\beta^*}) \geq
\left\{
\begin{array}{ll}
4 & \mbox{{\rm if $j \ne \frac{p}{2}$ and $k \ne \frac{q}{2}$}} \\
2 & \mbox{{\rm if either $j = \frac{p}{2}$ or $k =
\frac{q}{2}$}.}
\end{array}
\right. $$
\end{enumerate}
\end{prop}

\begin{proof}
  Suppose that $\beta$ is not a singular slope for a closed essential
  surface in $M$. Then Proposition \ref{4.10-4.12} implies that for
  each ideal point $x$ of $X(j,k)$ and for each slope $\alpha \not=
  \beta$, we have $f_\alpha(x) = \infty$.
  
  The natural surjection $\phi: \mathbb Z/2p * \mathbb Z/2q \to
  \mathbb Z/p * \mathbb Z/q$ induces an inclusion $\phi^*:
  X_{PSL_2}(\mathbb Z/p * \mathbb Z/q) \to X_{PSL_2}(\mathbb Z/2p *
  \mathbb Z/2q)$. Given a curve $X_0 \subset X_{PSL_2}(L_p \# L_q) =
  X_{PSL_2}(\mathbb Z/p * \mathbb Z/q)$, there is a curve $Y_0 \subset
  X_{SL_2}(\mathbb Z/2p * \mathbb Z/2q)$ whose image in
  $X_{PSL_2}(\mathbb Z/2p * \mathbb Z/2q)$ coincides with
  $\phi^*(X_0)$. The associated regular map $g: Y_0 \to X_0$ has
  degree $1$ if $j \ne \frac{p}{2}$ and $k \ne \frac{q}{2}$ and is of
  degree $2$ otherwise. Now $Y_0$ is also a complex line and so has a
  unique ideal point $y$. Extend $g$ to a map $\tilde
  g:\widetilde{Y_0}\to \widetilde X_0$ between the smooth projective
  models, and observe that $\tilde g(y) = x$. If $\tilde \beta^* \in
  \phi^{-1}(\beta^*)$ it is easy to see that $f_{\beta^*} \circ \tilde
  g = f_{\tilde \beta^*}$. It can be shown that $$\Pi_x(f_{\beta^*}) =
  \left\{
  \begin{array}{ll}
  \Pi_y(f_{\tilde \beta^*}) & \mbox{{\rm if $j \ne \frac{p}{2}$ and
  $k \ne \frac{q}{2}$}} \\ &\\ \frac12 \Pi_y(f_{\tilde \beta^*}) &
  \mbox{{\rm if either $j = \frac{p}{2}$ or $k = \frac{q}{2}$.}}
  \end{array} \right. $$
  We are reduced then to calculating $\Pi_y(f_{\tilde \beta^*})$.
  
  According to Dunfield (Proposition 2.2 of [Dn]), we may choose the
  simplicial tree $T_y$ associated to $y$ so that $\Pi_y(f_{\tilde
    \beta^*})$ equals the translation length $l(\tilde \beta^*)$ of
  the automorphism of $T_y$ associated to $\tilde \beta^*$. Now the
  action of $\mathbb Z/2p * \mathbb Z/2q$ on $T_y$ factors through an
  action of $\mathbb Z/p * \mathbb Z/q$, which in turn determines an
  action of $\pi_1(M)$ on $T_y$ via the surjection $\pi_1(M) \to
  \pi_1(M(\beta)) = \mathbb Z/p * \mathbb Z/q$. In particular
  $l(\tilde \beta^*) = l(\beta^*)$, where we have identified $\beta^*$
  with its image in $\pi_1(M)$ under one of the homomorphisms in the
  conjugacy class $\eta$ (see \ref{defs}).
  
  Consider now an essential surface $F$ properly embedded in $M$ which
  is dual to the action of $\pi_1(M)$.  The observation above implies
  that $F$ can be chosen so that $|\partial F| = l(\beta^*)$.  Let
  $F_0$ be a component of $F$ with non-empty boundary.  Note that $|\p
  F_0|$ is even since $F_0$ is separating in $M$. If $|\partial F_0| =
  2$, then the genus of $F_0$ is at least $1$ since $M$ is hyperbolic.
  The proof of Theorem 2.0.3 of [CGLS] then shows that $\beta$ is the
  singular slope for some closed essential surface, contrary to our
  hypotheses. Hence $$\Pi_y(f_{\tilde \beta^*}) = l(\tilde \beta^*) =
  l(\beta^*) = |\partial F| \geq 4$$
  and
  $$\Pi_x(f_{\beta^*}) \geq \left\{
\begin{array}{ll}
4 & \mbox{{\rm if $j \ne \frac{p}{2}$ and $k \ne \frac{q}{2}$}} \\
2 & \mbox{{\rm if either $j = \frac{p}{2}$ or $k =
\frac{q}{2}$.}}
\end{array} \right. $$
\end{proof}

\subsection*{Jumps in multiplicities of zeroes} 

Let $X_0$ be a non-trivial curve in $X_{PSL_2}(M)$ Recall that
$R(X_0)$ is the unique $4$-dimensional subvariety of $R_{PSL_2}(M)$
satisfying $t(R(X_0)) = X_0$.  Suppose that $\alpha$ is a slope on
$\partial M$ such that $f_\alpha|X_0 \not= 0$.  As a means to
estimate $\|\alpha\|_{X_0}$, we will be interested in the set
$$J_{X_0}(\alpha) = 
\{x \in \widetilde X_0 \; | \; Z_x(f_{\alpha}) > Z_x(f_\delta)
\mbox{ for some slope } \delta \mbox{ such that } f_\delta\not= 0\}.$$

\begin{lemma} \label{ordjump} 
Suppose that $x \in J_{X_0}$ is not an ideal point.
\begin{enumerate}
\item[(1)] {\rm [CGLS, Proposition 1.5.4]} If $\chi_\rho=\nu(x)$
  then $\rho(\alpha) = \pm I$.
\item[(2)] {\rm [B, Proposition 2.8]} If $b_1(M) = 1$, there
  exists a representation $\rho$, which is either irreducible or has
  non-Abelian image, such that $\chi_\rho=\nu(x)$.
\end{enumerate}
\end{lemma}
\qed

(Note that there exist irreducible $PSL_2(\mathbb C)$ representations
whose image is a Klein $4$-group, and hence is Abelian.)

\begin{lemma} \label{jumpnonideal}
  Let $X_0 \subset X_{PSL_2}(M(\beta)) \subset X_{PSL_2}(M)$ be a
  $\beta$-curve for a slope $\beta$ on $\partial M$.  Let
  $\beta^*=\{\pm b^*\}$ be a dual slope for $\beta$.  Suppose that
  $\alpha$ is a slope on $\partial M$ such that $\Delta(\alpha,\beta)
  > 1$.  For any non-ideal point $x\in J_{X_0}(\alpha)$ and any
  representation $\rho$ such that $\chi_\rho = \nu(x)$ we have that
  $\rho(b^*)$ is an elliptic element with order $d$ for some
  divisor $d > 1$ of $\Delta(\alpha,\beta)$.
\end{lemma}

\begin{proof}
  First observe that for any slope $\delta$ on $\partial M$ we have
  $f_\delta = f_{\Delta(\delta, \beta) \beta^*}$ and so $Z_x(f_\delta)
  = \Delta(\delta,\beta)Z_x(f_{\beta^*})$.  In particular, since $x\in
  J_{X_0}(\alpha)$, we must have $Z_x(f_{\beta^*}) > 0$.  Thus
  $Z_x(f_\alpha) = \Delta(\alpha,\beta)Z_x(f_{\beta^*}) >
  Z_x(f_{\beta^*})$.  It now follows from Proposition \ref{betacurve}
  that $\rho(\pi_1(\partial M))$ is a cyclic group of order $d >1$
  where $d$ divides $\Delta(\alpha,\beta)$.  Since this cyclic group
  is generated by $\rho(b^*)$, the lemma follows.  \end{proof}

\begin{prop} \label{jumpideal}
  Let $X_0 \subset X_{PSL_2}(M)$ be a non-trivial curve and let
  $\alpha$ be a slope on $\partial M$ such that $f_\alpha|X_0 \ne 0$.
  Suppose that there is no closed, essential surface in $M$ which
  remains essential in $M(\alpha)$. If $x \in J_{X_0}(\alpha)$ is an
  ideal point, then either
\begin{enumerate}
\item[(1)] $\alpha$ is a singular slope for a closed, essential
  surface in $M$, or
\item[(2)] for any slope $\beta \ne \alpha$, $f_\beta$ has a pole at
  $x$. In particular $\alpha$ is a boundary slope and $X_0$ is not a
  $\beta$-curve. Moreover if $b_1(M) = 1$ then $M(\alpha)$ is
  either a Haken manifold, $S^1 \times S^2$, or a connected
sum of two non-trivial lens spaces.
\end{enumerate}
\end{prop}

\begin{proof}
  Suppose that $\alpha$ is not a singular slope for a closed,
  essential surface in $M$. It then follows from Proposition
  \ref{4.10-4.12} that for any slope $\beta\not=\alpha$ the function
  $f_\beta$ has a pole at $x$.  Hence Proposition \ref{idealsurf}
  shows that $\alpha$ is a boundary slope. Finally if $b_1(M) = 1$, we
  can apply Theorem \ref{2.0.3} to deduce that $M(\alpha)$ is either
  Haken, $S^1 \times S^2$, or is a connected sum of two non-trivial
  lens spaces.
\end{proof}

\begin{prop} \label{jumpord}
  Let $X_0$ be a non-trivial curve in $X_{PSL_2}(M)$ and $\alpha$ a
  slope on $\partial M$ such that $f_\alpha|X_0 \not= 0$.  Suppose
  that $J_{X_0}(\alpha)$ contains an ordinary point $x$ of $X_0^\nu$
  and that there exists a representation $\rho$, which is either
  irreducible or has non-Abelian image, such that $\chi_\rho=\nu(x)$.
  If either
\begin{enumerate}
\item[(i)] $H^1(M(\alpha); sl_2(\mathbb C)_\rho) = 0$ and
  $\rho(\pi_1(\partial M)) \not = \{\pm I\}$, or
\item[(ii)] there is a slope $\beta$ such that $X_0 \subset X_{PSL_2}(M(\beta))$ and $H^1(M; sl_2(\mathbb C)_\rho) \cong \mathbb C$ (for instance the latter holds when $M(\beta) \cong L_p \# L_q$),  
\end{enumerate}
then
$$
Z_x(f_{\alpha}) = 
\begin{cases}
Z_x(f_\beta) + 1 & \mbox{ if $\rho$ is conjugate into $\mathcal N$;} \\
Z_x(f_\beta) + 2 & \mbox{ otherwise.}
\end{cases}
$$
Moreover, in case {\rm (i)} $\nu(x)$ is a simple point of
$X_{PSL_2}(M)$ and in case {\rm (ii)} $\nu(x)$ is a simple point of
$X_{PSL_2}(M(\beta))$.
\end{prop}

\begin{proof}
  If hypothesis (i) holds the conclusion follows from [Theorem 2.1,
  BB].
  
  Assume that hypothesis (ii) holds. Let $\beta^*$ be a dual slope to
  $\beta$ and fix simple closed curves $a$, $b$ and $b^*$ on $\partial
  M$ such that $\alpha = \{\pm [a]\}$, $\beta = \{\pm [b]\}$ and
  $\beta^*=\{\pm [b]^*\}$.  We also identify $[a]$, $[b]$ and $[b^*]$
  with their images under a homomorphism in the conjugacy class $\eta$
  (see \ref{defs}).
  
  Observe that Proposition \ref{betacurve} implies that
  $\rho(\pi_1(\partial M))$ is a non-trivial, finite cyclic group.
  Thus, $\rho(\pi_1(\partial M))$ is generated by $\rho([b^*])$.
  After possibly replacing $\rho$ by a conjugate representation, we
  may assume that $$\rho([b^*]) = \pm \left(
  \begin{array}{cc} t & 0 \\ 0 & t^{-1} \end{array} \right)$$ where $t
  \ne \pm 1$.
  
  Since $X_0 \subset X_{PSL_2}(M(\beta))$ and $H^1(M; sl_2(\mathbb
  C)_\rho) \cong \mathbb C$, Theorem A of [B] holds in our situation.
  In particular, the Zariski tangent space of $X_0$ at $\chi_\rho$ can
  be identified with $H^1(M(\beta); sl_2(\mathbb C)_\rho) \cong
  \mathbb C$.  We can therefore find a $1$-cocycle $u \in
  Z^1(M(\beta); sl_2(\mathbb C)_\rho)$ such that $\bar u \ne 0 \in
  H^1(M(\beta); sl_2(\mathbb C)_\rho)$ and an analytic curve
  $\chi_{\rho_s}$ in $X_0$ of the form $\rho_s = \exp(su +
  O(s^2))\rho$ defined for $|s|$ small.  Applying the arguments of \S
  1.1.1 and \S 1.2.1 of [BB] to this curve, modified to the
  $PSL_2(\mathbb C)$ setting (cf. \S 2 of [BB]), shows that the
  identities
$$
Z_x(f_{\alpha}) = 
\begin{cases}
Z_x(f_\beta) + 1 & \mbox{ if $\rho$ is conjugate into $\mathcal N$;} \\
Z_x(f_\beta) + 2 & \mbox{ otherwise.}
\end{cases}
$$
  hold as long as we can prove that $u([a]) \ne 0$. 
  
  Suppose that $u([a]) = 0$ in order to arrive at a contradiction.  We
  also have $u([b]) = 0$, since $u \in Z^1(M(\beta); sl_2(\mathbb
  C)_\rho)$, and thus $u(m[a]+n[b]) = 0$ for each pair of integers $m,
  n$.  Let $u([b^*]) = \left( \begin{array}{cc} p & q \\ r & -p
  \end{array} \right).$ We have assumed that $f_\alpha|X_0 \not= 0$,
and therefore $[a]$ and $[b]$ span a subgroup of index $k < \infty$ of
$H_1(\partial M)$. Then
\begin{eqnarray} 0 & = & u([b^*])^k = \sum_{j=0}^{k-1} \rho([b^*])^j u([b^*])
\rho([b^*])^{-j} \nonumber \\ 
& = & 
\left( \begin{array}{cc} kp & (1 + t^2 + \ldots + t^{2(k-1)})q \\ 
(1 + t^{-2} + \ldots + t^{-2(k-1)})r & -kp \end{array}
\right), \nonumber 
\end{eqnarray}
and therefore $p = 0$. Consider the coboundary $\delta^0 :
sl_2(\mathbb C) \to Z^1(M(\beta); sl_2(\mathbb C)_\rho)$ given by
$(\delta^0(A))(w) = A - \rho(w)A\rho(w)^{-1}$ and set
$u_1 = u - \delta^0(\left( \begin{array}{cc} 0 & \frac{q}{1 - t^2} \\
    \frac{r}{1 - t^{-2}} & 0 \end{array} \right))$. Since $\rho([b]) =
\pm I$ we have $u_1([b]) = u([b]) = 0$, while the fact that
$\rho([b^*]) = \pm \left( \begin{array}{cc} t & 0 \\ 0 & t^{-1}
  \end{array} \right)$ implies that $u_1([b^*]) = 0$ also. Hence $u_1
= 0$, which is impossible as $0 \ne \bar u = \bar u_1 = 0$.

Finally, if $M(\beta) \cong L_p * L_q$ we have $\pi_1(M(\beta)) \cong
\mathbb Z/p * \mathbb Z/q$. A simple calculation shows that the space
of $1$-cocycles $Z^1(M(\beta); sl_2(\mathbb C)_\rho)$ is isomorphic to
$\mathbb C^4$. Thus $H^1(M(\beta); sl_2(\mathbb C)_\rho) \cong \mathbb
C$ This completes the proof.
\end{proof}

\section{$PSL_2(\mathbb C)$-representations of fundamental groups of very small $3$-manifolds}\label{very small reps}

We begin by considering a $3$-manifold $W$ which fibres over $S^1$
with fibre a torus $T$ and monodromy $A$. It is known that $W$ is a
Sol manifold if and only if $|\mbox{tr}(A)| > 2$ and a Seifert fibred
space otherwise. Similarly if $W$ semi-fibres over the interval with
semi-fibre a torus $T$ and gluing map $A = \left( \begin{array}{cc} a
& b \\ c & d \end{array} \right) \in SL_2(\mathbb Z)$, then $W$ is
a Sol manifold if and only if $ad \ne 0,1$ and a Seifert fibred space
otherwise.

\begin{prop} \label{toroidalreps}
  Suppose that $W$ either fibres over the circle with torus fibre or
  semi-fibres over the interval with torus semi-fibre.  If $\rho:
  \pi_1(W) \to PSL_2(\mathbb C)$ is irreducible, then up to
  conjugation, the image of $\rho$ is $T_{12}$, or $O_{24}$, or lies
  in ${\mathcal N}$.  Moreover,
\begin{itemize}
\item if the image is $T_{12}$, then $\rho(\pi_1(T)) = \mathbb Z/2 \oplus \mathbb Z/2$ and
$W$ fibres over $S^1$;
\item if the image is $O_{24}$, then $\rho(\pi_1(T)) = \mathbb Z/2
  \oplus \mathbb Z/2$ and $W$ semi-fibres over the interval.
\end{itemize}
\end{prop}

\begin{proof}
  Let $T$ denote the (semi-)fibre and consider the normal subgroup $G
  = \rho(\pi_1(T))$ of $\rho(\pi_1(W))$. We can conjugate $G$ so that
  it equals $\mathbb Z/2 \oplus \mathbb Z/2 \subset {\mathcal N}$, or
  it is contained in either ${\mathcal P}$, the group of
  upper-triangular parabolic matrices, or ${\mathcal D}$, the group of
  diagonal matrices.
  
  If $G = \mathbb Z/2 \oplus \mathbb Z/2$, a simple calculation
  implies that $\rho(\pi_1(W))$ is finite.  The only finite subgroups
  of $PSL_2(\mathbb C)$ which contain such a normal subgroup are
  $T_{12}, O_{24}$, and the dihedral group $D_{2} \subset {\mathcal
    N}$.  The first possibility is ruled out when $T$ separates
  $M(\alpha)$ into two twisted $I$-bundles over the Klein bottle,
  since otherwise $\rho$ would induce a surjection of $\mathbb Z/2 *
  \mathbb Z/2 = \pi_1(W)/\pi_1(T)$ onto $T_{12}/(\mathbb Z/2 \oplus
  \mathbb Z/2) = \mathbb Z/3$, which is impossible. Similarly if $T$
  does not separate, then the image of $\rho$ cannot be $O_{24}$.
  
  Next we can rule out the possibility that $\{\pm I\} \ne G \subset
  {\mathcal P}$ since if this case did arise, the normality of $G$ in
  $\rho(\pi_1(W))$ would then imply that $\rho$ is reducible.
  
  Finally assume that $G \subset {\mathcal D}$. If $G = \{\pm I\}$,
  then $\rho$ factors through $\pi_1(W)/\pi_1(T)$ which is isomorphic
  to either $\mathbb Z$ or $\mathbb Z/2 * \mathbb Z/2$. The
  irreduciblity of $\rho$ excludes the former possibility while the
  lemma clearly holds in the latter. If $\{\pm I\} \ne G \subset
  {\mathcal D}$ is non-trivial, then its normality in $\rho(\pi_1(W))$
  implies that the latter is a subset of ${\mathcal N}$.
\end{proof}

\begin{prop} \label{bundlereps}
  Let $W$ be a torus bundle over $S^1$ with monodromy $A \in
  SL_2(\mathbb Z)$ and fibre $T$.  Consider a representation $\rho:
  \pi_1(W) \to PSL_2(\mathbb C)$ which is either
  irreducible or has non-Abelian image.
\begin{enumerate}
\item[(1)] If $\rho$ is irreducible, then $H^1(W; sl_2(\mathbb C)_{Ad
    \rho}) = 0$ as long as $tr(A) \ne -2$.
\item[(2)] If $\rho$ is reducible and $W$ fibres over the circle and
  the image of $\rho$ contains non-trivial torsion, then it is Seifert
  fibred.  Moreover, if there is torsion of order greater than $2$,
  then $|tr(A)| \leq 1$.
\end{enumerate}
\end{prop}

\begin{proof}
  Write $A = \left( \begin{array}{cc} a & b \\ c & d \end{array}
  \right)$ and recall that there is a presentation of $\pi_1(W)$ of
  the form
  $$\langle x,y,t \; | \; [x,y] = 1, txt^{-1} = x^ay^c, tyt^{-1} =
  x^by^d \rangle$$
  where $x, y$ generate $\pi_1(T)$ and $t$ projects
  to a generator $\bar t$ of $\pi_1(S^1) \cong \mathbb Z$.
  
  (1) Consider the exact sequence $1 \to \pi_1(T) \to \pi_1(W) \to
  \mathbb Z \to 1$. The Lyndon-Serre spectral sequence yields an
  associated exact sequence in cohomology
  $$
  0 \to H^1(\mathbb Z; (sl_2(\mathbb C)_{Ad \rho})^{\pi_1(T)}) \to
  H^1(\pi_1(W); sl_2(\mathbb C)_{Ad \rho}) \to H^1(T; sl_2(\mathbb
  C)_{Ad \rho})^{\mathbb Z} \to 0$$
  Since $\rho$ is irreducible, we
  have either $\rho(\pi_1(T)) = \mathbb Z/2 \oplus \mathbb Z/2$, or
  $\{\pm I\} \ne \rho(\pi_1(T)) \subset {\mathcal D}$ and
  $\rho(\pi_1(W)) \subset {\mathcal N}$.
  
  If $\rho(\pi_1(T)) = \mathbb Z/2 \oplus \mathbb Z/2$, then
  $(sl_2(\mathbb C)_{Ad \rho})^{\pi_1(T)} = 0$.  On the other hand,
  using duality with twisted coefficients and the fact that $\chi(T;
  Ad \rho) = 3 \chi(T) = 0$, we see that the associated Betti numbers
  satisfy $b_1(T; sl_2(\mathbb C)_{Ad \rho}) = 2b_0(T; sl_2(\mathbb
  C)_{Ad \rho})$. But since $\rho|\pi_1(T)$ is irreducible, we have
  $b_0(T; sl_2(\mathbb C)_{Ad \rho}) = 0$. Thus $H^1(T; sl_2(\mathbb
  C)_{Ad \rho})^{\mathbb Z} = 0$ which implies the desired result.
  
  Next suppose that $\{\pm I\} \ne \rho(\pi_1(T)) \subset {\mathcal
    D}$ and $\rho(\pi_1(W)) \subset {\mathcal N}$.  In this case
  $$(sl_2(\mathbb C)_{Ad \rho})^{\pi_1(T)} = \{ \left(
    \begin{array}{cc} z & 0 \\ 0 & -z \end{array} \right) \; | \; z
  \in \mathbb C\} \cong \mathbb C.$$
  The irreducibility of $\rho$
  implies that up to conjugation we may suppose that $\rho(t) = \pm
  \left(\begin{array}{cc} 0 & 1 \\ -1 & 0 \end{array} \right)$ and
  therefore $\mathbb Z$ acts on $(sl_2(\mathbb C)_{Ad
    \rho})^{\pi_1(T)}$ by multiplication by $-1$.  Thus the set of
  invariants of this action, which is isomorphic to $H^0(\mathbb Z;
  (sl_2(\mathbb C)_{Ad \rho})^{\pi_1(T)})$, is $0$.  Duality then
  yields $H^1(\mathbb Z; (sl_2(\mathbb C)_{Ad \rho})^{\pi_1(T)}) = 0$.
  
  On the other hand, it is easy to see that $H^1(\pi_1(T);
  sl_2(\mathbb C)_{Ad \rho})$ may be identified with the set of
  homomorphisms of $\pi_1(T)$ into $\mathbb C$ in such a way that if
  $f$ is such a homomorphism, then $\bar t$ acts on $f$ as
  $$(\bar t \cdot f)(x^my^n) = - f(x^{am + bn}y^{cm + dn}) = -(am +
  bn)f(x) - (cm + dn)f(y).$$
  Hence $f$ is invariant under the action
  of $\bar t$ if and only if $(f(x), f(y))$ is a (-1)-eigenvector of
  the transpose of $A$. It follows that $H^1(\pi_1(W); sl_2(\mathbb
  C)_{Ad \rho}) \cong H^1(T; sl_2(\mathbb C)_{Ad \rho})^{\mathbb Z}
  \ne 0$ if and only if $tr(A) = -2$.
  
  (2) Write $A = \left( \begin{array}{cc} a & b \\ c & d \end{array}
  \right)$.  As $\rho$ is reducible with non-Abelian image, we must
  have $\{\pm I\} \ne \rho(\pi_1(T)) \subset {\mathcal P} \cong
  \mathbb C$.  Then the image of $\rho$ lies in ${\mathcal U}$.
  Suppose that
  $$\rho(x) = \pm \left( \begin{array}{cc} 1 & \sigma \\ 0 & 1
\end{array} \right),
\rho(y) = \pm \left( \begin{array}{cc} 1 & \tau \\ 0 & 1 \end{array} \right),
\rho(t) = \pm \left( \begin{array}{cc} u & v \\ 0 & u^{-1} 
\end{array} \right).$$
Since the kernel of the projection ${\mathcal U} \to {\mathcal D}$ is
${\mathcal P} \cong \mathbb C$, any torsion element in the image of
$\rho$ is sent to an element of the same order in ${\mathcal D}$ under
this projection. On the other hand, since any element of $\pi_1(W)$
can be written as a product of the form $x^ly^mt^n$, the image of
$\rho(\pi_1(W))$ under the projection to ${\mathcal D}$ is isomorphic
to $\{u^n \; | \; n \in \mathbb Z\} \subset \mathbb C^*$.  Thus
$\rho(\pi_1(W))$ contains a non-trivial torsion element if and only if
$u$ is a non-trivial root of unity.  Assume this occurs.  The
relations in the presentation for $\pi_1(W)$ imply that
$$u^2 \sigma = a \sigma + c \tau, \;\;\; u^2 \tau = b \sigma + d \tau.
$$
Thus $u^2$ is an eigenvalue of $A$. It is well known that these
eigenvalues are roots of unity if and only if $|tr(A)| \leq 2$.
Moreover when $tr(A) = 2$ we have $u = \pm 1$, when $tr(A) = -2$ we
have $u = \pm i$.  Thus the proposition holds.  \end{proof}

\begin{prop} \label{semibundlereps}
  Let $W$ semi-fibre over the interval with semi-fibre $T$. If there
  is a representation $\rho: \pi_1(W) \to PSL_2(\mathbb C)$ which is
  reducible and has non-Abelian image, then the torsion elements in
  the image of $\rho$ have order $2$.
\end{prop}

\begin{proof}
  Now $W$ splits along $T$ into two twisted $I$-bundles over the Klein
  bottle. Thus there is a presentation of $\pi_1(W)$ of the form
  $$\langle x_1, y_1, x_2, y_2 \; | \; x_1y_1x_1^{-1} = y_1^{-1},
  x_2y_2x_2^{-1} = y_2^{-1}, x_1^2 = x_2^{2a}y_2^c, y_1 =
  x_2^{2b}y_2^d \rangle$$
  where $x_1, y_1$ generate the fundamental
  group of one of the twisted $I$-bundles, $x_2, y_2$ generate the
  fundamental group of the other, and $A = \left( \begin{array}{cc} a
      & b \\ c & d
  \end{array} \right)$ is the gluing matrix. Note that $\pi_1(T)$ is
generated by either pair $x_1^2, y_1$ and $x_2^2, y_2$.

We can suppose that either $\{\pm I\} \ne \rho(\pi_1(T)) \subset
{\mathcal P}$ or $\rho(\pi_1(T)) \subset {\mathcal D}$. In the latter
case, $\rho(\pi_1(T)) = \{\pm I\}$ as otherwise the normality of
$\pi_1(T)$ in $\pi_1(W)$ and the reducibility of $\rho$ imply that
$\rho(\pi_1(W)) \subset {\mathcal D}$.

Assume first that $\rho(\pi_1(T)) = \{\pm I\}$. Then $\rho(x_1)^2 =
\rho(y_1) = \rho(x_2)^2 = \rho(y_2) = \pm I$. Note that neither
$\rho(x_1) = \pm I$ nor $\rho(x_2) = \pm I$ as otherwise the image of
$\rho$ would be Abelian. Thus up to conjugation we have $\rho(x_1) =
\pm \left( \begin{array}{cc} i & 0 \\ 0 & -i \end{array} \right)$ and
$\rho(x_2) = \pm \left( \begin{array}{cc} i & 1 \\ 0 & -i \end{array}
\right)$. Thus the only torsion elements in the image of $\rho$ have
order $2$.

Next assume that $\{\pm I\} \ne \rho(\pi_1(T)) \subset {\mathcal P}$.
The relation $x_1y_1x_1^{-1} = y_1^{-1}$ implies that exactly one of
$x_1^2, y_1$ is sent to $\pm I$ by $\rho$. If $\rho(x_1^2) = \pm I$,
then up to conjugation, $\rho(y_1) = \pm \left( \begin{array}{cc} 1 &
    1 \\ 0 & 1 \end{array} \right)$. Hence, as $x_2^2 =
x_1^{2d}y_1^{-c}$ and $y_2 = x_1^{-2b}y_1^{a}$, we have $\rho(x_2^2) =
\pm \left( \begin{array}{cc} 1 & -c \\ 0 & 1
\end{array} \right),
\rho(y_2) = \pm \left( \begin{array}{cc} 1 & a \\ 0 & 1 \end{array} \right).$
Thus the image of $\rho$ is generated by the images of $x_1, y_1$ and 
$x_2$. Projecting into ${\mathcal D}$
then shows that the only non-trivial torsion elements in the image of 
$\rho$ must have order $2$.
If $c = 0$, then $ad = 1$ and so $W$ is Seifert fibred. On the other 
hand, if $c \ne 0$,
then $\rho(x_2) = \pm \left( \begin{array}{cc} 1 & -\frac{c}{2} \\ 0 
& 1 \end{array} \right)$
and so the relation $x_2y_2x_2^{-1} = y_2^{-1}$ implies that
$a = 0$. Therefore $ad = 0$ and $W$ is Seifert fibred. A similar argument shows that the proposition holds when $\rho(y_1) = \pm I$.
\end{proof}

\begin{lemma} \label{ngdrep}
  Let $W$ be a closed, connected, orientable, irreducible, very small
  $3$-manifold which is not virtually Haken.  Then the image of any
  representation $\rho: \pi_1(W) \to PSL_2(\mathbb C)$ is a finite
  group.
\end{lemma}

\begin{proof}
  Let $\rho: \pi_1(W) \to PSL_2(\mathbb C)$ be a representation.  The
  Tits alternative implies that there is a finite index subgroup $G$
  of $\rho(\pi_1(W))$ which is solvable. It suffices to show that $G$
  is finite.
  
  If $G = \{\pm I\}$ we are done so assume otherwise. Then since $G$
  is solvable it contains a non-trivial normal subgroup $A$ which is
  Abelian. Up to conjugation $A$ is either contained in ${\mathcal
    D}$, or in ${\mathcal P}$, or is the Klein $4$-group $\mathbb Z/2
  \oplus \mathbb Z/2$ realized in $PSL_2(\mathbb C)$ as
$$D_2 = \{\pm I, \pm \left(\begin{array}{cc} i & 0 \\ 0 & -i 
\end{array}\right),
\pm \left(\begin{array}{cc} 0 & 1 \\ -1 & 0 \end{array}\right), \pm
\left(\begin{array}{cc} 0 & i \\ i & 0 \end{array}\right)\}.$$
Since
$A \ne \{\pm I\}$ is normal in $G$, it follows that $A \subset
{\mathcal N}$ if the first or third possibilities arise. In these
cases let $A_0 = G \cap {\mathcal D}$ and observe that $A_0$ is
Abelian and has index at most $2$ in $G$.  Then $A_0$ has finite index
in $\rho(\pi_1(W))$ and, since $W$ is not virtually Haken, must
therefore be finite. But then $\rho(\pi_1(W))$ is finite and we are
done.

On the other hand, suppose that $A \subset {\mathcal P}$. Then the
non-triviality of $A$ and its normality in $G$ imply that $G \subset
{\mathcal U}$, the group of upper-triangular matrices in
$PSL_2(\mathbb C)$. Since each finite degree cover $\tilde W$ of $W$
is irreducible but not Haken, it has zero first Betti number. Thus the
projection of $G$ in ${\mathcal D}$ is finite and so the kernel of
this projection is of finite index in $G$. But this kernel lies in
${\mathcal P}_+$, the subgroup of ${\mathcal P}$ consisting of
matrices of trace $2$. Since this group is isomorphic to $\mathbb C$,
and again using the fact that $W$ is not virtually Haken we see that
the kernel is trivial. Thus $G$ is finite.
\end{proof}

We now apply the results above to 

\begin{prop} \label{fin-nvh}
  Suppose that $X_0 \subset X_{PSL_2}(M)$ is a non-trivial curve and
  $\alpha$ is a slope on $\partial M$ which is not a singular slope
  for any closed, essential surface in $M$.  If $M(\alpha)$ either has
  a finite fundamental group, or is an irreducible, very small
  $3$-manifold which is not virtually Haken, then
\begin{enumerate}
\item[(1)] $J_{X_0}(\alpha) \subset X_0^\nu$;
\item[(2)] for each $x \in J_{X_0}(\alpha)$, there is an irreducible
  representation $\rho$ with finite image such that $\chi_\rho =
  \nu(x)$, $\rho(\pi_1(\partial M)) \ne \{\pm I\}$ and $H^1(M(\alpha);
  sl_2(\mathbb C)_\rho) = 0$;
\item[(3)] if $x \in J_{X_0}(\alpha)$ then $\nu(x)$ is a simple
point of $X_{PSL_2}(M)$.  
\end{enumerate}
\end{prop}

\begin{proof}
  Our hypotheses imply that $b_1(M) = 1$.  Thus Proposition
  \ref{jumpideal} implies that $J_{X_0}(\alpha) \subset X_0^\nu$.
  Consider $x \in J_{X_0}(\alpha)$ and suppose that $\chi_\rho =
  \nu(x)$.
  
  If $\pi_1(M(\alpha))$ is finite, then so is the image of $\rho$. The
  same conclusion holds when $\pi_1(M(\alpha))$ is not finite by Lemma
  \ref{ngdrep}.
  
  Suppose next that $\rho$ is reducible. Since its image is finite, it
  is conjugate to a diagonal representation and as this is true for
  each representation in $t^{-1}(\nu(x)$, any two representations in
  $t^{-1}(\nu(x))$ are conjugate.  Hence the dimension of
  $t^{-1}(\nu(x))$ is at most $2$, contrary to Corollary 1.5.3 of
  [CS]. This shows that $\rho$ is irreducible. The fact that
  $\rho(\pi_1(\partial M)) \ne \{\pm I\}$ can now be proven in exactly
  the same way as Lemma 4.2 of [BZ1].
  
  Next we show that $H^1(M(\alpha), sl_2(\mathbb C)_\rho) = 0$.  Let
  $G = \rho(\pi_1(M(\alpha))$ and consider the left
  $\pi_1(M(\alpha))$-module $\mathbb C[G]_\rho$. It is well known that
  $\mathbb C[G]$ splits as a direct sum $\oplus_\sigma V_\sigma$ of
  irreducible $\mathbb C G$-modules $V_\sigma$ and each irreducible
  $\mathbb C G$-module appears at least once in this decomposition
  [Ser]. On the other hand if $W \to M(\alpha)$ is the finite cover
  corresponding to the kernel of $\rho$, our hypotheses imply that
  $H^1(W; \mathbb C) = 0$. This is obvious if $\pi_1(M(\alpha))$ is
  finite and follows from the fact that $W$ is irreducible and
  non-Haken otherwise. Thus
  $$0 = H^1(W; \mathbb C) = H^1(M(\alpha); \mathbb C[G]_\rho) =
  \oplus_\sigma H^1(M(\alpha); (V_\sigma)_\rho).$$
  This shows that for
  any irreducible $\mathbb C[G]$-module $V$, $H^1(M(\alpha); V_\rho) =
  0$ and therefore, $H^1(M(\alpha), sl_2(\mathbb C)_\rho) = 0$ as
  claimed.
  
  Finally, we note that, according to [BZ4, Theorem 3], conditions (1)
  and (2) imply that $\nu(x)$ is a simple point of $X_{PSL_2}(M)$.
\end{proof}

\begin{prop} \label{pssms}
  Let $X_0$ be a non-trivial curve in $X_{PSL_2}(M)$ and $\alpha$ a
  slope on $\partial M$ such that $f_\alpha|X_0 \not= 0$.  Suppose
  that $\alpha$ is not a singular slope for a closed, essential
  surface in $M$. Assume as well that either
\begin{enumerate}
\item[(i)] $\pi_1(M(\alpha))$ is finite or $M(\alpha)$ is an
  irreducible very small $3$-manifold which is not virtually Haken, or
\item[(ii)] $M(\alpha)$ is a non-Haken Seifert manifold with base
  orbifold of the form $S^2(r,s,t)$ and there is a slope $\beta$ on
  $\partial M$ such that $M(\beta) \cong S^1 \times S^2$, or
\item[(iii)] $X_0 \subset X_{PSL_2}(M(\beta))$ where $\beta$ is a
  slope on $\partial M$ such that $M(\beta) \cong L_p \# L_q$.
\end{enumerate}
Then
$$\|\alpha\|_{X_0} = m_0 + 2|J_{X_0}(\alpha)| - A$$
where $m_0 =
\sum_{x \in \widetilde X_0} \mbox{min}\{ Z_x(\tilde f_\beta) \; | \;
\tilde f_\beta| \widetilde X_0 \not= 0\}$, and $A$ is the number of
irreducible characters $\chi_\rho \in \nu(J_{X_0}(\alpha))$ of
representations $\rho$ which are conjugate into ${\mathcal N}$.
\end{prop}

\begin{proof}
  Case (ii) is done in [Theorem 2.3, BB] while the proof in case (i)
  is handled analogously. The idea is that by combining
  (\ref{normformula}), Proposition \ref{4.10-4.12}, and the previous
  two propositions, the calculation of $\|\alpha\|_{X_0}$ reduces to a
  weighted count of characters of representations $\pi_1(M(\alpha))
  \to PSL_2(\mathbb C)$.  Note that under our assumptions,
  $J_{X_0}(\alpha) \subset X_0^\nu$ and $\nu|J_{X_0}(\alpha)$ is
  injective.
  
  Finally, for case (iii), Proposition \ref{jumpideal} implies that
  $J_{X_0}(\alpha) \subset X_0^\nu$ and a calculation similar to that
  used in case (i) yields the desired conclusion.
\end{proof}

\section{Proof of Proposition \ref{TLF}}\label{PTLF}\label{proof one}

We suppose in this section that $b_1(M) = 1$, that neither $\alpha$
nor $\beta$ is a singular slope for a closed, essential surface in
$M$, that $M(\alpha)$ has a finite fundamental group, and that
$M(\beta)$ is either a connected sum of two lens spaces or $S^1 \times
S^2$. Theorem \ref{2.0.3} implies that $\alpha$ is not a boundary
slope.

A finite filling slope $\alpha$ is either of $C$-type or $D$-type or
$Q$-type or $T(k)$-type ($1 \leq k \leq 3$) or $O(k)$-type ($1 \leq k
\leq 4$) or $I(k)$-type ($1 \leq k \leq 5, k \ne 4$). We refer to
[BZ5, Pages 93-94 and 98] for these definitions. We will show
$$\Delta(\alpha, \beta) \leq \left\{ \begin{array}{ll}
2 & \mbox{{\rm if $M(\beta) \cong L_2 \# L_3$, $H_1(M) \cong \mathbb 
Z \oplus \mathbb Z/2$ and $\alpha$ is of type
$O(2)$;}} \\  1 & \mbox{{\rm otherwise.}}\end{array} \right.$$
The key relationships between Culler-Shalen seminorms and finite filling
classes is contained in the following result from [BZ5].

\begin{prop} \label{fin}
  Suppose that $X_0$ is a non-trivial curve in $X_{PSL_2}(M)$ and that
  $\alpha$ is a finite or cyclic filling slope which is not a boundary
  slope associated to an ideal point of $X_0$.
\begin{enumerate}
\item[(1)] \mbox{ {\rm ([CGLS])} } If $\alpha$ is a
  cyclic filling slope then $\|\alpha\|_{X_0}=s_{X_0}$.
\item[(2)] If $\alpha$ is a $D$-type or a $Q$-type filling slope and
  $X_0$ is not virtually trivial then
\begin{enumerate}
  \item[(i)] $\|\alpha\|_{X_0}\leq 2s_{X_0}$;
  \item[(ii)] $\|\alpha\|_{X_0}\leq \|\beta\|_{X_0}$ for any slope
  $\beta$ such that $\Delta(\alpha, \beta) \equiv 0 \;
  \mbox{{\rm (mod $2$)}}$.
\end{enumerate}
\item[(3)] If $\alpha$ is a $T(k)$-type filling slope, then $k \in
  \{1, 2, 3\}$ and
  \begin{enumerate}
  \item[(i)]  $\|\alpha\|_{X_0}\leq s_{X_0}+2$;
  \item[(ii)] $\|\alpha\|_{X_0}\leq \|\beta\|_{X_0}$ for any slope
    $\beta$ such that $\Delta(\alpha, \beta)\equiv 0 \; \mbox{{\rm
        (mod $k$)}}$.
  \end{enumerate}
\item[(4)] If $\alpha$ is an $O(k)$-type filling slope, then $k \in
  \{1, 2, 3, 4\}$ and
  \begin{enumerate}
  \item[(i)] $\|\alpha\|_{X_0}\leq s_{X_0}+3$;
  \item[(ii)] $\|\alpha\|_{X_0}\leq \|\beta\|_{X_0}$ for any slope
    $\beta$ such that $\Delta(\alpha, \beta) \equiv 0 \; \mbox{{\rm
        (mod $k$)}}$.
  \end{enumerate}
\item[(5)] If $\alpha$ is an $I(k)$-type filling slope, then $k \in
  \{1, 2, 3, 5\}$ and
\begin{enumerate}
  \item[(i)]  $\|\alpha\|_{X_0}\leq s_{X_0}+4$;
  \item[(ii)] $\|\alpha\|_{X_0}\leq \|\beta\|_{X_0}$ for any slope
    $\beta$ such that $\Delta(\alpha,\beta) \equiv 0 \; \mbox{{\rm
        (mod $k$)}}$.
  \end{enumerate}
\end{enumerate}
\end{prop}
\qed

We split the proof of Proposition \ref{TLF} into three cases.

\begin{ccase}\label{two lens spaces}
$M(\beta) \ne P^3 \# P^3$  is a connected sum of two lens spaces.
\end{ccase}

Recall that $X_{PSL_2}(M(\beta)) \subset X_{PSL_2}(M)$ contains
exactly $[\frac{p}{2}][\frac{q}{2}]$ non-trivial curves $X(j,k)$,
where $1 \leq j \leq \frac{p}{2}$ and $1 \leq k \leq \frac{q}{2}$.
Let $X$ be the union of these curves and observe that since $\beta$ is
not a singular slope for any closed essential surface in $M$,
Proposition \ref{lb} implies that
\begin{equation} s_X \geq s_0 = \left\{ \begin{array}{ll} (p-1)(q-1) + 1 &
\mbox{if $p,q$ even} \\ (p-1)(q-1) & \mbox{otherwise.}
\end{array} \right. \end{equation}
By (\ref{minformula}), $\|\alpha\|_X = \Delta(\alpha, \beta)s_X$.  If
$\alpha$ is a $C$-type filling slope, then $\|\alpha\|_X \leq s_X$ by
Proposition \ref{fin} (recall that $\alpha$ is not a boundary slope)
and therefore $\Delta(\alpha, \beta) \leq 1$.  If it is a $D$ or
$Q$-type filling slope, then all irreducible representations of
$\pi_1(M(\alpha))$ conjugate into ${\mathcal N}$. Thus Propositions
\ref{fin-nvh} and \ref{pssms} show that for each $x \in
J_{X}(\alpha)$, $\nu(x)$ is an irreducible character and $\|\alpha\|_X
\leq s_X + |\nu(J_{X}(\alpha))|$.  On the other hand, Lemma
\ref{charcount} shows that if $X(j,k)$ is a component of $X$ with $j$
and $k$ relatively prime to $p$ and $q$ respectively, then it contains
the character of an irreducible representation with image in
${\mathcal N}$ if and only if $p = 2$, and if $p = 2$, there is a
unique such character. Hence $\Delta(\alpha, \beta)s_{X} =
\|\alpha\|_{X} \leq s_{X} + [\frac{p}{2}][\frac{q}{2}] < 2 s_{X}$, and
therefore $\Delta(\alpha, \beta) \leq 1$.

Next assume that $\alpha$ is either a $T$ or $O$ or $I$-type
filling slope. Then (\ref{minformula}) and Proposition \ref{fin} show that
$$\Delta(\alpha, \beta) \leq \left\{
\begin{array} {ll}
1 + \frac{2}{s_X} \leq 1 + \frac{2}{s_0} &
\mbox{if $\alpha$ is $T$-type} \\ 1 +
\frac{3}{s_X} \leq 1 + \frac{3}{s_0} & \mbox{if
$\alpha$ is $O$-type} \\ 1 + \frac{4}{s_X} \leq
1 + \frac{4}{s_0} & \mbox{if $\alpha$ is $I$-type.}
\end{array} \right. $$
Combining this inequality with (8.1.1) shows that $\Delta(\alpha, 
\beta) \leq 1$ unless, perhaps,
\begin{itemize}
\item  $(p,q) = (2,4), (2,5), (3,3), \Delta(\alpha, \beta) \leq 
2$ and $\alpha$ is $I$-type, or
\item  $(p, q) = (2,3), \Delta(\alpha, \beta) \leq 3$ and 
$\alpha$ is $I$-type, or
\item  $(p, q) = (2,3), \Delta(\alpha, \beta) \leq 2$ and 
$\alpha$ is either $T$ or $O$-type.
\end{itemize}
Assume first that one of these cases arises and $\alpha$ is either of
type $T$ or $I$. It is well-known that
\begin{equation} H_1(M(\alpha)) \cong \left\{\begin{array}{l}
\mathbb Z/3^kj \mbox{ $k\geq 1$ and $j$ relatively prime to 
$6$, if $\alpha$ is $T$-type}\\
\mathbb Z/j \mbox{ where $j$ is relatively prime to $30$, if $\alpha$ 
is $I$-type}
\end{array}\right. \end{equation} 
(see [BZ1] for instance) so that in each of these cases,
$H_1(M(\alpha))$ is cyclic. This implies that $H_1(M) \cong \mathbb Z
\oplus \mathbb Z/n$ where $n \geq 1$. Then $n$ divides
$|H_1(M(\delta)|$ for each primitive $\delta \in H_1(\partial M)$.
Taking $\delta = \alpha$ we see that $n$ divides $3^kj$ where
$\gcd(j,6) = 1$ if $\alpha$ is $T$-type, and divides $j$ where
$\gcd(j,30) = 1$ if $\alpha$ is $I$-type.  On the other hand, $n$ also
divides $|H_1(M(\beta))| \cong \mathbb Z/p \oplus \mathbb Z/q$, so
given the constraints we have imposed on $(p,q)$ we see that $n = 1$.
Thus $H_1(M) \cong \mathbb Z$, and so each Dehn filling of $M$ has a
cyclic first homology group.  This rules out the possibility that
$(p,q) = (2,4)$ or $(3,3)$.  Consider then the cases where $(p,q) =
(2,3)$ or $(2,5)$. There is a basis $\{\mu, \lambda\}$ of
$H_1(\partial M)$ such that $\lambda$ is zero homologically in $M$ and
$\mu$ generates $H_1(M)$.

If $\alpha$ is a $T$-type filling slope, then $(p,q) = (2,3)$ and so
by our choice of $\mu$ and $\lambda$, (8.1.2) implies that there are
integers $a, b$ such that up to sign, $\alpha = \{\pm(3^kj \mu + a
\lambda)\}$, and $\beta = \{\pm(6 \mu + b \lambda)\}$. Then $b$ is odd
and the constraints on $j,k$ show that $\Delta(\alpha, \beta) = |6a -
3^kjb| \equiv 1 \mbox{ {\rm (mod 2)}}$. As $\Delta(\alpha, \beta) \leq
2$, we have $\Delta(\alpha, \beta) = 1$.

Next suppose that $\alpha$ is an $I$-type filling class.  Then $(p,q)
= (2,3)$ or $(2,5)$. By (8.1.2) there are integers $a, b$ such that
$\alpha = \{\pm(j \mu + a \lambda)\}$, and $\beta = \{\pm(2q \mu + b
\lambda)\}$. Then $b$ is relatively prime to $2q$ and since $\gcd(j,
30) = 1$, $\Delta(\alpha, \beta) = |2qa - jb|$ is relatively prime to
$2q$ as well.  When $q = 5$, this shows that $\Delta(\alpha, \beta)$
is odd, and therefore as $\Delta(\alpha, \beta) \leq 2$ in this case,
we have $\Delta(\alpha, \beta) = 1$.  Finally when $q = 3$, it shows
that $\Delta(\alpha, \beta)$ is relatively prime to $6$, and therefore
as $\Delta(\alpha, \beta) \leq 3$, we have $\Delta(\alpha, \beta) =
1$.

Finally suppose that $\Delta(\alpha, \beta) = 2$, $(p,q) = (2,3)$, and
$\alpha$ has type $O$.  Now $H_1(M(\alpha)) \cong \mathbb Z/2j$ where
$j$ is relatively prime to $6$ ([BZ1]) and we can argue as above to
see that either $H_1(M) \cong \mathbb Z$ or $H_1(M) \cong \mathbb Z
\oplus \mathbb Z/2$. When $H_1(M) \cong \mathbb Z$, we can find, as
above, a basis $\mu, \lambda$ of $H_1(\partial M)$ such that $\lambda$
is zero homologically in $M$ and $\mu$ generates $H_1(M)$.  There are
integers $a, b$ such that $\alpha = \{\pm(2j \mu + a \lambda)\}$, and
$\beta = \{\pm(6 \mu + b \lambda)\}$ where $a$ and $b$ are odd. Since
$j$ is odd as well, we have $2 \geq \Delta(\alpha, \beta) = |6a - 2bj|
\equiv 0 \mbox{ {\rm (mod $4$)}}$. This contradiction shows that this
case does not arise.  If $H_1(M) \cong \mathbb Z \oplus \mathbb Z/2$
and $\Delta(\alpha, \beta) = 2$, the fact that $\alpha$ has type
$O(2)$ follows from Proposition \ref{fin} (4).

Thus in all cases, $\Delta(\alpha, \beta) \leq 2$ and $\Delta(\alpha,
\beta) \leq 1$ unless, perhaps, $H_1(M) \cong \mathbb Z \oplus \mathbb
Z/2$, $M(\beta) \cong L_2 \# L_3$ and $M(\alpha)$ has type $O(2)$.
This completes the proof in Case \ref{two lens spaces}.
\qed

\begin{ccase} \label{Pthree plus Pthree}
$M(\beta)= P^3\# P^3$
\end{ccase}
We show that in this case, $\Delta(\alpha, \beta) \leq 1$. First we
need some auxilliary results.  

Note that there is a $2$-fold cover $p: \widetilde M_\beta \to M$
obtained by restricting the cover $S^1 \times S^2 \to P^3 \# P^3\cong
M(\beta)$.  Let $\phi_\beta : \pi_1(M) \to \mathbb Z/2$ be the
associated homomorphism.  Note also that $|\partial \widetilde
M_\beta| \in \{1,2\}$.

\begin{prop} \label{index1-2}
  Suppose that $M(\beta) \cong P^3 \# P^3$ and that $\beta=\{\pm b\}$
  is not a strict boundary slope. Suppose that $X_0 \subset X(M)$ is a
  curve which is not virtually trivial and that $\|\beta\|_{X_0} \ne
  0$. Then there is an index $\frac{2}{|\partial \widetilde M_\beta|}$
  sublattice $\widetilde L$ of $H_1(\partial M)$ containing $b$ such
  that $\|\beta\|_{X_0} \leq \|\alpha\|_{X_0}$ for each slope
  $\alpha=\{\pm a\}$ where $a \in \widetilde L$ and $\|\alpha\|_{X_0}
  \ne 0$. In particular $\|\beta\|_{X_0} \leq
  \frac{2s_{X_0}}{|\partial \widetilde M_\beta|}$.
\end{prop}

\begin{proof}
 The proof is identical to the proof of Theorem 2.1(a) of [BZ1]. In
that result a non-strict boundary slope $\beta_0$ on $\partial M$ was
given along with a cover $\widetilde{M(\beta_0)} \to M(\beta_0)$ where
$\pi_1(\widetilde{M(\beta_0)})$ is a finite cyclic group. Let $p_0:
\widetilde M \to M$ be the associated cover of $M$ and $T$ be a
boundary component of $\widetilde M$.  It was shown in [BZ1] that if
$\widetilde L = (p_0|T)_*(H_1(T))$.  that for any slope $\{\pm a\}$
such that $ a \in \widetilde L$ and $\|\alpha|_{X_0} \ne 0$, we have
$\|\beta_0\|_{X_0} \leq \|\alpha\|_{X_0}$. The reader can readily
verify that the proof works equally well in the case where
$\pi_1(\widetilde{M(\beta_0)})$ is an infinite cyclic group, the
situation we are considering. Let $T$ be a boundary component of the
double cover $p: \widetilde M_\beta \to M$.  If we now set $\widetilde
L = (p|T)_*(H_1(T)$ then for any slope $\alpha = \{\pm a\}$ such that
$a\in\widetilde L$ and $\|\alpha|_{X_0} \ne 0$, we have
$\|\beta\|_{X_0} \leq \|\alpha\|_{X_0}$. The index of $p_*(H_1(T))$ in
$H_1(\partial M)$ is $\frac{2}{|\partial \widetilde M_\beta|}$, so the
conclusions of the proposition hold.
\end{proof}

\begin{cor} \label{singormin}
  Suppose that $M(\beta) \cong P^3 \# P^3$, $\beta$ is not a singular
  slope for a closed essential surface in $M$, and let $C \subset
  \mathbb Z/2 * \mathbb Z/2 = \pi_1(M(\beta))$ be the unique cyclic
  subgroup of index $2$. Then $\pi_1(\partial M)$ is sent to a
  non-trivial subgroup of $C$ under the natural homomorphism $\pi_1(M)
  \to \pi_1(M(\beta))$. Moreover for any curve $X_0 \subset
  X_{PSL_2}(M)$ which {\it is not virtually trivial}, we have
  $\|\beta\|_{X_0} \leq s_{X_0}$.
\end{cor}

\begin{proof}
 Let $\beta^*$ be a dual class to $\beta$ and choose elements $b$
and $b^*$ of $H_1(\partial M)$ with $\beta=\{\pm b\}$ and $\beta^* =
\{\pm b^*\}$.  Identify $\pi_1(\partial M)$ with $H_1(\partial M)$,
and let $\gamma$ denote the image of $b^*$ in $\pi_1(M(\beta))$.  If
$\gamma^2 = 1$, then $f_{b + 2nb^*} = 0$ and so $\|b + 2nb^*\|_{X_0} =
0$ for each $n \in \mathbb Z$. It follows that $\|\cdot\|_{X_0} = 0$,
and so Proposition \ref{4.10-4.12} implies that $\beta$ is a singular
slope for a closed essential surface in $M$, contrary to our
hypotheses. Thus $\gamma$ has infinite order in $\mathbb Z/2 * \mathbb
Z/2 = \pi_1(M(\beta))$. It follows that $\gamma \in C$ and since $b\in
\pi_1(\partial M)$ maps to the identity in $\pi_1(M(\beta))$ we see
that $\pi_1(\partial M)$ is sent to $C$.  Now $C$ is the kernel of the
homomorphism $\pi_1(M(\beta)) \to \mathbb Z/2$ defining the cover $S^1
\times S^2 \to P^3 \# P^3$, and thus $\pi_1(\partial M) \subset
\hbox{ker}(\phi_\beta)$.  It follows that $|\partial \widetilde
M_\beta| = 2$. As $\beta$ is not a strict boundary slope (cf.
Corollary \ref{strictsing}), the previous proposition shows that
$\|\beta\|_{X_0} \leq s_{X_0}$.  \end{proof}

\begin{lemma} \label{dist1}
  Let $X_M \subset X_{SL_2}(M)$ \footnote{$SL_2(\mathbb C)$-character
    varieties and $SL_2(\mathbb C)$ Culler-Shalen seminorms are
    defined in a manner similar to their $PSL_2(\mathbb C)$
    counterparts and possess similar properties. We refer the reader
    to [CGLS].}  be the canonical curve and suppose that the slope
  $\delta$ is not a strict boundary class and satisfies
  $\|\delta\|_{X_M} = s_M$.  Suppose that $\alpha$ is a slope such
  that $\pi_1(M(\alpha))$ is either finite or cyclic or $\mathbb Z/2 *
  \mathbb Z/2$. Then either
  \begin{enumerate}
  \item[(1)] $\alpha$ is a singular slope for a closed essential surface
    in $M$, or
  \item[(2)] $\Delta(\alpha, \delta) \leq 2$ and if $\Delta(\alpha,
    \delta) = 2$, then $\alpha$ is of $T(k), O(k)$ or $I(k)$-type where
    $k \geq 3$.
  \end{enumerate}
\end{lemma}

\begin{proof}
 Suppose that $\alpha$ is not a singular slope for a closed
essential surface in $M$. Then Theorem \ref{2.0.3} and Corollary
\ref{strictsing} imply that it is not a strict boundary slope and
therefore we can apply Proposition 7.2 of [BZ1] to see that
$\Delta(\alpha, \delta) \leq 2$ when $\pi_1(M(\alpha))$ is either
finite or cyclic.  When it is $\mathbb Z / 2 * \mathbb Z / 2$, an
$SL_2(\mathbb C)$ version of Proposition \ref{index1-2} shows that
$\|\alpha\|_{X_M} \leq 2s_M$ and it follows from the basic properties
of $\|\cdot\|_M$ ([BZ1]) that $\Delta(\alpha, \delta) \leq 2$.
 
Suppose then that $\Delta(\alpha, \delta) = 2$ and let $\tau=\{\pm
t\}$ be a dual slope to $\delta=\{\pm d\}$. Then $\alpha = \{\pm(nd +
2t)\}$ for some $n \in \mathbb Z$. Hence $\Delta(\alpha,\delta) \equiv
0 \mbox{ (mod 2)}$ and thus if $\alpha$ is of type $D$ or $Q$, or
$T(k), O(k), I(k)$ where $k \leq 2$, or $\pi_1(M(\alpha)) \cong
\mathbb Z$ or $\mathbb Z/2 * \mathbb Z/2$, then $\|\alpha\|_M \leq
\|\delta\|_M = s_M$ (cf.  Propositions \ref{fin} and \ref{index1-2}).
But it was shown in \S 1.1 of [CGLS] that if the distance between two
slopes of minimal non-zero Culler-Shalen norm is $2$, then both are
strict boundary slopes.  Hence these cases do not arise and so
$\alpha$ is of type $T(k), O(k)$ or $I(k)$-type where $k \geq 3$.
\end{proof}

\noindent {\bf Proof of Proposition \ref{TLF} when $M(\beta)= P^3\# P^3$.} 
Since neither $\alpha$ nor $\beta$ is a singular slope for a closed
essential surface in $M$, they are not strict boundary slopes (see
Theorem \ref{2.0.3}, Corollary \ref{strictsing}). Thus Corollary
\ref{singormin} and Lemma \ref{dist1} show that $\|\beta\|_M = s_M$
and $\Delta(\alpha, \beta) \leq 2$ with equality implying that
$\alpha$ has type $T(k), O(k), I(k)$ where $k \geq 3$. Since
$H_1(M(\beta)) \cong \mathbb Z/2 \oplus \mathbb Z/2$, $H_1(M; \mathbb
Z/2) \supseteq \mathbb Z/2 \oplus \mathbb Z/2$ and so $H_1(M(\alpha);
\mathbb Z/2) \ne 0$. Hence $\alpha$ is neither $T$ or $I$ type (cf.
(8.1.2)). We must consider the possibility that it is of type $O(k)$
where $k = 3, 4$.

Let $X_0 \subset X_{PSL_2}(\mathbb Z/2 * \mathbb Z/2) =
X_{PSL_2}(M(\beta)) \subset X_{PSL_2}(M)$ be the unique non-trivial
curve. According to Proposition \ref{lb}, $\|\cdot\|_{X_0} \ne 0$ and
further $s_{X_0} \geq 2$.  It is easy to verify that the only
irreducible representations of $\pi_1(M) \to PSL_2(\mathbb C)$ with
finite image whose character lies in $X_0$ are ones with dihedral
image. Since $O$-type groups admit only one such character (Lemma 5.3
of [BZ1]), it follows from Proposition \ref{pssms} that
$\Delta(\alpha, \beta)s_{X_0} = \|\alpha\|_{X_0} \leq s_{X_0} + 1$.
Thus $\Delta(\alpha, \beta) \leq 1$ as claimed, which completes the
proof in Case \ref{Pthree plus Pthree}.  \qed

\begin{ccase}
$M(\beta) = S^1 \times S^2$
\end{ccase}
We prove $\Delta(\alpha, \beta) = 1$.

By Theorem \ref{2.0.3}, $\beta$ is not a strict
boundary slope and so Proposition \ref{fin} implies $\|\beta\|_M = s_M$. Thus Lemma 
\ref{dist1} shows that
$\Delta(\alpha, \beta) \leq 2$, and if it equals $2$, then $\alpha$ 
has type $T(q), O(q)$ or $I(q)$
where $q \geq 3$. We assume below that $\Delta(\alpha, \beta) = 2$ in 
order to arrive at a contradiction.
Let $i: \partial M \to M$ be the inclusion.

\begin{obs} \label{basis}
  Let $\beta=\{\pm b\}$.  There is an integer $n \geq 1$ such that
  $H_1(M) \cong \mathbb Z \oplus \mathbb Z/n$ in such a way that
  $i_*(b) = (0,1)$. Moreover there is a dual slope $\beta^* = \{\pm
  b^*\}$ for $\beta$ such that $i_*(b^*) = (n,0)$.
\end{obs}

\begin{proof}
 Since $\alpha$ is a finite filling slope, the first Betti number
of $M$ is $1$.  Since $H_1(M(\beta)) \cong \mathbb Z$, we have $H_1(M)
\cong \mathbb Z \oplus \mathbb Z/n$ where $n \geq 1$ and $i_*(b)$
generates $\mathbb Z/n$, say $i_*(b) = (0, \bar 1) \in \mathbb Z
\oplus \mathbb Z/n$.  Let $b^*_1$ be any dual class for $b$ and
observe that since $i_*$ has rank $1$, we have $i_*(b^*_1) = (d,\bar
k)$ for some integers $d \ne 0$ and $k$. Then $b^* = b^*_1 - kb$ is
also dual to $b$ and satisfies $i_*(b^*) = (d, \bar 0)$. Let $\xi \in
H_1(M)$ correspond to $(1, \bar 0)$.  By our assumptions, there is a
generator $\eta \in H_2(M, \partial M)$ such that $\partial(\eta) = n
b$.  Lefschetz duality implies that $|\xi \cdot \eta| = 1$. Hence
$|d| = |i_*(b^*) \cdot \eta| = |b^* \cdot \partial(\eta)| =
n$.  It follows that $i_*(b^*) = \pm n \xi$, which completes the
proof of Observation \ref{basis}.
\end{proof}

Since $\Delta(\alpha, \beta) = 2$, we can write $a = 2 b^* +
mb$ (up to sign) for some $m \in \mathbb Z$. A homological
calculation now shows that $|H_1(M(\alpha))| = 2n^2$ and so $\alpha$
cannot have type $T$ or $I$. Thus it has type $O$ and so
$\pi_1(M(\alpha)) \cong O^* \times \mathbb Z/j$ where $O^*$ is the
binary octahedral group and $j$ is an integer relatively prime to $6$.
Then $\mathbb Z/2j \cong H_1(M(\alpha)) \cong \mathbb Z/2n^2$. It
follows that $n^2 = j$ and therefore $n$ is odd. Lemma 3.1 (4) of
[BZ5] now shows that $\alpha$ has type $O(4)$.  Thus the image of
$\pi_1(\partial M)$ under the representation $\rho$, given by
composition $\pi_1(M) \to \pi_1(M(\alpha)) \to O_{24} \subset
PSL_2(\mathbb C)$, has image $\mathbb Z/4$. As $\rho(\alpha) = \pm I$
and $\Delta(\alpha, \beta) = 2$, $\rho(\beta)$ is the square of an
element of order $4$ in $O_{24}$.  Thus it lies in the kernel of the
surjective homomorphism $\phi: O_{24} \to D_3$, which sends any
element of order $4$ to an element of order $2$.  Then $\phi \circ
\rho$ induces a surjective homomorphism of $\pi_1(M(\beta)) \cong
\mathbb Z$ onto the non-Abelian group $D_3$, which is impossible.
Thus it must be that $\Delta(\alpha, \beta) \leq 1$.
\qed

\section{Proof of Proposition \ref{weakTLVS}}\label{weakPTLVS}
\label{proof two}
  
Here we suppose that $\beta$ is a strict boundary slope but is not a
singular slope for a closed, essential surface in $M$. It follows from
Theorem \ref{2.0.3} and Corollary \ref{strictsing} that $M(\beta)$ is
not homeomorphic to $P^3\# P^3$ or $S^1 \times S^2$.  The proof of
Proposition \ref{weakTLVS} is therefore a consequence of the
following result which, unlike Proposition \ref{TLVS}, does not assume
that $M(\alpha)$ admits a geometric decomposition.

\begin{prop} \label{TLVSLS}
Suppose that $M(\beta)$ is a connected sum $L_p \# L_q$ of two lens spaces
where $2 \leq p \leq q$ and $2<q$, and that $M(\alpha)$ is an irreducible 
very small $3$-manifold. Then
$$\Delta(\alpha, \beta) \leq
\begin{cases}
3 & \mbox{if $(p,q) \in \{(2,3), (2,5), (3,5)\}$;} \\
2 & \mbox{if $(p,q) \in \{ (2,4), (3,3), (3,4), (5,5)\}$;} \\
1 & \mbox{otherwise}.
\end{cases} $$
\end{prop}

\begin{proof}
 Let $X_0$ be one of the curves $X(j,k) \subset X_{PSL_2}(\mathbb
Z/p * \mathbb Z/q) = X_{PSL_2}(M(\beta)) \subset X_{PSL_2}(M)$ where
$j,k$ are relatively prime to $p, q$ respectively. Suppose that $x \in
J_{X_0}(\alpha)$.  Proposition \ref{fin-nvh} shows that $x \in
X_0^\nu$ and $\nu(x)$ is a simple point of $X_{PSL_2}(M)$ which is the
character of an irreducible representation $\rho$ whose image is a
finite subgroup of $PSL_2(\mathbb C)$.  In particular, this
implies that if $\nu(x) = \chi_\rho$, where $\rho\in \mathcal N$ then
$\rho$ must have dihedral image.

Let $X \subset X_{PSL_2}(M)$ be the union of the curves
$X(j,k) \subset X_{PSL_2}(M(\beta)) \subset X_{PSL_2}(M)$
where $j,k$ are relatively prime to $p, q$. If $d$ is the number of 
components of $X$, then
Proposition \ref{lb} shows that
$$s_{X} \geq \left\{
\begin{array}{ll}
2d & \mbox{{\rm if }}  p = 2\\
4d  & \mbox{{\rm if }}  p > 2. \end{array} \right. $$
Recall from (\ref{minformula}) that $\Delta(\alpha, \beta) = 
\|\alpha\|_{X}/s_X$. On the other hand,
Proposition \ref{pssms} and our discussion 
above show that 
$\|\alpha\|_{X} = s_{X} + 2|J_X(\alpha)| - A$ where $A$ is the number of 
dihedral characters in $\nu(J_{X}(\alpha))$.
According to Lemma \ref{charcount}(2) we have $A = d$ if $p=2$ and
$A = 0$ if $p>2$.  If we set $n =  |J_X(\alpha)|$ then we have
\begin{equation}\label{6.6}
\Delta(\alpha, \beta) = 1 + \frac{2n - A}{s_{X}} \leq
\begin{cases}
1 + \frac{2n - d}{2d} & \mbox{if $p = 2$;} \\
1 + \frac{2n}{4d}  & \mbox{if $p > 2$.}
\end{cases}
\end{equation}

We have $d=[\frac{p}{2}][\frac{q}{2}]$, and $n$ is determined 
by Lemma \ref{charcount} since $\nu(x)$ is the character of an
irreducible representation with finite image for each $x\in J_X(\alpha)$.
Checking each case, we see that
$$
\Delta(\alpha, \beta) \leq
\begin{cases}
5 & \mbox{if $(p,q) = (2,3)$;}\\
3 & \mbox{if $(p,q) \in \{(2,5),  (3,5)\}$;} \\
2 & \mbox{if $(p,q) \in \{(2, 4), (3, 3), (3,4), (5,5)\}$;}\\
1 & \mbox{otherwise.}
\end{cases}
$$
Thus it will suffice to prove that $\Delta(\alpha, \beta) \leq 3$
when $(p,q) = (2,3)$.

Suppose that $\Delta(\alpha, \beta) = 5$ and $(p, q)=(2,3)$. Lemma
\ref{charcount} and Inequality (\ref{6.6}) imply that $s_X = 2$.
Proposition \ref{betacurve} (1) shows that for every point $\chi_\rho
\in J_X(\alpha)$, $\rho(\pi_1(\p M)) = \mathbb Z/5$. In particular,
$\rho(\pi_1(M))$ has an element of order $5$.  The only finite,
non-cyclic subgroups of $PSL_2(\mathbb C)$ which have such elements
are $I_{60}$ and $D_k$ where $k \equiv 0 \mbox{ {\rm (mod $5$)}}$.
Therefore Lemma \ref{charcount} shows that $10 = 5s_X \leq s_X + 5 =
7$, which is impossible.

Suppose next that $\Delta(\alpha, \beta) = 4$. Lemma \ref{charcount}
and Inequality (\ref{6.6}) imply that $s_X \leq 3$. Let $\beta^*$ be a
dual slope to $\beta$ and recall that $\|\beta^*\|_X=s_X$.

If $s_X = 2$, then $8 = \Delta(\alpha, \beta) s_X = \|\alpha\|_X = 2 +
2n - A$ where $A \in \{0,1\}$.  Thus $n = 3, m = 0$ and so
$\nu(J_X(\alpha))$ consists of $3$ elements where at most two are
$I_{60}$-characters, at most one is an $O_{24}$-character, and at most
one is a $T_{12}$-character. Proposition \ref{betacurve}(1) shows
that for every point $\chi_\rho \in J_X(\alpha)$, $\rho(\pi_1(\p M)) =
\mathbb Z/2$ or $\mathbb Z/4$.  Since only the $O_{24}$-character has
elements of order $4$, there are at least two characters $\chi_\rho$ in
$J_X(\alpha)$ such that $\rho(\pi_1(\p M))=\mathbb Z/2$. This implies
that $4 = 2s_X = \|2\beta^*\|_X \geq s_X + 4 = 6$, which is
impossible.

Finally suppose that $s_X = 3$. Then $12 = \Delta(\alpha, \beta) s_X =
\|\alpha\|_X = 3 + 2n - A$ where $A \in \{0,1\}$.  Hence $n = 5$, and
$\nu(J_X(\alpha))$ consists of $5$ elements -- two
$I_{60}$-characters, one $O_{24}$-character, one $T_{12}$-character,
and one $D_3$-character. A similar argument to that of the previous
paragraph shows that $6 = 2s_X = \|2\beta^*\|_X \geq s_X + 7 = 10$,
which is impossible.  This completes the proof.  \end{proof}

\section{Proof of Proposition \ref{TLVS}} \label{PTLVS}\label{proof three}

In this section we suppose that $b_1(M) = 1$, neither $\alpha$ nor
$\beta$ is a singular slope for a closed essential surface in $M$,
$M(\alpha)$ is an irreducible, very small $3$-manifold which admits a
geometric decomposition, and $M(\beta)$ is either $S^1 \times S^2$ or
a connected sum of lens spaces $L_p \# L_q$ where $2 \leq p \leq q$.
We must show $\Delta(\alpha, \beta) \leq 2$.

The reader will verify that given our assumptions on $M(\alpha)$, 
one of the following possibilities holds. Either $M(\alpha)$
\begin{itemize}
\item is a torus bundle over $S^1$ with monodromy $A \in SL_2(\mathbb
Z)$ such that $|tr(A)| \geq 2$; or
\item semi-fibres over $I$ with semi-fibre a torus; or
\item admits a Seifert structure with base orbifold $S^2(3,3,3),
S^2(2,4,4)$,  or $S^2(2,3,6)$.
\end{itemize}
We treat these cases separately.
\setcounter{case}{0}
\begin{case}
  $M(\alpha)$ fibres over the circle with monodromy $A$ for which
  $|tr(A)| \geq 2$.
\end{case}

Note that $\alpha$ is the rational longitudinal class in this case so
that $M(\beta) \ne S^1 \times S^2$. Thus $M(\beta) \cong L_p \# L_q$
for some $2 \leq p \leq q$. According to Proposition \ref{TLVSLS} we
may assume that either $\Delta(\alpha, \beta) = 3$ and $(p,q) \in
\{(2,3), (2,5), (3,5)\}$ or $M(\beta) \cong P^3 \# P^3$. We consider
the former case first.

Let $X_0$ be a curve in $X_{PSL_2}(M(\beta)\subset X_{PSL_2}(M)$.
Since $X_0$ is a $\beta$-curve, it follows from Lemma
\ref{jumpnonideal} that, for each $x \in J_{X_0}(\alpha)$ and $\rho
\in R(X_0) \cap t^{-1}(\nu(x))$, we have that $\rho(\beta^*)$ has
order $3$.  Proposition \ref{bundlereps}(2) implies that there are no
reducible characters in $\nu(J_{X_0}(\alpha))$.  Hence if $\chi_\rho
\in \nu(J_{X_0}(\alpha))$, then the image of $\rho$ is either
contained in ${\mathcal N}$ or is $T_{12}$ (Proposition
\ref{toroidalreps}). Since $q>2$ it follows from (\ref{minformula}),
Lemma \ref{charcount}, and Proposition \ref{lb} that
$$\Delta(\alpha, \beta) \leq 
\begin{cases}
1 & \mbox{when $(p, q) \ne (2,3), (3,3)$} \\
2 & \mbox{when $(p, q) = (2,3), (3,3)$}
\end{cases} ,$$
contradicting our assumption that $\Delta(\alpha, \beta) =3$. 

Next suppose that $M(\beta) \cong P^3 \# P^3$. It follows that $H_1(M)
\cong \mathbb Z \oplus A$ where $A$ is either (i) $\mathbb Z/2$ or
(ii) $\mathbb Z/2 \oplus \mathbb Z/2$. Now $H_1(M(\alpha))$ is
infinite, so $\alpha$ is the slope of the rational longitude in
$H_1(\partial M)$, say $\alpha = \{\pm a\}$ and $i_*(a) = \sigma \in
A$ where $i: \partial M \to M$ is the inclusion.  If $\alpha^* = \{\pm
a^*\}$ is any dual slope to $\alpha$ we have $i_*(a^*) = d \xi +
\tau$ where $d \geq 1, \xi$ generates a free factor of $H_1(M)$ and
$\tau \in A$.  Write $\beta = \{\pm(ma + na^*)\}$ and observe that
$\Delta(\alpha, \beta) = |m|$.  A simple computation shows that since
$H_1(M(\beta)) \cong \mathbb Z/2 \oplus \mathbb Z/2$, we must have $m
= \pm 1$ in case (ii) and therefore $\Delta(\alpha, \beta) = 1$.
Similarly in case (i) we must have $\Delta(\alpha, \beta) \leq 2$.
Both cases contradict our hypotheses, so we also have $\Delta(\alpha,
\beta) \leq 2$ when $q = 2$.

\begin{case}
$M(\alpha)$ semi-fibres over the interval
\end{case}
        
\begin{subcase}
$M(\beta) \cong L_p \# L_q \ne P^3 \# P^3$
\end{subcase}
Again, according to Proposition \ref{TLVSLS}, we may assume that
either $\Delta(\alpha, \beta) = 3$ and $(p,q)\in\{(2,3), (2,5),
(3,5)\}$.

Let $\beta^*$ be a dual class to $\beta$.  According to Lemma
\ref{jumpnonideal}, for each $x \in J_{X_0}(\alpha)$ and $\rho \in R(X_0)
\cap t^{-1}(\nu(x))$, we have that $\rho(\beta^*)$ has order $3$.
Proposition \ref{semibundlereps} shows that there are no reducible
characters in $\nu(J_{X_0}(\alpha))$. Thus if $\chi_\rho \in
\nu(J_{X_0}(\alpha))$, the image of $\rho$ is either contained in
${\mathcal N}$ or is $O_{24}$ by Proposition \ref{toroidalreps}. Since
$q \geq 3$, Lemma \ref{charcount} and Proposition \ref{pssms} show that
$\Delta(\alpha, \beta) \leq 2$, contradicting our assumption that
$\Delta(\alpha, \beta) = 3$.

\begin{subcase}
$M(\beta) \cong P^3 \# P^3$
\end{subcase}

This case follows from Theorem 1.2 of [Lee]. 

\begin{subcase}
$M(\beta) = S^1 \times S^2$
\end{subcase}

There is an exact sequence $1 \to \mathbb Z \oplus \mathbb Z \to
\pi_1(M(\alpha)) \to \mathbb Z/2 * \mathbb Z/2 \to 1$ and therefore a
non-trivial curve $X_0 \subset X_{PSL_2}(\mathbb Z/2 * \mathbb Z/2)
\subset X_{PSL_2}(M(\alpha)) \subset X_{PSL_2}(M)$.  As we have
assumed that $\alpha$ is not a singular slope for a closed, essential
surface in $M$, Proposition \ref{4.10-4.12} implies that
$\|\cdot\|_{X_0} \ne 0$. Since we have assumed that $\beta$ is not a
singular slope for a closed, essential surface in $M$, the same
proposition implies that $J_{X_0}(\beta) \subset X_0^\nu$.  Thus
Proposition \ref{fin} (1) shows that $\Delta(\alpha, \beta) = 1$.

\begin{case}
$M(\alpha)$ admits a Seifert structure with base orbifold $S^2(3,3,3),
S^2(2,4,4)$, or $S^2(2,3,6)$
\end{case}

Our proof in this case depends on obtaining good estimates for the
value of a Culler-Shalen seminorm on $\alpha$. To that end, let $X_0
\subset X_{PSL_2}(M)$ be a non-trivial curve and suppose that $\chi$
is a character contained in $\nu(J_{X_0}(\alpha))$. Since $b_1(M) =
1$, $\chi = \chi_\rho$ where $\rho \in R(X_0)$ is either irreducible
or has a non-Abelian image by Lemma \ref{ordjump} and further,
$\rho(\alpha) = \pm I$.  Thus $\rho$ factors through
$\pi_1(M(\alpha))$. Now apply Lemma 3.1 of [BB] to see that $\rho$
factors through $\Delta(r,s,t)$, the orbifold fundamental group of the
base orbifold $S^2(r,s,t)$ of $M(\alpha)$ ($a \leq b \leq c$). The
irreducible characters $\Delta(r,s,t)$ were calculated in Propositions
5.2, 5.3, 5.4 of [B]. If $\chi_\rho$ is reducible, $\rho$ induces a
representation $\sigma: \Delta(r,s,t) \to PSL_2(\mathbb C)$ whose
image is upper-triangular and non-Abelian.  Write $\Delta(r,s,t) =
\langle x, y : x^a, y^b, (xy)^c \rangle$ and observe that up to
conjugation, $\sigma(x)$ is diagonal of order $a$ and the $(1,2)$
entry of $\sigma(y)$ is $1$. The reader will verify that as
$\sigma(xy)$ is of finite order, there is at most one possibility for
the character of $\sigma$. Thus we have proven,

\begin{lemma} \label{pqrchars}\ 
\begin{enumerate}
\item[(1)] $\Delta(3,3,3)$ has exactly one irreducible $PSL_2(\mathbb
  C)$-character and it is the character of a representation with image
  $T_{12}$. It has exactly one reducible $PSL_2(\mathbb C)$-character
  which can lie on a non-trivial curve in $X_{PSL_2}(M)$.
\item[(2)] $\Delta(2,4,4)$ has exactly three irreducible
  $PSL_2(\mathbb C)$-characters and they are the characters of
  representations with dihedral images $D_2, D_4$ and $D_4$. It has
  exactly one reducible $PSL_2(\mathbb C)$-character which can lie on
  a non-trivial curve in $X_{PSL_2}(M)$.
\item[(3)] $\Delta(2,3,6)$ has exactly two irreducible $PSL_2(\mathbb
  C)$-characters, one corresponding to a representation with image
  $D_3$, and the other to a representations with image $T_{12}$. It
  has exactly one reducible $PSL_2(\mathbb C)$-character which can lie
  on a non-trivial curve in $X_{PSL_2}(M)$.
\end{enumerate}
\end{lemma}
\qed

\noindent Proposition \ref{pssms} now yields the estimates we need. 

\begin{prop} \label{vsseif}
  Suppose that $X_0$ is a non-trivial curve in $X_{PSL_2}(M)$ and that
  $\alpha$ is a slope on $\partial M$ such that $M(\alpha)$ admits a
  Seifert structure with base orbifold $S^2(3,3,3)$, $S^2(2,4,4)$, or
  $S^2(2,3,6)$. If $\alpha$ is not a boundary slope associated to an
  ideal point of $X_0$, then
$$ \|\alpha\|_{X_0} \leq \left\{\begin{array}{ll} s_{X_0} + 4 & \hbox{ if } M(\alpha) \hbox{ has base orbifold } S^2(3,3,3) \\
s_{X_0} + 5 & \hbox{ if } M(\alpha) \hbox{ has base orbifold } S^2(2,3,6) \hbox{ or } S^2(2,4,4) 
\end{array} \right.$$ 
\end{prop} 
\qed

\begin{subcase}
$M(\beta) \cong L_p \# L_q$ where $2 \leq p \leq q$
\end{subcase}

Let $X_0 = X(1,1) \subset X_{PSL_2}(\mathbb Z/p * \mathbb Z/q) = 
X_{PSL_2}(M(\beta))
\subset X_{PSL_2}(M)$. Since we have assumed that $\beta$ is not a 
singular slope for any
closed essential surface in $M$, (\ref{minformula}) and Proposition \ref{lb} imply that
\begin{equation} \Delta(\alpha, \beta) = \frac{\|\alpha\|_{X_0}}{s_{X_0}} \;\;\;\;\; \mbox{ where } \;\;\;\;\;\;\;
s_{X_0} \geq \left\{ \begin{array}{ll}
2 & \mbox{if } p = 2 \\ 4 & \mbox{if } p > 2. \\
\end{array} \right.   
\end{equation}
Hence Proposition \ref{vsseif} yields $\Delta(\alpha, \beta) \leq 2$ when $p > 2$. Similarly, if there are no irreducible characters in $\nu(J_{X_0}(\alpha))$, then $\|\alpha\|_{X_0} \leq s_{X_0} + 2$ (cf. Lemma \ref{charcount}), which yields the desired distance estimate. Assume then that $p = 2$ and $\nu(J_{X_0}(\alpha))$ contains at least one irreducible character. 

\begin{subsubcase}
$2 = p = q$
\end{subsubcase}

In this case, all irreducible characters in $X_0$ are characters of
representations which conjugate into ${\mathcal N}$ and therefore the
base orbifold of $M(\alpha)$ cannot be $S^2(3,3,3)$ (Proposition
\ref{vsseif}). When it is $S^2(2,3,6)$, we obtain $\|\alpha\|_{X_0}
\leq s_{X_0} + 3$ and so $\Delta(\alpha, \beta) \leq 2$ by (10.2.1). When it is
$S^2(2,4,4)$, Corollary \ref{singormin} implies that the natural
homomorphism $\pi_1(M) \to \pi_1(M(\beta))$ sends $\pi_1(\partial M)$
to the unique index $2$ cyclic subgroup $C$ of $\mathbb Z/2 * \mathbb
Z/2$ (since $\beta$ is not a singular slope for a closed essential
surface in $M$). Thus $\pi_1(\partial M)$ is sent to $\pm I$ under the
diagonal representation whose character lies on $X_0$. It follows that
$\nu(J_{X_0}(\alpha))$ does not contain a reducible character (cf.
Proposition \ref{betacurve}). Thus Lemmas \ref{charcount} and
\ref{pqrchars} show that $\|\alpha\|_{X_0} \leq s_{X_0} + 3$, which
yields the desired result.

\begin{subsubcase}
$2 = p < q$
\end{subsubcase}

In this case, $X_0$ contains exactly one character of an irreducible
representation with image contained in ${\mathcal N}$ (Lemma
\ref{charcount}). Thus, when the base orbifold of $M(\alpha)$ is
$S^2(2,4,4)$ we have $\|\alpha\|_{X_0} \leq s_{X_0} + 3$ and therefore
$\Delta(\alpha, \beta) \leq 2$. When it is $S^2(3,3,3)$ we have
$\|\alpha\|_{X_0} \leq s_{X_0} + 4$ so that $\Delta(\alpha, \beta)
\leq 3$. If this distance is $3$, then $X_0$ contains the character of
a representation with image $T_{12}$ and therefore $q = 3$ (Lemma
\ref{charcount}). Then $H_1(M(\beta)) \cong \mathbb Z/6$ so that
$H_1(M) \cong \mathbb Z \oplus \mathbb Z/n$ where $n$ divides $6$.
There is a primitive element $\lambda \in H_1(\partial M)$, unique up
to sign, which is sent to a torsion element of $H_1(M)$. Let $d$ be
its order.  The argument used in the proof of Observation \ref{basis}
shows that there is a dual class $\mu \in H_1(\partial M)$ to
$\lambda$ which is sent to $(d, \bar j) \in \mathbb Z \oplus \mathbb
Z/n = H_1(M)$. If $\beta = a \mu + b \lambda$ in $H_1(\partial M)$,
then a homological calculation shows that $6 = |H_1(M(\beta))| =
|dan|$. As $d$ divides $n$ and $6$ is square-free, we have $d = 1$.
Hence $\lambda$ is homologically trivial in $M$ and therefore if
$\alpha = s\mu + t\lambda$, $H_1(M(\alpha)) \cong \mathbb Z/s \oplus
\mathbb Z/n$. Since this group surjects onto $H_1(\Delta(3,3,3)) \cong
\mathbb Z/3 \oplus \mathbb Z/3$, both $s$ and $n$ are divisible by
$3$. Thus $t$ is relatively prime to $3$ and the same holds for $a$ as
$6 = |H_1(M(\beta))| = |dan| = |an|$. Hence $\Delta(\alpha, \beta) =
|at - bs| \not \equiv 0 \mbox{ {\rm (mod 3)}}$, and we are done in
this case.

Finally assume that the base orbifold of $M(\alpha)$ is $S^2(2,3,6)$.
Since $\nu(J_{X_0}(\alpha))$ contains the character of an irreducible
representation, Lemmas \ref{pqrchars} and \ref{charcount} imply that
$q = 3$.  From Proposition \ref{vsseif} we have $\|\alpha\|_{X_0} \leq
s_{X_0} + 5$ so that $\Delta(\alpha, \beta) \leq 3$. Note moreover
that if $\Delta(\alpha, \beta) = 3$, then $s_{X_0} = 2$ and
$\nu(J_{X_0}(\alpha))$ consists of a $T_{12}$ character and a
reducible character (cf. Proposition \ref{pssms} and Lemma
\ref{vsseif}).  Now $H_1(M(\beta)) \cong \mathbb Z/6$, so $H_1(M)
\cong \mathbb Z \oplus \mathbb Z/n$ where $n$ divides $6$.  The
argument of the last paragraph shows that there is a basis $\mu,
\lambda$ for $H_1(\partial M)$ such that if $i: \partial M \to M$ is
the inclusion, then $i_*(\mu)$ generates a $\mathbb Z$-summand of
$H_1(M)$, while $i_*(\lambda) = 0$. Thus for a primitive class $\delta
= s \mu + t \lambda$ we have $H_1(M(\delta)) \cong \mathbb Z/s \oplus
\mathbb Z/n$. In particular taking $\beta = p \mu + q \lambda$ we have
$\mathbb Z/6 \cong \mathbb Z/p \oplus \mathbb Z/n$, so $\gcd(p,n) = 6$
and $pn = 6$.

There is a presentation $\pi_1(M(\alpha)) \cong \langle x, y, h \; | 
\; x^2 = h^{-i}, y^3 = h^{-j}, (xy)^6 =
h^{-k},  h \mbox{ {\rm central}} \rangle$ where $i,j,k$ are 
relatively prime to $2,3,6$ respectively.
Thus $H_1(M(\alpha))$ is presented by the matrix
$$A = \left( \begin{array}{ccc} 2 & 0 & 6 \\ 0 & 3 & 6 \\ i & j & k 
\end{array} \right).$$
Since the gcd of the minors of size $1$ of $A$ is $1$, as are those 
of size $2$, while the determinant
of $A$ is $6(k - 2j - 3i) \equiv 0 \mbox{ {\rm (mod 12)}}$, we have
$H_1(M(\alpha)) \cong \mathbb Z/12l$ where $l \geq 0$. On the other 
hand if $\alpha = s \mu + t \lambda$,
then $H_1(M(\alpha)) \cong \mathbb Z/s \oplus \mathbb Z/n$, so 
$\gcd(s,n) = 1$ and $sn = 12l$.
These two conditions are not mutually compatible when $n \in \{2, 
6\}$, so $n \in \{1, 3\}$ is odd. But then 
$3 = \Delta(\alpha, \beta) = |at - sb| \equiv 0 \mbox{ {(mod 2)}}$, 
which is impossible.
Hence we must have $\Delta(\alpha, \beta) \leq 2$.

\begin{subcase}
$M(\beta) = S^1 \times S^2$
\end{subcase}

In this case, $\beta$ is the slope of the rational longitude in
$H_1(\partial M)$ and therefore $b_1(M(\alpha)) = 0$. It follows that
$M(\alpha)$ is not Haken [Ja, VI.13] and therefore Theorem \ref{2.0.3}
implies that $\alpha$ is not a boundary slope. Note, moreover, that as
the Euler number $e(M(\alpha)) \in \mathbb Q$ is the obstruction to
the existence of a horizontal surface in $M(\alpha)$, and since a
Seifert manifold of the form we are considering admits a horizontal
surface if and only if its first Betti number is $1$ [Ja, VI.15], we
have $e(M(\alpha)) \ne 0$.

Consider the canonical curve $X_M \subset X_{PSL_2}(M)$ defined by a
complete hyperbolic structure [BZ5, \S9]. Denote by $B_M$ the largest
$\|\cdot\|_M$-ball which contains no non-zero elements of
$H_1(\partial M)$ in its interior and recall that $s_M$ is the radius
of $B_M$. We have assumed that $\beta$ is not a singular slope
associated to a closed, essential surface in $M$, and therefore
Theorem \ref{2.0.3} implies that $\beta$ is not a strict boundary
slope. It follows from Proposition \ref{fin} (1) that $\|\beta\|_M =
s_M$. Indeed, [CGLS, \S 1] implies that $Z_x(f_\beta) \leq
Z_x(f_\delta)$ for each $x \in \widetilde X_0$ and $\delta \in
H_1(\partial M)$. According to Proposition \ref{vsseif} we have
\begin{equation} \|\alpha\|_{X_M} \leq \left\{\begin{array}{ll} s_{M} + 4 & \hbox{ if } M(\alpha) \hbox{ has base orbifold } S^2(3,3,3) \\
s_{M} + 5 & \hbox{ if } M(\alpha) \hbox{ has base orbifold } S^2(2,3,6) \hbox{ or } S^2(2,4,4) 
\end{array} \right. 
\end{equation}

\begin{lemma} \label{height} Let $\beta=\{\pm b\}$ and $\beta^*=\{\pm b^*\}$
Then 
\begin{enumerate}
\item[(1)] $b \in \partial B_M$ but is not a vertex. No class of 
distance $2$ from $b$ lies on
$\partial B_M$.
\item[(2)] If $\pm(c_1 b + d_1 b^*), \pm(c_2 b + d_2 b^*), 
\ldots , \pm(c_k b + d_1 b^*) \in
H_1(\partial M)$ are the primitive classes associated to the vertices 
of $B_M$, then
$\sum_{i=1}^k |d_i| \leq s_M$.  
\item[(3)] If $s_M = 2$, then $\Delta(\alpha, \beta) \leq 
\frac{\|\alpha\|_M}{s_M}$. Further,  if
$s_M \geq 3$ then $\Delta(\alpha, \beta) < t\frac{\|\alpha\|_M}{s_M}$ where
$$t =\left\{
\begin{array}{ll}
\frac{6}{5} & \mbox{if } s_M = 3 \\
\frac{4}{3}   & \mbox{if } s_M = 4 \\
2  & \mbox{if } s_M \geq 5 \\
\end{array} \right.$$
\end{enumerate}
\end{lemma}

\begin{proof}
 As $\beta$ is not a strict boundary slope and $M(\beta)$ has a
cyclic fundamental group, $\|\beta\|_M = s_M$ (Proposition
\ref{fin}(1)) and $\beta$ is not a vertex of $B_M$. It is shown in
Lemma 6.4 of [BZ1] that if there is a class of distance $2$ from
$\beta$ lies on $\partial B_M$, then $\beta$ would be a vertex of
$B_M$.  This proves part (1).

It was shown in \S 1.4 of [CGLS] that there is a homomorphism $\phi_x:
H_1(\partial M) \to \mathbb Z$ such that $\Pi_x(f_\gamma) =
|\phi_x(\gamma)|$. Since $|\phi_x(\delta_x)| = 0$, it is simple to see
that for each $\gamma \in H_1(\partial M)$, $|\phi_x(\gamma)| = e
\Delta(\gamma, \delta_x)$ for some fixed integer $e \geq 1$. In
particular we have $e = \Pi_x(f_\beta)/\Delta(\beta, \delta_x)$. Hence
$d_x = \Delta(\beta, \delta_x)$ divides $\Pi_x(f_\beta)$ and for each
$\gamma$ we have $\Pi_x(f_\gamma) = \frac{\Delta(\gamma,
  \delta_x)}{\Delta(\beta, \delta_x)} \Pi_x(f_\beta)$.  Summing over
all the ideal points yields
\begin{equation}\|\gamma\|_M = \sum_{x} \frac{\Delta(\gamma, 
\delta_x)}{\Delta(\beta, \delta_x)} \Pi_x(f_\beta). 
\end{equation}
In particular, $s_M = \|\beta\|_M = \sum_{x} \Pi_x(f_\beta) \geq 
\sum_{x} |d_x|$. This proves part (2) of the
lemma.

It follows from part (1) that if $xb + yb^* \in B_M$, then
$|y| < 2$, and therefore part (3) of the lemma holds for $s_M \geq 5$.
Let
$$t_0 = \mbox{ {\rm sup} } \{ y \; | \; xb + yb^* \in B_M \}$$
and and
observe that$\alpha = \{\pm(pb + qb^*)\}$ where $q = \Delta(\alpha,
\beta)$. Since $\frac{s_M}{\|\alpha\|_M} \alpha \in B_M$ we have
$\Delta(\alpha, \beta) = q \leq t_0 \frac{ \|\alpha\|_M}{s_M}$.
Furthermore we have strict inequality if there is a unique $xb + yb^*
\in B_M$ with $y = t_0$, since in this case equality would imply that
$\alpha$ is the slope of a vertex of $B_M$ and therefore a strict
boundary slope.  To complete the proof of (3), we must show that $t_0$
is given as in the statement of the lemma when $s_M \in \{2,3,4\}$.

First note that there is a vertex of $B_M$ of the form $x_0b +
t_0b^*$.  Let $z$ be an ideal point of $X_M$ associated to a
strict boundary class $c_xb + d_xb^* = v_z \in
H_1(\partial M)$.  The vertex of $B_M$ associated to $v_z$ is
given by $\frac{s_M}{\|v_z\|_M} v_z$.  We explain below
how to calculate the maximum value taken on by the
$b^*$-coordinate of $\frac{s_M}{\|v_{z}\|_M} v_{z} \in
\partial B_M$ where $z$ varies over all ideal points of $X_M$.

If $X_M$ has $k$ ideal points $z_1, z_2, \ldots , z_k$, then
$\Pi_{z_1}(f_\beta), \Pi_{z_2}(f_\beta), \ldots , \Pi_{z_k}(f_\beta)$ 
gives a partition of
$s_M = \|\beta\|_M$ into $k$ positive integers. Let $v_{z_i} = 
c_ib + d_ib^*$ and recall that
$d_i$ divides $\Pi_{z_i}(f_\beta)$. If we have prior knowledge of the 
integers $\Pi_{z_i}(f_\beta), c_i, d_i$,
then we can calculate the values $\|v_{z_i}\|_M$ using (10.3.1), 
and therefore we can determine the vertices of
$B_M$. In general though, we are not given these values, so we 
proceed as follows.

Fix an integer $k \geq 2$, a partition $(\Pi_1, \Pi_2, \ldots ,
\Pi_k)$ of $s_M$, and a sequence of classes $v_i = c_ib + d_i
b^*$ where $d_i \geq 1$ is a divisor of $\Pi_i$. Set $\|v_i\| =
\sum_{j\ne i} \frac{\Delta(v_i, v_j)}{\Delta(b,
  v_j)} \Pi_j$ and $v_i = \frac{s_M}{\|v_i\|} v_i$.
Next we consider the polygon in $H_1(\partial M; \mathbb R)$ whose
vertices are $\pm v_1, \pm v_2, \ldots , \pm v_k$. We discard all
polygons which are not convex, or which contain a non-zero element of
$H_1(\partial M)$ in their interior, or whose maximal
$b^*$-coordinates are at least $2$, since such polygons cannot be the
boundaries of a possible $B_M$. In this way we obtain a list of the
possibilities for $B_M$ for each value of $s_M$. In particular we can
determine an upper bound for their maximal $b^*$-coordinate. For
instance when $s_M = 2$ or $3$, an $SL_2(\mathbb C)$ version of the
calculation is contained in Lemma 6.5 of [BZ1]. The case $s_M = 4$ is
handled similarly from this one observes that part (3) of the lemma
holds. This completes the proof of the lemma.  \end{proof}

Note that Inequality (10.2.2) and part (3) of the previous lemma show
that
\begin{equation}\Delta(\alpha, \beta) \leq 3 \end{equation}
We must show that this inequality is strict. Denote the base orbifold
of $M(\alpha)$ by $S^2(r,s,t)$.  By Observation \ref{basis}, there is
a dual slope $\beta^*=\{\pm b^*\}$ for $\beta=\{\pm b\}$, an integer $n
\geq 1$, and an isomorphism $H_1(M) \cong \mathbb Z \oplus \mathbb
Z/n$ such that if $i: \partial M \to M$ is the inclusion, then
$$i_*(b^*) = (n, \bar 0), \;\;\; i_*(b) = (0, \bar 1).$$
Let $\xi \in H_1(M)$ correspond to $(1, 0)$, so that $i_*(b^*) = n \xi$.
Choose integers $t, u$ such that
$$\alpha = \{\pm(t b^* + u b)\}.$$
Then $\Delta(\alpha, \beta) = |t|$.

\begin{lemma} \label{h1}
  There is an isomorphism $H_1(M(\alpha)) \cong \mathbb Z/(u,n) \oplus
  \mathbb Z/\frac{tn^2}{(u,n)}$ where $\mathbb Z/(u,n)$ and $\mathbb
  Z/\frac{tn^2}{(u,n)}$ are generated, respectively, by the images of
  $\frac{tn}{(u,n)} \xi + \frac{u}{(u,n)} i_*(\beta)$ and $\xi$.
  Furthermore,
\begin{enumerate}
\item[(1)] if $(r,s,t) = (3,3,3)$, then $(u,n) = 3$. Hence
  $\Delta(\alpha, \beta) = |t| \not \equiv 0 \mbox{ {\rm (mod 3)}}$.
\item[(2)] if $(r,s,t) = (2,4,4)$, then $(u,n) = 2$ and $n \equiv 0
  \mbox{ {\rm (mod 4)}}$. Hence
  $\Delta(\alpha, \beta) = |t|$ is odd.
\item[(3)] if $(r,s,t) = (2,3,6)$, then $\gcd(u,n) = 1$ and $tn^2$ is
  divisible by $12$.
\end{enumerate}
\end{lemma}

\begin{proof}
 Since $e(M(\alpha)) \ne 0$, $H_1(M(\alpha))$ is finite. Moreover
it follows from our conventions that it is presented by the matrix
$\left(\begin{array}{cc} tn & 0 \\ u & n\end{array} \right).$ Thus
$H_1(M(\alpha)) \cong \mathbb Z/(u,n) \oplus \mathbb
Z/\frac{tn^2}{(u,n)}$ where the factors are generated as claimed.
Comparison of this isomorphism with the calculations of the previous
lemma yields the remaining conclusions of this one.  \end{proof}

Part (1) of the previous lemma and Inequality (10.3.2) show that
$\Delta(\alpha, \beta) \leq 2$ when $(r,s,t) = (3,3,3)$. In order to
deal with the remaining two cases we suppose that $\Delta(\alpha,
\beta) = 3$ in order to derive a contradiction. Setting $\beta=\{\pm b\}$
and $\beta^*=\{\pm b^*\}$ we have,
$$\alpha = \pm(3 b^* + u b)\}$$
so that $\gcd(3,u) = 1$. 

Assume that $(r,s,t) = (2,3,6)$. Then Lemma \ref{h1} (3) implies that
$3n^2$ is even and $\gcd(u,n) = 1$, so $n$ is even and $u$ is odd.
Thus $\gcd(u,6) = 1$. Consider the representation $\rho: \pi_1(M) \to
PSL_2(\mathbb C)$ with image $D_3$ constructed as a composition of
surjective homomorphisms $\pi_1(M) \to \pi_1(M(\alpha)) \to
\pi_1^{orb}(S^2(2,3,6)) = \Delta(2,3,6) \to D_3 \subset PSL_2(\mathbb
C)$.  Now $\rho(\pi_1(\partial M)) \subset D_3$ is Abelian, hence
cyclic of order $1,2$ or $3$. It cannot have order $1$, as otherwise
it would factor through $\pi_1(M(\beta)) \cong \mathbb Z$. Thus it has
order $2$ or $3$.  Since $\rho(\alpha) = \pm I$ and $|u| =
\Delta(\alpha, \beta^*)$ is relatively prime to $6$,
$\rho(\pi_1(\partial M)))$ is generated by $\rho(b^*)$. Thus the
image of $b^*$ generates the image of $\pi_1(\partial M)$ under
the composition of $\rho$ with the Abelianization homomorphism $\phi:
D_3 \to \mathbb Z/2$. Now $\phi \circ \rho$ factors through $H_1(M)$
and $b^*$ is divisible by $2$ in this group (Observation
\ref{basis}).  Thus $\phi \circ \rho (b^*) = 0$ and therefore
$\rho(b^*) \in [D_3, D_3] = \mathbb Z/3$. But then as
$\Delta(\alpha, \beta) = 3$, $\rho$ factors through $\pi_1(M(\beta))
\cong \mathbb Z$, which is impossible.  We conclude that
$\Delta(\alpha, \beta) \leq 2$.

Finally assume that $(r,s,t) = (2,4,4)$. 
There is a dual class $\beta^*_0 = \{\pm b_0^*\}$ for $\beta$
such that
$$\alpha = \{\pm (b + 3b^*_0)\}.$$ 
Set
$$b^*_1 = b + b^*_0.$$
\begin{lemma} \label{244chars}\ 
\begin{enumerate}
\item[(1)] $b$ is sent to a generator of a $\mathbb Z/2$ factor of
  $H_1^{orb}(S^2(2,4,4)) \cong \mathbb Z/2 \oplus \mathbb Z/4$ under
  the composition $H_1(\partial M) \to H_1(M) \to
  H_1^{orb}(S^2(2,4,4))$.
\item[(2)] If $x \in J_{X_M}(\alpha)$, then $f_{\beta^*_1}(x) = 0$. 
\end{enumerate}
\end{lemma}

\begin{proof}
 (1) Lemma \ref{h1} implies that in our situation, $H_1(M(\alpha))
\cong \mathbb Z/2 \oplus \mathbb Z/ \frac{3n^2}{2}$ where $\mathbb
Z/2$ and $\mathbb Z/\frac{3n^2}{2}$ are generated, respectively, by
the images of $\omega = \frac{3n}{2} \xi + \frac{u}{(u,n)} i_*(b)$ and
$\xi$.  It follows that $\omega$ is sent to an element of order $2$ in
$H_1^{orb}(S^2(2,4,4)) \cong \mathbb Z/2 \oplus \mathbb Z/4$, and
$\xi$ is sent to an element of order $4$.  Lemma \ref{h1} also shows
that $\frac{u}{(u,n)}$ is odd, so the image of $b$ in
$H_1^{orb}(S^2(2,4,4))$ coincides with that of $\omega - \frac{3
  \epsilon n}{2} \xi$ for some $\epsilon = \pm 1$.  It follows that
$\xi$ and $b$ generate $H_1^{orb}(S^2(2,4,4))$, so the image of $b$ is
non-zero there, and since $n$ is divisible by $4$, the image of $2 b$
in $H_1^{orb}(S^2(2,4,4))$ is zero.  Thus (1) holds.

(2) Let $x \in J_{X_M}(\alpha)$ and set $\nu(x) = \chi_\rho$ where
$\rho \in R(X_M)$. As $\alpha = \{\pm(-2b + 3b^*_1)\}$, we see that
\begin{equation}\rho(b^*_1)^{-3} = \rho(b)^2. \end{equation}
We observed in the opening paragraph of \S 6.3 that $\rho$ factors
through a representation $\sigma: \Delta(2,4,4) \to PSL_2(\mathbb C)$.
If $\rho$ is reducible, there is a diagonal representation $\sigma_0:
\Delta(2,4,4) \to PSL_2(\mathbb C)$ with the same character as
$\sigma$. Since $\sigma_0$ factors through $H_1(\Delta(2,4,4)) \cong
\mathbb Z/2 \oplus \mathbb Z/4$, (10.5.1) shows that $\sigma_0$ sends
the image of $b_1^*$ to $\pm I$. It follows that (2) holds in this
case.

Assume next that $\rho$ is irreducible.  Lemma \ref{pqrchars} shows
that the image of $\rho$ is either $D_2$ or $D_4$ and so as
$\pi_1(M(\beta)) \cong \mathbb Z$, we have $\rho(b) \ne \pm I$. On
the other hand, by (10.5.1) it suffices to show that $\rho(b)^2
=\pm I$. This is obvious if $\chi_\rho$ is a $D_2$-character so
suppose that it is a $D_4$-character. Write $\Delta(2,4,4) = \langle
x,y \; | \; x^2 = y^4 = (xy)^4 = 1 \rangle$ and $D_4 = \langle z,w \;
| \; z^2 = w^4 = (zw)^2 = 1 \rangle \subset PSL_2(\mathbb C)$. There
are two characters of representations $\Delta(2,4,4) \to PSL_2(\mathbb
C)$ with image $D_4$ and they are represented by the homomorphisms
$\phi_1, \phi_2: \Delta(2,4,4) \to D_4$ where $\phi_1(x) = z,
\phi_1(y) = w$ and $\phi_2(x) = zw, \phi_2(y) = w$.  As these two
representations differ by an automorphism of $D_4$, it suffices to
prove that the image of $b$ in $D_4$ under $\phi_1$ has order $2$.

Suppose otherwise. Then its image has order $4$ and so $b$ is sent
to $y^\epsilon v \in \Delta(2,4,4)$ where $\epsilon \in \{\pm 1\}$ and
$v \in \mbox{ ker}(\phi_1)$. Now $\mbox{ ker}(\phi_1)$ is normally
generated in $\Delta(2,4,4)$ by $(xy)^2$, so $v = \Pi_{i=1}^k u_i
(xy)^{2\theta_i} u_i^{-1}$ where $u_i \in \Delta(2,4,4)$ and $\theta_i
\in \{\pm 1\}$. Now $x$, resp. $y$, projects to an element $\bar x$,
resp. $\bar y$, of order $2$, resp. $4$, in $H_1(\Delta(2,4,4))$ and
therefore as $b$ is sent to $(\epsilon + 2\sum_i \theta_i) \bar y$
in this group, it also has order $4$ there. But this contradicts part
(1) of the lemma. Therefore $b$ must be sent to an element of
order $2$ in $D_4$.  \end{proof}

Now we complete the proof of our current case.  We set $a = b+3b^*$,
so $\alpha = \{\pm a\}$.

Suppose first that $s_M = 2$. Then the only roots of $\tilde f_\beta$
on $\widetilde X_M$ are the two discrete, faithful characters of
$\pi_1(M)$. It follows that $\tilde f_\beta(x) \ne 0$ for each $x \in
J_{X_M}(\alpha)$. But then part (2) of the previous lemma shows that
$J_{X_M}(\alpha) \subset J_{X_M}(\beta_1^*)$ and so putting these
observations together with Lemma \ref{height} (3) we conclude that
$\|\beta_1^*\|_M \geq \|\alpha\|_M \geq 3 s_M$.  Then $\frac13
\|\beta_1^*\|_M \geq s_M$. Now $\frac17a$ lies on the line in
$H_1(\partial M; \mathbb R)$ which passes through $b$ and $\frac13
b_1^*$ and consideration of its position there shows that $\|\frac17
\alpha\|_M \geq s_M = 2$. But this contradicts $\|\alpha\|_M \leq s_M
+ 5 = 7$, so $s_M \ne 2$.

Next suppose that $s_M = 3$. If $x \in J_{X_M}(\alpha)$ is such that
$\chi_\rho = \nu(x)$ is irreducible, the image of $\rho$ is finite and
non-Abelian (Lemma \ref{pqrchars}), from which we deduce
$f_\beta(\chi_\rho) \ne 0$. Thus part (2) of the previous lemma shows
that $x \in J_{X_M}(\beta_1^*)$. It follows that $\|\beta_1^*\|_M \geq
s_M + 3 = 2 s_M$. Thus $[\frac12 b_1^*, b_0^*] \cap \hbox{int}
(B_M) = \emptyset$. But $\frac14 a \in [\frac12 b_1^*,
b_0^*]$ so that $\frac14 \|\alpha\| \geq s_M = 3$. But then $8 =
s_M + 5 \geq \|\alpha_M\| \geq 12$, which is impossible. Hence $s_M \ne
3$.

Next note that $s_M \ne 4$ since Lemma \ref{height} (3) shows that
$\Delta(\alpha, \beta) < \frac43 (\frac{9}{4}) = 3$.

Suppose then that $s_M \geq 5$ so that $\|\alpha\|_M \leq s_M + 5 \leq
2 s_M$, or equivalently, $\frac{a}{2} \in B_M$. The line segment
$[-b, \frac{a}{2}]$, which passes through $b_0^*$, is
contained in $B_M$. Therefore it is contained in $\partial B_M$ and
hence $\|\frac{a}{2}\|_M = s_M$. But then $2 s_M = \|\alpha\|_M
\leq s_M + 5$. It follows that $s_M = 5$. We noted above that $\alpha$
is not a boundary slope, so $\frac{a}{2}$ is not a vertex of
$B_M$, nor is $b + 2 b_0^*$ by Lemma \ref{height}(1). Thus
there is a vertex $v_0 = x_0 b + y_0 b_0^*$ of the edge of
$\partial B_M$ containing $[-b, \frac{a}{2}]$ with $2 <
\frac{y_0}{x_0} < 2$. Let $c_0 b + d_0 b_0^* \in H_1(\partial
M)$ be the boundary class which is a rational multiple of $v_0$. As
$s_M = 5$, part (2) of Lemma \ref{height} shows that $|d_0| \in \{3,
4\}$ and therefore since $\frac{d_0}{c_0} = \frac{y_0}{x_0}$, either
$\frac{3}{|c_0|} \in (2, 3)$ or $\frac{4}{|c_0|} \in (2, 3)$, which is
impossible. This final contradiction shows that $\Delta(\alpha, \beta)
\leq 2$ when $(r,s,t) = (2,4,4)$. (Lemma \ref{h1} (2) then shows that
we have $\Delta(\alpha, \beta) = 1$ in this case).

\section{Characteristic subsurfaces associated to
  a reducible Dehn filling} \label{CharacterSurfaces}

In this section we develop the background results needed to prove
Theorem \ref{TLS}.  We assume that $M$ is a simple manifold,
$M(\beta)$ is a connected sum of two non-trivial lens spaces one of
which is not $P^3$.

Recall that an embedded $2$-sphere in a 3-manifold is called essential
if it does not bound a $3$-ball. Since $M(\beta)$ is a connected sum
of two non-trivial lens spaces, a standard cut-paste argument shows
that there is an essential $2$-sphere $\hat F$ in $M(\beta)$ such that
$F=M\cap \hat F$ is a connected properly embedded essential planar
surface $F$ in $M$ with boundary slope $\beta$.  Any such surface $F$
is separating in $M$ since $M(\beta)$ has zero first Betti number.
Any such surface $F$ is not a semi-fibre since otherwise $M(\beta)$
would be a connected sum of two $P^3$'s. Among all such surfaces, we
assume that $F$ has been chosen to have the minimal number of boundary
components. Set $m=|\p F|$. Note that $m$ is an even number since $F$
is separating. Since $M$ is a simple manifold, we have $m\geq 4$. The
planar surface $F$ splits $M$ into two components, $X^+$ and $X^-$,
and $\hat F$ separates $M(\beta)$ as $\hat X^+$ and $\hat X^-$ each of
which is a punctured lens space. We may and shall assume that $\hat
X^+$ is not $P^3$. We use $\e$ to denote an element in $\{\pm\}$.

We call a properly embedded annulus $(A, \p A) \subset
(X^\epsilon, F)$ {\it essential} if its inclusion is not homotopic
rel $\partial A$ to a map whose image lies in $F$. The minimality
of $m = |\partial F|$ has the following useful consequence.

\begin{lemma}\label{equal} Suppose that $(A, \p A)\subset (X^\epsilon,
  F)$ is a properly embedded essential annulus. The boundary of $A$
  splits $\hat F$ into an annulus $B$ and two disks $N, N'$.  Then the
  number of boundary components of $F$ which lie in $N$ equals the
  number of boundary components of $F$ which lie in $N'$.
\end{lemma}

\begin{proof}
 Since $\hat X^\e$  has zero first Betti number, the annulus $A$
separates $\hat X^\epsilon$ into two pieces $W$ and $V$ where
$\partial W$ is a $2$-sphere and $\partial V$ is a torus.

Let $n, n'$ and $b$ be the number of boundary components of $F$
which lie in $N, N'$ and $B$ respectively. We may suppose that $n
\leq n'$. If $b = 0$, then $\partial V \subset M$ and so $V$ is a
solid torus in which the winding number of $B$ is at least $2$
(since $A$ is an essential annulus and thus not parallel to $B$).
It follows that a regular neighborhood in $M(\beta)$ of $N \cup V$
is a punctured lens space whose boundary $S$ is an essential
$2$-sphere in $M(\beta)$. Hence the number of components of $S
\cap \partial M$ is at least $m$. That is, $2n \geq m = n + n'$.
Hence $n \geq n'$, which implies that desired result.

On the other hand if $b > 0$, then $\partial W$ is inessential in
$M(\beta)$ and thus $W$ is a $3$-ball. If the $2$-sphere boundary
$S_1$ of a regular neighborhood $U$ in $X^\epsilon$ of $N \cup V$
is inessential in $M(\beta)$, it follows that $U$ is also a
$3$-ball. But this is impossible as it would imply that $\hat
X^\epsilon$ is a $3$-ball. Hence $S_1$ is essential in $M(\beta)$.
Since it intersects $\partial M$ in $2n + b$ components we have
$2n + b \geq m = n + b + n'$, i.e. $n \geq n'$. This completes the
proof. \end{proof}

Each essential annulus $A$ properly embedded in $(X^\epsilon, F)$
separates the punctured lens space $\hat X^\epsilon$, and hence
$X^\epsilon$. Let $V(A)$ be the component of $\hat X^\epsilon_A$
such that $V(A) \cap \hat F$ is an annulus $E(A)$. We call a pair
of disjoint essential annuli $A$ and $A'$ properly embedded in
$(X^\epsilon, F)$ and $(X^{\epsilon'}, F)$ {\it nested}
if  either $\partial A' \subset E(A)$ or $\partial A \subset
E(A')$.

The only Seifert fibred spaces contained in a simple manifold are
solid tori. This fact has the following useful application. 

\begin{lemma} \label{nested} If $A$ and $A'$ are disjoint essential annuli
  properly embedded in $(X^\epsilon, F)$ and $(X^{\epsilon'}, F)$,
  then they are nested.
\end{lemma}

\begin{proof}
 Let $c_0, c_1$ be the boundary components of $A$ and $c_0',
c_1'$ those of $A'$. We assume first that $\epsilon = \epsilon'$.
If $A, A'$ are not nested, then $V(A) \cap V(A') = \emptyset$ and
we can number the boundary components of $A$ and $A'$ in such a
way that they divide the $2$-sphere $\hat F$ into five components
whose interiors are pairwise disjoint: a disk $N$ bounded by
$c_1$; the annulus $B = E(A)$ bounded by $c_1$ and $c_0$; an annulus $E$
bounded by $c_0$ and $c'_0$; an annulus $B' = E(A')$ bounded by $c'_0$ and
$c'_1$; and a disk $N'$ bounded by $c'_1$. Let $n = |N \cap
\partial F|, b = |B \cap \partial F|$ and define $e, b', n'$
similarly. According to Lemma \ref{equal} we have $n = e + b' +
n'$ and $n' = n + b + e$. It follows that $b = e = b' = 0$ and
therefore $V(A), E, V(A') \subset M$. Since $M$ is simple, both
$V(A)$ and $V(A')$ are solid tori and as $A$ and $A'$ are
essential in $(X^\epsilon, F)$, the winding numbers of $B$ in
$V(A)$ and $B'$ in $V(A')$ are at least $2$ in absolute value. It
follows that a regular neighbourhood of $V(A) \cup E \cup V(A')$
in $M$ is Seifert fibred with an incompressible torus for
boundary. But the simple manifold $M$ does not contain such a
Seifert fibred space. Thus $A, A'$ must be nested.

Assume then that $\epsilon \ne \epsilon'$. The case where $E(A)
\cap E(A') = \emptyset$ can be shown to be impossible as in the
previous paragraph. Next suppose that $E(A) \cap E(A') \ne
\emptyset$ but neither $E(A) \subset E(A')$ nor $E(A') \subset
E(A)$. We number the boundary components of $A$ and $A'$ in such a
way that they divide the $2$-sphere $\hat F$ into five components
whose interiors are pairwise disjoint: a disk $N$ bounded by
$c_1$; an annulus $B$ bounded by $c_1$ and $c_0'$; an annulus $E$
bounded by $c_0'$ and $c_0$; an annulus $B'$ bounded by $c_0$ and
$c'_1$; and a disk $N'$ bounded by $c'_1$. Let $n = |N \cap
\partial F|$ and define $b, e, b', n'$ similarly. Lemma
\ref{equal} implies that $b = b' = 0$ and thus $A'$ may be
isotoped in $(X^{\epsilon'}, F)$ so that $\partial A' = \partial
A$. Then $T = A \cup A'$ is a torus in $M$ which must be
compressible as $m > 2$. As $T$ is not contained in a $3$-ball, it
bounds a solid torus $V$ in $M$. It is easy to see that $V = V(A)
\cup V(A')$ and so $E = V(A) \cap V(A') \subset F$, i.e. $e = 0$.
But then $E$ is isotopic through $V$ to either $A$ or $A'$, which
contradicts the essentiality of these two annuli. Hence it must be
that either $\partial A' \subset E(A)$ or $\partial A \subset
E(A')$ and thus $A, A'$ are nested. \end{proof}

\begin{lemma} \label{parallel}
  If $A$ and $A'$ are disjoint essential annuli properly embedded in
  $(X^\epsilon, F)$ and $(X^{\epsilon'}, F)$ such that a boundary
  component of $A$ is isotopic in $F$ to a boundary component of $A'$,
  then $\epsilon = \epsilon'$ and $A$ and $A'$ are parallel in
  $X^\epsilon$.
\end{lemma}

\begin{proof}
 By the previous lemma, $A$ and $A'$ are nested. Without loss
of generality we may suppose that $\partial A \subset E(A')$. Let
$c_0, c_1$ be the boundary components of $A$, and  $c'_0, c'_1$
those of $A'$, where the indices are chosen in such a way that the
four curves $c_0$, $c_1$, $c'_0$ and $c'_1$ divide $\hat F$ into
five components whose interiors are pairwise disjoint: a disk $N$
bounded by $c_0'$; an annulus $E \subset F$ bounded by $c_0$ and
$c'_0$; an annulus $B$ bounded by $c_0$ and $c_1$; an annulus $E'$
bounded by $c_1$ and $c'_1$; and a disk $N'$ bounded by $c'_1$.
Let $n$ be the number of components of $N \cap \partial M$. Define
$b, e', n'$ similarly so that $n + b + e' + n' = m$. Lemma
\ref{equal} shows that $e' = 0$. Now it must be that $\epsilon =
\epsilon'$ as otherwise $A'$ can be isotoped in $(X^{\epsilon'},
F)$ so that its boundary equals that of $A$. The argument in the
last paragraph of the proof of the previous lemma shows that this
situation cannot arise. Thus $\epsilon = \epsilon'$. If $b = 0$,
then $A' \cup E \cup B \cup E'$ is a torus bounding a solid torus
$V$ in $M$.  Since $A' \subset V$ and is not parallel into $F$, it
must be parallel into $A$. Thus the lemma holds. On the other
hand, if $b \ne 0$, then $S_1 = N \cup E \cup A \cup E' \cup N'$
is an inessential $2$-sphere in $M(\beta)$ and therefore bounds a
$3$-ball $W$ in $\hat X^{\epsilon}$. It follows that  $A$ and $A'$
are parallel in $X^\epsilon$ through $W$. This completes the
proof. \end{proof}

Let $(\Sigma^\e_*, \Phi^\e_*)\subset (X^\e,F)$ be the
characteristic Seifert pair of $(X^\epsilon, F)$ and $(\Sigma^\e,
\Phi^\e)\subset (X^\e,F)$ be the characteristic $I$-bundle pair it
contains. We shall use $\tau_\e$ to denote the free involution on
$\Phi^\e$ induced by $I$ fibres of $\Sigma^\e$. Let  $\Phi_j^\e$
denote the $j$-th characteristic subsurface with respect to the
pair $(M,F)$ as defined in \S 5 of [BCSZ]. Note that $\Phi_1^\e$
is the large part of $\Phi^\e$ and that the involution $\tau_\e$
restricts to a free involution on $\Phi_1^\e$, which will still be
denoted as $\tau_\e$. Let $(\Sigma_1^\e, \Phi_1^\e)$ be the
corresponding $I$-bundle pair.

\begin{lemma}\label{Phi+product} $(\Sigma^+, \Phi^+)$
  is a product $I$-bundle pair, i.e. there is no embedded M\"{o}bius
  band $(B,\p B)\subset (\Sigma^+, \Phi^+)$.  In particular,
  $\Phi^+\ne F$.
\end{lemma}

\begin{proof}
 Suppose otherwise that $(B,\p B)\subset (\Sigma^+, \Phi^+)$ is an
embedded M\"{o}bius band. Then $\p B$ bounds a disk $N$ in $\hat F$.
The union of $N$ and $B$ is an embedded projective plane in $\hat
X^+$. A regular neighborhood of this projective plane in $\hat X^+$ is
a punctured $P^3$. This implies that $\hat X^+$ itself is a punctured
$P^3$, contrary to our assumptions. \end{proof}

\begin{lemma} \label{NoMobAnn}
  Suppose that $(B,\partial B) \subset (X^-, F)$ is a properly
  embedded M\"{o}bius band. Then $\p B$ cannot be isotoped into
  $\Phi_1^+$.
\end{lemma}

\begin{proof}
 Let $A'$ be the essential annulus in $(X^-, F)$ which is the
frontier of a regular neighbourhood of $B$ in $X^-$. If $\partial B$
can be isotoped into $\Phi_1^+$ then the previous lemma shows that
there is an essential annulus $A$, properly embedded in $(X^+, F)$,
whose boundary contains $\partial B$. After a small isotopy of $A$ rel
$\partial B$ we can assume that $A$ and $A'$ are disjoint. But this
contradicts Lemma \ref{parallel} since a boundary component of $A$ is
isotopic to a boundary component of $A'$. Thus $\p B$ cannot be
isotoped into $\Phi_1^+$. \end{proof}

A {\it root torus} in $(X^\epsilon, F)$ is a solid torus $\Theta
\subset X^\epsilon$ such that $\Theta \cap F$ is an incompressible
annulus in $\p \Theta$ whose winding number in $\Theta$ is at least
$2$ in absolute value.

\begin{lemma} \label{notI}\ 
\begin{enumerate}
\item[(1)] Let $\Theta$ be a component of $\Sigma^\e_*$ and set $\Phi
  = \Theta \cap F$. If $(\Theta, \Phi)$ is not an $(I, \partial
  I)$-bundle, then $\Theta$ is a root torus
\item[(2)] Let $\Phi_1$ and $\Phi_2$ be distinct components of
  $\Phi^\e_*$ and $E \subset F$ an annulus whose boundary consists of
  a component $c_1$ of $\partial \Phi_1$ and a component $c_2$ of
  $\partial \Phi_2$. Then after possibly renumbering $\Phi_1, \Phi_2$,
  there are an annulus $E' \supseteq E$ in $F$ with $c_1 \subset
  \partial E'$ and components $\Sigma_1, \Sigma_2$ of
  $\Sigma_*^\epsilon$ such that $\Sigma_1$ is a product $I$-bundle
  component of $\Sigma_1^\epsilon$ containing $\Phi_1$ and $\Sigma_2$
  is a root torus such that $\Sigma_2 \cap F \subset E'$. Moreover,
  either
\begin{enumerate}
\item[(i)] $E = E'$, $\Sigma_1 \cap F = \Phi_1 \cup \Phi_2$,
  $\tau_\epsilon(c_1) = c_2$ and $\Sigma_2 \cap F \subset
  \mbox{int}(E)$, or
\item[(ii)] $E \ne E'$, $\Sigma_2 \cap F = \Phi_2 \subset
  \mbox{int}(E')$
\end{enumerate}
\end{enumerate}
In particular, there is a root torus in $X^\epsilon$ whose
intersection with $F$ lies in $E$.
\end{lemma}

\begin{proof}
 (1) Since simple manifolds contain no Seifert submanifolds
with incompressible boundaries, $\Theta$ is a solid torus. Now
$\Phi$ is a disjoint union of essential annuli $B_1, B_2, \ldots ,
B_n$. If $n > 1$ then Lemma \ref{parallel} shows that $n = 2$ and
$(\Theta, \Phi) \cong (S^1 \times I \times I, S^1 \times I \times
\partial I)$, contrary to our hypotheses. Thus $n = 1$ and from
the defining properties of the characteristic Seifert pair we see
that the winding number of $\Phi = B_1$ in $\Theta$  is at least
$2$ in absolute value.

(2) Let $\Sigma_1, \Sigma_2'$ be the components of
$\Sigma_*^\epsilon$ which contain $\Phi_1, \Phi_2$ respectively.
For $j = 1, 2$ there is a unique annulus $(A_j, \partial A_j)
\subset (\mbox{fr}_{X^\epsilon}(\Sigma^j), \partial \Phi^j)$ which
is essential in $(X^\epsilon, F)$ and which contains $c_j$. If
$A_1 = A_2$ then $\Sigma_1 = \Sigma_2'$ and so $\Sigma_1 \cap F
\supseteq \Phi_1 \cup \Phi_2$ has at least two components. It
follows that $\Sigma_1$ is a product $I$-bundle with $\Sigma_1
\cap F = \Phi_1 \cup \Phi_2$ (cf. part (1) of the lemma). Clearly
$\tau_\epsilon(c_1) = c_2$. Moreover $A_2 \cup E$ is a torus in
$M$ which bounds a solid torus $V \subset X^\epsilon$. Since $A_2$
is essential, it is isotopic to a component $\Sigma_2$ of
$\Sigma_*^\epsilon$ with $\Sigma_2 \cap F \subset \mbox{int}(E)$.
Thus (i) holds.

Assume then that $A_1 \ne A_2$.  According to Lemma
\ref{parallel}, $A_1$ and $A_2$ are parallel in $X^\epsilon$.
Hence there is another annulus $E^*$ in $F$ such that $\partial
(A_1 \cup A_2) = \partial (E \cup E^*)$. Lemma \ref{nested}
implies that at least one of the $(\Sigma_j, \Sigma_j \cap F)$,
say $(\Sigma_1, \Sigma_1 \cap F)$, is an $(I, \partial I)$-bundle.
Then $(\Sigma_2, \Sigma_2 \cap F)$ cannot be an $(I, \partial
I)$-bundle as otherwise the product region $N$ between $A_1$ and
$A_2$ could be used to build an $(I, \partial I)$-bundle structure
on $\Sigma_1 \cup N \cup \Sigma_2$, contrary to the defining
properties of $\Sigma_*^\epsilon$. Thus $(\Sigma_2, \Phi_2)$ is a
root torus. Set $E' = E \cup (\Sigma_2 \cap F) \cup E^*$ and
observe that (ii) holds. \end{proof}

A boundary component of $\Phi^\e$ or $\Phi_j^\e$ is called an
{\it inner} boundary component if it is not isotopic in $F$ to a
component of $\p F$,  otherwise it is called an {\it outer}
boundary component. Note that every boundary component $c$ of
$\Phi_1^\e$ is a boundary component of an essential annulus in
$(\Sigma_1^\e,\Phi_1^\e)\subset (X^\e,F)$ whose boundary is $c$
and $\tau_\e(c)$. The following result is a consequence of Lemma
\ref{equal}.

\begin{lemma}\label{BdryToBdry}
  A simple closed curve $c$ in $F$ is an inner, resp. outer, boundary
  component of $\Phi^\e$ if and only if $\tau_\e(c)$ is an inner, resp.
  outer, boundary component of $\Phi^\e$.
\end{lemma}
\qed

By Lemmas \ref{BdryToBdry} and \ref{parallel}, we can and shall
normalize $\Phi_j^\e$ to have the property that if a component of
$\p F$ is isotopic to a boundary component of $\Phi_j^\e$, then it
is already contained in $\Phi_j^\e$.

Recall from [BSCZ, Section 7] that a subsurface $T$ of $F$ is said
to be {\it tight} if the frontier of $T$ in $F$ is a connected
simple closed curve. Thus a  component of $\Phi_1^\e$ is tight if
and only if it has exactly one inner boundary component. It
follows from  Lemma \ref{BdryToBdry} that $\tau_\epsilon$ permutes the tight
components of $\Phi_1^\epsilon$. Note also that a component
$\Phi_0$ of $\Phi_1^\epsilon$ left invariant by the free
involution $\tau_\epsilon$ has an even number of inner boundary
components since $\tau_\epsilon|\Phi_0$ reverses orientation. In
particular, no tight component of $\Phi_1^\epsilon$ is invariant
under $\tau_\epsilon$. Thus they are paired by this involution.

\begin{lemma}\label{TwoTight}
  If $\Phi_1^\e\ne F$ and $\chi(F)=\chi(\Phi_1^\e)$, then $\Phi_1^\e$
  consists of a pair of tight components $T_1$, $T_2$ and it contains
  $\partial F$. Moreover $\tau_\e(T_1)=T_2$.
\end{lemma}

\begin{proof}
 Note that we also have $\chi(F)=\chi(\Phi^\e)$ and $\Phi^\e \ne
F$.  Obviously $\Phi^\e$ has at least two tight components $T_1, T_2$
with $\tau_\epsilon(T_1) = T_2$. If $c_j$ denotes the inner boundary
component of $T_j$, then we also have $\tau_\epsilon(c_1) = c_2$.
Since $\chi(F)=\chi(\Phi^\e)$, there is an annulus $E \subset
\mbox{int}(F)$ such that $E \cap \Phi^\epsilon = \partial E$ and $E
\cap T_1 = c_1$. According to Lemma \ref{notI} (2), there is a product
$I$-bundle component of $\Sigma^\epsilon$ which intersects $F$ in $T_1
\cup \tau_\epsilon(T_1) = T_1 \cup T_2$ and $\partial E = c_1 \cup
\tau_\epsilon(c_1) = c_1 \cup c_2$. It follows that $F = T_1
\cup_{c_1} E \cup_{c_2} T_2$ as claimed by the lemma.  \end{proof}

Suppose that $c$ is a simple closed curve in $F$.  We will say
that $c$ {\it sweeps out} an essential annulus in $(X^\epsilon,
F)$ if there is an essential annulus in $(X^\epsilon, F)$ having a
boundary component isotopic to $c$.

\begin{lemma}\label{tauann}
  Let $c$ be an essential simple closed curve contained in
  $\Phi^\epsilon$. If $c$ sweeps out an essential annulus $A$ in
  $(X^\epsilon, F)$, then $A$ is isotopic in $(X^\epsilon, F)$ to an
  essential annulus in the component of $\Sigma^\epsilon$ which
  contains $c$. In particular, $\partial A$ is isotopic in $F$ to $c
  \cup \tau_\epsilon(c)$.
\end{lemma}

\begin{proof}
 Let $\Phi_0$ be the component of $\Phi^\epsilon$ which
contains $c$ and $\Sigma_0$ the component of $\Sigma^\e$
containing $\Phi_0$. The annulus $(A, \partial A)$ is homotopic in
$(X^\epsilon, F)$ into a component $\Theta$ of the characteristic
Seifert pair $(\Sigma_*^\epsilon, \Phi_*^\epsilon)$. If $\Theta =
\Sigma_0$ then it is easy to see that  the lemma holds. On the
other hand if $\Theta \ne\Sigma_0$, then $c$ is isotopic in $F$ to
the core of an annulus $E \subset F$ whose boundary consists of a
component of $\partial \Phi_0$  and a component of $\partial
(\Theta \cap F)$. Without loss of generality we can suppose that
$c = \partial E \cap \Phi_0$. Then $c$ sweeps out an annulus $A_1
\subset \mbox{fr}_{X^\epsilon}( \Sigma_0)$ which is essential in
$(X^\epsilon, F)$. Set $c' = \partial E \setminus c \subset
\partial (\Theta \cap F)$ and let $A_2$ be the essential annulus
contained in $\mbox{fr}_{X^\epsilon}(\Theta)$ which is swept out
by $c'$. By Lemma \ref{parallel}, $A_1$ is parallel to $A_2$ in
$X^\epsilon$ and by Lemma \ref{notI}, $(\Theta, \Theta \cap F)$ is
a root torus. Since $A$ is homotopic into $\Theta$ but not into
$F$, it is isotopic to $A_2$, and therefore to $A_1 \subset
\Sigma^\epsilon$. This completes the proof of the lemma. \end{proof}

\begin{lemma}\label{SweepAnn}
Let $c$ be an essential
  simple closed curve in $F$. The following conditions are equivalent:
\begin{enumerate}
\item[(1)] $c$ sweeps out an essential annulus in $(X^\epsilon, F)$;
\item[(2)] $c$ is isotopic in $F$ to a simple closed curve $c'$ in
  $\Phi^\epsilon$ such that the geometric intersection number of $c'$
  and $\tau_\e(c')$ is 0.
\end{enumerate}
\end{lemma}

\begin{proof}
 From Lemma \ref{tauann} it is clear  that (1) implies (2).

If condition (2) holds for $c$ then, by choosing a negatively curved
metric on $ F$, we may assume that either $c'$ and $\tau_\e(c')$ are
disjoint or that $c'$ is invariant under $\tau_\epsilon$.  In the
first case, there is an essential annulus in $(X^\epsilon, F)$ with
boundary curves $c'$ and $\tau_\epsilon(c')$.  In the second case
there is an embedded M\"obius band $(B,\p B)\subset (X^\epsilon, F)$
with boundary curve $c'$. The frontier of a regular neighborhood of
$B$ in $X^\e$ is an essential annulus with both boundary curves
isotopic to $c'$, and hence to $c$. \end{proof}

\begin{lemma}\label{Boundary}
  If $c$ is an inner boundary component of $\Phi_j^\epsilon$ which is
  isotopic to a simple closed curve in $\Phi_{j+1}^{-\epsilon}$, then
  $c$ sweeps out an essential annulus in $(X^{-\epsilon}, F)$.
\end{lemma}

\begin{proof}
 Let $c'$ be a simple closed curve in $\Phi_{j+1}^{-\epsilon}$,
which is isotopic to $c$.  Since $\tau_{-\epsilon}(c')$ lies in
$\Phi_{j}^{\epsilon}$, and since $c'$ is isotopic to a boundary
curve of $\Phi_{j}^{\epsilon}$, it follows  that the geometric
intersection number of $c'$ and $\tau_{-\epsilon}(c')$ is zero.
Now apply Lemma \ref{SweepAnn}.\end{proof}

Recall from [BCSZ, Proposition 5.3.1] that for each $\e\in \{\pm\}$
and $j\geq 0$, there is a homeomorphism 
$h_j^\e:\Phi_j^\e\ra \Phi_j^{(-1)^{j+1}\e}$, unique up
to isotopy, which satisfies some useful properties. In particular
\begin{equation} h_{2j}^\e: \Phi_{2j}^\e \stackrel{\cong}{\lra} \Phi_{2j}^{-\e}
\mbox{ for each } \e\in \{\pm\}
   \mbox{ and each } j\geq 0. \end{equation}
Moreover 
$$ h_{2j+1}^\e: \Phi_{2j+1}^\e \stackrel{\cong}{\lra}
\Phi_{2j+1}^{\e} \mbox{ is a free involution for each } \e\in
\{\pm\} \mbox{ and each } j\geq 0. $$

For any compact surface $S$, $\chi(S)$ denotes the Euler
characteristic of $S$.

\begin{prop}\label{IndexTwo}
  Suppose that $j \geq 2$ and that $\chi(\Phi_j^\epsilon) =
  \chi(\Phi_{j+1}^\epsilon)$. Then $\Phi_j^\epsilon =
  \Phi_{j+1}^\epsilon$.
\end{prop}

\begin{proof}
 If $\Phi_j^\epsilon \ne \Phi_{j+1}^\epsilon$, there is an annulus
$(E, \partial E) \subset (\Phi_j^\epsilon \setminus {\rm
  int}(\Phi_{j+1}^{\epsilon}), \partial \Phi_{j+1}^\epsilon)$. We show
that this leads to a contradiction.

Consider the homeomorphism $h_j^\e:\Phi_j^\e\ra
\Phi_j^{(-1)^{j+1}\e}$. The image of $\Phi_{j+1}^\epsilon$ under this
homeomorphism is $\Phi_{1}^{(-1)^j\epsilon} \wedge
\Phi_{j}^{(-1)^{j+1}\epsilon}$ (Proposition 5.3.5 of [BCSZ]). Thus the
image $E_{0}$ of $E$ under this map satisfies $$E_0 \subset (F
\setminus {\rm int}(\Phi_1^{(-1)^j\epsilon})) \wedge
\Phi_{j}^{(-1)^{j+1}\epsilon}.$$
Let $c_0$ be a boundary component of
$E_0$. Then $c_0$ is a boundary component of $\Phi_1^{(-1)^j\epsilon}$
and thus is a boundary component of an annulus $A$ which is properly
embedded and essential in $(X^{(-1)^j\epsilon}, F)$. On the other
hand, $c_0$ is isotopic in $F$ to a curve in
$\Phi_{j}^{(-1)^{j+1}\epsilon}$ and so since $j \geq 2$, Lemma
\ref{Boundary} implies that $c_0$ is a boundary component of an
essential annulus $(A_1, \partial A_1) \subset
(X^{(-1)^{j+1}\epsilon}, F)$. But this contradicts Lemma
\ref{parallel}. Hence $\Phi_j^\epsilon = \Phi_{j+1}^\epsilon$ \end{proof}

\begin{cor} \label{JumpHigh}
  Fix $\e\in \{\pm1\}$ and suppose that $\chi(\Phi_{2k + 1}^\epsilon)
  < 0$ for some $k \geq 1$.  Then
$$\chi(\Phi_3^\epsilon) < \chi(\Phi_5^\epsilon) < \ldots <
\chi(\Phi_{2k + 3}^\epsilon).$$
\end{cor}

\begin{proof}
 Apply Proposition \ref{IndexTwo} and Proposition 5.3.9 of [BCSZ].
\end{proof}

\begin{lemma} \label{JumpLowNot}
  Suppose that $(X^\epsilon, F)$ is not a twisted $I$-bundle pair.
  Then $\chi(F) < \chi(\Phi_3^\epsilon)$.
\end{lemma}

\begin{proof}
 Suppose otherwise that $\chi(F) = \chi(\Phi_3^\epsilon)$.
According to the previous lemma we have $\Phi_2^\epsilon =
\Phi_3^\epsilon$ and therefore Proposition 5.3.9 of [BCSZ] implies
that $\Phi_1^\epsilon \ne \Phi_2^\epsilon$. But since
$\chi(\Phi_1^\epsilon) = \chi(\Phi_2^\epsilon)$, there is an annulus
$(E, \partial E) \subset (\Phi_1^\epsilon \setminus {\rm
  int}(\Phi_2^\epsilon), \partial \Phi_2^{\epsilon})$. Let $E_1 =
\tau_\e(E) \subset \Phi_1^\epsilon$ and observe that $E_1 \subset F
\setminus {\rm int}(\Phi_1^{-\epsilon})$ while $\partial E_1 \subset
\partial \Phi_1^{-\epsilon}$. By Lemma \ref{notI} there is a root
torus $V_1 \subset X^{-\epsilon}$ such that $V_1 \cap F \subset E_1$.
Let $A_1$ be the essential annulus in $(X^{-\epsilon}, F)$ given by
$\partial V_1 \setminus (V_1 \cap F)$.

Next observe that since $X^\epsilon$ is not a twisted $I$-bundle
but $\chi(F) = \chi(\Phi_1^\epsilon)$, there is an annulus $E_2
\subset  F \setminus {\rm int}(\Phi_1^{\epsilon})$ such that
$\partial E_2 \subset \partial \Phi_1^{\epsilon}$. Another
application of Lemma \ref{notI} produces a root torus $V_2 \subset
X^\epsilon$ such that $V_2 \cap F \subset E_2$. Let $A_2$ be the
essential annulus in $(X^{-\epsilon}, F)$ given by $\partial V_2
\setminus (V_2 \cap F)$. Since $V_1 \cap F \subset E_1 \subset
\Phi_1^\epsilon$ and $V_2 \cap F \subset E_2 \subset F \setminus
{\rm int}(\Phi_1^{\epsilon})$, we may suppose that $V_1 \cap V_2$
is empty. But then $A_1, A_2$ are disjoint essential annuli which
are not nested, contrary to Lemma \ref{nested}. Thus we must have
$\chi(F) < \chi(\Phi_3^\epsilon)$. \end{proof}

\begin{lemma} \label{Length}
  Any reduced homotopy in $(M, F)$ has length at most $m-2$ if $(X^-,
  F)$ is not a twisted $I$-bundle pair or has length at most $m-1$ if
  $(X^-, F)$ is a twisted $I$-bundle pair. Furthermore if a reduced
  homotopy in $(M, F)$ has length $m-1$ then it starts and ends on the
  $X^-$ side.
\end{lemma}

\begin{proof}
  First note that by ([BCSZ, Corollary 5.3.8]),
$\chi(\Phi_{j}^\epsilon)$ is even for each $j\geq 1$ odd. Applying
this together with  Corollary \ref{JumpHigh} and Lemma
\ref{JumpLowNot} we see that  $\Phi_{m-1}^\e$ is the empty set if
$(X^\e,F)$ is not a twisted $I$-bundle pair. So the length of a
reduced homotopy in $(M,F)$ is at most $- \chi(F)=m-2$ if the
homotopy starts on a side which is not a twisted $I$-bundle pair,
and at most $ 1 - \chi(F)=m-1$ if the homotopy starts on a side
which is a twisted $I$-bundle. In the latter case, the homotopy starts 
on the $X^-$ side by Lemma \ref{Phi+product}, and finishes there since $m$ is even.
   \end{proof}

It follows from the definition of $\Phi_j^\e$ that if $(X^-, F)$
is a twisted $I$-bundle pair, then $\Phi_{2j}^-=\Phi_{2j+1}^-$ and
$\Phi_{2j+1}^+=\Phi_{2j+2}^+$ for each $j\geq 0$.

\begin{lemma} \label{OneThree}If  $(X^-, F)$ is a twisted $I$-bundle pair
  and $\Phi_1^+$ is not empty, then $\chi(\Phi_1^+)<\chi(\Phi_3^+)$.
\end{lemma}

\begin{proof}
 Suppose otherwise. Then $\chi(\Phi_2^+)=\chi(\Phi_3^+)$. By
Proposition \ref{IndexTwo}, we have $\Phi_2^+=\Phi_3^+$. Thus
$\Phi_1^+=\Phi_2^+=\Phi_3^+$. But this is impossible as it
contradicts Proposition 5.3.9 of [BCSZ]. \end{proof}

\begin{prop} \label{Length3}
  If $\chi(F)<\chi(\Phi_1^+)$, then any reduced homotopy in $(M, F)$
  has length at most $m-3$.
\end{prop}

\begin{proof}
 First assume that $(X^-, F)$ is a twisted $I$-bundle pair. It
follows from  Lemma \ref{OneThree}, Corollary \ref{JumpHigh} and
the assumption $\chi(F)<\chi(\Phi_1^+)$ that $\Phi_{m-3}^+$ is the
empty set. Hence a reduced homotopy in $(M, F)$ has length  at
most $m-4$ if it starts on $X^+$ side and length at most $m-3$ if
it starts on $X^-$ side.

Suppose then that $(X^-,F)$ is not a twisted $I$-bundle pair. If
$\chi(\Phi_1^+)<\chi(\Phi_3^+)$, then arguing as in the previous
paragraph  yields the desired conclusion. Suppose then that
$\chi(\Phi_1^+)=\chi(\Phi_3^+)$. Then $\Phi_2^+=\Phi_3^+$ by
Proposition \ref{IndexTwo}. It follows from the definition of the
characteristic subsurfaces that $\Phi_{2j+1}^-=\Phi_{2j+2}^-$ and
$\Phi_{2j}^+ =\Phi_{2j+1}^+$ for each $j\geq 1$. Now Corollary
\ref{JumpHigh} and the condition $\chi(F)<\chi(\Phi_1^+)$ imply
that  $\Phi_{m-1}^+$ is the empty set. Since $m-1$ is an odd
number, $\Phi_{m-2}^+=\Phi_{m-1}^+$ is the empty set. But
$\Phi_{m-2}^-$ is homeomorphic to $\Phi_{m-2}^+$ (cf. (11.11.1)) and
thus $\Phi_{m-3}^-=\Phi_{m-2}^-$ is the empty set.
   Therefore the length of a
reduced homotopy in $(M,F)$ is at most $m-4$ if the homotopy
starts on the $X^-$ side and therefore at most $m-3$ in general.
\end{proof}

\begin{cor}\label{Length4} If
  there is a reduced homotopy in $(M, F)$ with length at least $m-2$,
  then $\Phi_1^+$ consists of a pair of tight components and contains
  $\p F$. Further, $\Phi_1^-$ is either a twisted $I$-bundle or
  consists of a pair of tight components and contains $\p F$.
\end{cor}

\begin{proof}
 By Proposition \ref{Length3}, we have $\chi(F)=\chi(\Phi_1^+)$. By
Lemma \ref{Phi+product} $\Phi_1^+\ne F$.  Now apply Lemma
\ref{TwoTight} to see that $\Phi_1^+$ consists of a pair of tight
components and contains $\p F$.

If $\Phi_1^-$ is not a twisted $I$-bundle, we may exchange $X_+$ and $X_-$
\end{proof}

\section{Proof of Proposition \ref{TLS}}\label{PTLS}\label{proof four}

Recall that we are assuming that $\beta$ is a strict boundary slope,
$M(\beta)$ is a connected sum of two non-trivial lens spaces, one of
which is not $P^3$, and $M(\alpha)$ admits a $\pi_1$-injective
immersion of a torus. We will use the method of [BCSZ] to show that
$\Delta(\alpha, \beta) \leq 4$.

Let $V_\alpha$ be the filling solid torus used in forming $M(\alpha)$. As
in [BCSZ] we obtain a map $h: T\ra M(\alpha)$ from a torus $T$ to
$M(\alpha)$ such that
\begin{enumerate}
\item[(1)] $h^{-1}(V_\alpha)$ is a non-empty set of embedded disks in
  $T$ and $h$ is an embedding when restricted on $h^{-1}(V_\alpha)$;
\item[(2)] $h^{-1}(F)$ is a set of arcs or circles properly embedded
  in the punctured torus $Q=T\setminus h^{-1}(V_\alpha)$, where $F$ is
  the planar surface given in Section \ref{CharacterSurfaces};
\item[(3)] If $e$ is an arc component of $h^{-1}(F)$, then $h:e\ra F$
  is an essential (immersed) arc;
\item[(4)] If $c$ is a circle component of $h^{-1}(F)$, then $c$ does
  not bound a disk in $Q$ and $h:c\ra F$ is an essential (immersed)
  $1$-sphere.
\end{enumerate}

For any subset $s$ of $T$, we use $s^*$ denote its image under the map
$h$.  Denote the components of $\p(h^{-1}( V_\alpha))$ by $a_1, ...,
a_{n}$ so that $a_1^*, ..., a_n^*$ appear consecutively on $\p M$.
Note again that $a_1, ..., a_n$ are embedded in $\p M$ and each of
these curves has slope $\alpha$. Denote the components of $\p F$ by
$b_1, ..., b_{m}$ so that they appear consecutively in $\p M$.  We fix
an orientation on $Q$ and let each component $a_i$ of $\p Q$ have the
induced orientation. Two components $a_i$ and $a_j$ are said to have
the same orientation if $a_i^*$ and $a_j^*$ are homologous in $\p M$.
Otherwise, they are said to have different orientations. Similar
definitions are defined for the components of $\p F$. Since $Q$, $F$
and $M$ are all orientable, one has the following

\noindent {\bf Parity rule}: An arc component $e$ of $h^{-1}(F)$
in $Q$ connects components of $\p Q$ with the same orientation
(resp. opposite orientations) if and only if the corresponding
$e^*$ in $F$ connects components of $\p F$ with opposite
orientations (resp. the same orientation).

We define a graph $\G$ on the torus $T$ by taking $h^{-1}(V_\alpha)$
as (fat) vertices and taking arc components of $h^{-1}(F)$ as edges.
Note that $\G$ has no trivial loops, i.e.  no $1$-edge disk faces.
Also note that each $a_i^*$ intersects each component $b_j$ in $\p M$
in exactly $\D(\alpha,\beta)$ points. If $e$ is an edge in $\G$ with
an endpoint at the vertex $a_i$, then the corresponding endpoint of
$e^*$, is in $a_i^*\cap b_j$ for some $b_j$, and the endpoint of $e$
is thus given the label $j$. So when we travel around $a_i$ in some
direction, we see the labels of the endpoints of edges appearing in
the order $1,...,m,..., 1,...,m$ (repeated $\D(\alpha,\beta)$ times).
It also follows that each vertex of $\G$ has valence
$m\D(\alpha,\beta)$.

Suppose that $e$ and $e'$ are two adjacent parallel edges of $\G$.
Let $R$ be the bigon face between them, realizing the parallelism.
Then $(R, e\cup e')$ is mapped into $(X^\e, F)$ by the map $h$ for
some $\e$.  Moreover $h|_R$ provides a basic essential homotopy
between the essential paths $h|_e$ and $h|_{e'}$ (cf. [BCSZ]). We may
and shall assume that $R^*=h(R)$ is contained in the characteristic
$I$-bundle pair $(\Sigma_1^\e, \Phi_1^\e)$ of $(X^\e, F)$.  We may
consider $R$ as $e\times I$ and assume that the map $h: R\ra
\Sigma_1^\e$ is $I$-fibre preserving.

A face $f$ of $\G$ is said to lie on the $X^\e$ side if $f^*$ is
contained in $X^\e$. Every face of $\G$ lies on either the $X^+$ side
or the $X^-$ side. Since $F$ separates $M$, if two faces of $\G$ share
a common edge, then the two faces will lie on different sides of $F$.

The torus $\p M$ is cut by $\p F$ into $m=2g$ parallel annuli. We
denote these annuli by $B_1,....,B_m$ so that $\p B_i=b_i\cup b_{i+1}$
for $i=1,...,m-1$ and $\p B_m=b_m\cup b_{1}$. We may assume that $B_1$
is contained in $X^-$. Then for each odd $i$, $B_i$ is contained in
$X^-$ and for each even $i$, $B_i$ is contained in $X^+$. So $\p
X^-=F\cup B_1\cup B_3\cup...\cup B_{2g-1}$ and $\p X^+=F\cup B_2\cup
B_4\cup...\cup B_{2g}$, both being closed surface of genus $g$.

The complement of the interior of $F$ in the essential 2-sphere
$\hat F$ is a set of $m$ disjoint meridian disks of the attached
solid torus $V_\alpha$. These disks cut the solid torus $V_\alpha$ into
$m$ pieces, denoted $H_1,...,H_m$, such that each $H_i$ is a
2-handle attached to $X^-$ (when $i$ odd) or to $X^+$ (when $i$
even) along $B_i$.

Suppose that  the characteristic $I$-bundle pair $(\Sigma_1^+,
\Phi_1^+)\subset (X^+,F)$ is a connected trivial $I$-bundle
containing all $B_i$ with $i$ even, i.e. $\Phi_1^+$ is a pair of
tight components $T_1$ and $T_2$ including all components of $\p
F$. This happens when the length of a reduced homotopy in $(M,F)$
is at least $m-2$ by Corollary \ref{Length4}. Let $\hat
\Sigma_1^+$ be $\Sigma_1^+$ with all the $2$-handles $H_i$, $i$
even, attached along $B_i$. Then $\hat \Sigma_1^+$ is an
$I$-bundle over the disk $\hat T_1$, where $\hat T_1$ is the disk
in $\hat F$ whose intersection with $F$ is the tight component
$T_1$ of $\Phi_1^+$. Write $\hat \Sigma_1^+=\hat T_1\times [0,1]$.
Let $D_{1/2}=\hat T_1\times\{1/2\}$. Let $U$ be the union of
$\hat\Sigma_1^+$ and a regular neighborhood of $\hat F$ in $\hat
X^+$. Obviously $U$ is a once punctured solid torus with $D_{1/2}$
as a meridian disk. The torus boundary of $U$
must bound a solid torus $V$ in $X^+$. That is, the once punctured
lens space $\hat X^+$ is the union of $U$ and $V$ along their torus
boundary. Hence the core curve of $U$ carries a generator
of the first homology group of $\hat X^+$ when given an
orientation. We record this property in the following lemma which
will be used later in the proof of Lemma \ref{Corner}.

\begin{lemma}\label{homology}If the length
  of a reduced homotopy in $(M,F)$ is at least $m-2$, then the core
  curve of the punctured solid torus $U$ given in the proceeding
  paragraph carries a generator of the first homology group of the
  non-trivial punctured lens space $\hat X^+$ and the disk $D_{1/2}$
  is a meridian disk of $U$.
\end{lemma}

\begin{defn}
  A pair of adjacent parallel edges $\{e,e'\}$ of $\G$ is called an
  $S$-cycle if
\begin{itemize}
\item the two edges connect two vertices $v$ and $v'$ with the same
  orientation;
\item the label of the endpoint of $e$ at $v$ is $j$ and the label of
  the endpoint of $e$ at $v'$ is $j+1$ (note all calculations
  concerning labels are defined mod (m));
\item the label of the endpoint of $e'$ at $v$ is $j+1$ and the label
  of the endpoint of $e'$ at $v'$ is $j$.
\end{itemize}
An $S$-cycle $\{e,e'\}$ is called an extended $S$-cycle if the
two edges $e$ and $e'$ are the two middle edges in a family of four
adjacent parallel edges of $\G$. (Figure 1).
\end{defn}

\begin{figure}
\centerline{\includegraphics[width=5truein]{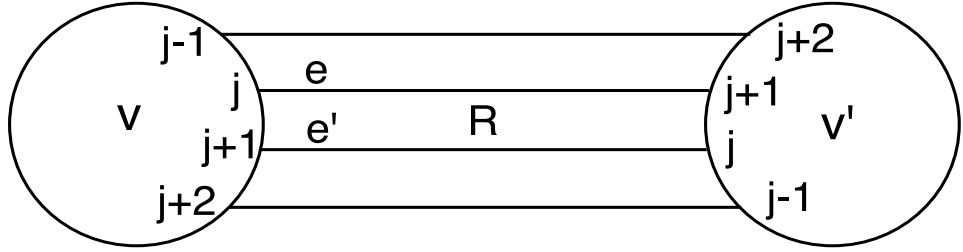}}
\caption{An extended $S$-cycle }
\end{figure}

\begin{lemma}\label{Scycle} Suppose that two vertices $v$ and $v'$
   of $\G$ have the same orientation and
   are connected by a family of $n$ parallel consecutive edges
   $e_1,...,e_n$ of $\G$.
   \begin{enumerate}
   \item[(1)] If $n>m/2$, then there is an $S$-cycle in this family of
     edges.
   \item[(2)] If $n>\frac{m}{2}+1$, then either there is an extended
     $S$-cycle in this family of edges or both $\{e_1,e_2\}$ and
     $\{e_{n-1},e_n\}$ are $S$-cycles.
   \item[(3)] If $n>\frac{m}{2}+2$, then there is an extended
     $S$-cycle in this family of edges.
   \end{enumerate}
\end{lemma}

\begin{proof}
 Part (1) is [CGLS, Corollary 2.6.7]. Parts (2) and (3) follow from
part (1) directly.
\end{proof}

\begin{lemma}\label{No+Scycle} If $\{e,e'\}$ is an $S$-cycle in $\G$,
  then the bigon face $R$ between them is mapped into $(\Sigma_1^-,
  \Phi_1^-)$ under the map $h$. Moreover there is properly embedded
  M\"obius band $B\subset X^-$ such that $\p B$ is contained in
  $\Phi_1^-$.
\end{lemma}

\begin{proof}
 Assume $R^*$ is contained in $\Sigma_1^\e$ and that the $S$-cycle
has labels $j$ and $j+1$.  Then $e^*$ and $e'^*$ are paths in $F$
connecting the two components $b_j$ and $b_{j+1}$ of $\p F$. Recall
that $B_j$ denotes the annulus in $\p M$ with boundary $b_j\cup
b_{j+1}$, and $\tau_\e$ denotes the involution of $\Phi_1^\e$.  We
have $\tau_\e(b_j)=\tau_\e(b_{j+1})$ and $\tau_\e(e^*)=\tau_\e(e'^*)$
and hence the connected set $b_j\cup e^*\cup e^{'*}\cup b_{j+1}$ is
invariant under $\tau_\e$. There is a $\tau_\e$-invariant regular
neighborhood $N$ of $b_j\cup e^*\cup e^{'*}\cup b_{j+1}$ contained in
$\Phi_1^\e$ and it is simple to see that there is a
$\tau_\e$-invariant essential simple closed curve in $N$.  Thus there
is a properly embedded M\"obius band $B\subset X^\e$ such that $\p B$
is contained in $N$.  Therefore by Lemma \ref{Phi+product}, we have
$\e=-$. \end{proof}

\begin{lemma}\label{NoEScycle}
  There is no extended $S$-cycle in $\G$.
\end{lemma}

\begin{proof}
 Suppose that $\{e, e'\}$ is an extended $S$-cycle of $\G$ as shown
in Figure 1. If $R$ denotes the bigon face between $e$ and $e'$, Lemma
\ref{No+Scycle} shows that $R^*$ is contained in $\Sigma_1^-$.  From
Figure 1, one easily sees that the set $b_j\cup e^*\cup e'^*\cup
b_{j+1}$ is contained in $\Phi_1^+$ and so the same may be assumed
true for its regular neighborhood $N$ used in the proof of Lemma
\ref{No+Scycle}. It follows that the boundary of the M\"obius band in
$X^-$ constructed in the proof of Lemma \ref{No+Scycle} is contained
in $\Phi_1^+$. But Lemma \ref{NoMobAnn} prohibits this possibility.
Thus $\G$ contains no extended $S$-cycles.  \end{proof}

\begin{lemma}\label{NoSameOri} Suppose that $m\geq 6$.
   If two vertices $v$ and $v'$
   of $\G$ have the same orientation, then they cannot be
   connected by $ 5m/6$ parallel edges.
\end{lemma}

\begin{proof}
 By Lemma \ref{NoEScycle} and Lemma \ref{Scycle} (3), $5m/6\leq
\frac{m}{2}+2$, i.e. $m\leq 6$.  So suppose $m=6$ and there are
$5m/6=5=\frac{m}{2}+2$ parallel consecutive edges $e_1,...,e_5$
connecting two vertices with the same orientation. Then by Lemma
\ref{Scycle} (2) and Lemma \ref{NoEScycle}, we may assume that both
$\{e_1,e_2\}$ and $\{e_4,e_5\}$ are $S$-cycles.  But the bigon face
$R$ between $e_1,e_2$ and the bigon face $R'$ between $e_4,e_5$ are on
different sides of $F$. This contradicts Lemma \ref{No+Scycle}. \end{proof}

Suppose that $\G$ has $m$ consecutively parallel edges $e_1,...,
e_{m}$ connecting two vertices $v$ and $v'$ with different
orientations. The existence of the $m$ parallel edges implies that
there is a length $m-1$ reduced homotopy in $(M,F)$.  Let $R_i$ denote
the bigon face between the adjacent parallel edges $e_i$ and
$e_{i+1}$, $i=1,...,m-1$. Then $R_1^*,...,R_{m-1}^*$ are contained
alternatively in $X^-$ and $X^+$ starting and ending on the $X^-$ side
of $F$ by Lemma \ref{Length}. Thus each of the bigon faces $R_1$,
$R_3$, ..., $R_{m-1}$ is mapped in $(\Sigma_1^-, \Phi_1^-)$ and each
of $R_2$, $R_4$, ..., $R_{m-2}$ is mapped in $(\Sigma_1^+, \Phi_1^+)$.

Orient all the edges $e_1, ..., e_{m}$ in the same direction such
that their tails are in $v$ and their heads are in $v'$. Up to
renumbering, we may assume that the labels of the tails of
$e_1,..., e_{m}$ are $1,...,m$ respectively. The labels of the
heads of $e_1, ...,e_{m}$ are $\s(1), ..., \s(m)$ for some
permutation $\s$ of $\{1,...,m\}$. (Note that the indices are
defined modulo $m$.)

Since $F$ separates $M$, $b_i$ and $b_{i+1}$ have different
orientations, for all $i$.  Also $b_i$ and $b_j$ have the same
orientation if and only if $i\equiv j$ (mod $2$).  By the parity rule,
for each $i\in \{1,...,m\}$, the components $b_i$ and $b_{\s(i)}$ of
$\p F$, connected by $e_i^*$, have the same orientation.  (Note that
if $b_i$ and $b_{\s(i)}$ are the same component of $\p F$ for some
$i$, then they are the same component for all $i=1,...,m$, i.e. $\s$
is the trivial permutation.)  It follows that $b_i$ is different from
$b_{\s(i+1)}$ and that $b_{\s(i)}$ is different from $b_{i+1}$, for
all $i$.

\begin{figure}
\centerline{\includegraphics{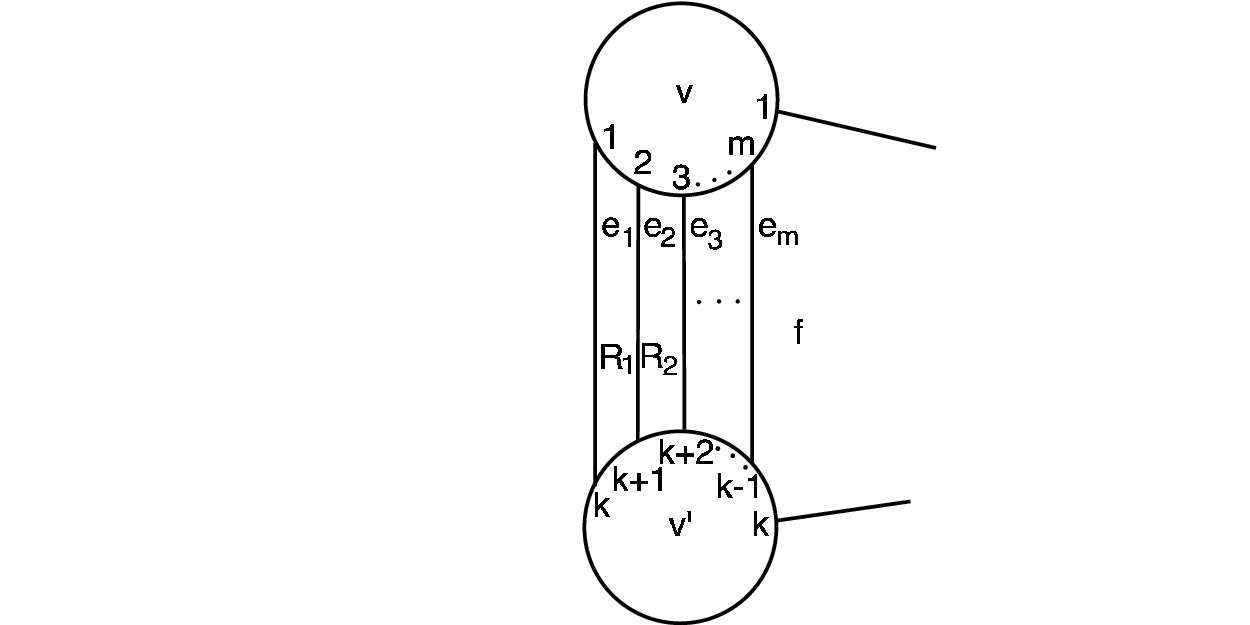}}
\caption{A pair of vertices of opposite orientations connected by $m$ parallel edges}
\end{figure}

Let $d$ be the number of orbits of the action of the
permutation $\s$ on the set $\{b_1,...,b_m\}$,
each of $m/d$ elements.
We may assume that indices are given as shown in Figure 2. By the
parity rule, the index $k$ in Figure 2 must be an odd number. From
Figure 2, we see obviously that $b_1$ and $b_k$ are in the same
orbit, and $b_m$ and $b_{k-1}$ are in another orbit. By Corollary
\ref{Length4},  $\Phi_1^+$ is a pair of tight components, $T_1$
and $T_2$, which include all boundary components of $F$.

\begin{lemma}\label{orbit} Suppose that $e_1,..., e_{m}$ are  $m$
  consecutively parallel edges of $\G$ connecting two vertices $v$ and
  $v'$ with different orientations. We may assume that the permutation
  $\s$ given in the preceding paragraph is as shown in Figure 2. Then
  $b_1\cup b_k$ and $b_{k-1}\cup b_m$ are contained in different
  components of $\Phi_1^+$; i.e. one in $T_1$ and the other in $T_2$.
\end{lemma}

\begin{proof}
 Recall that the annulus $B_{k-1}$ in $\p M$ has boundary
$b_{k-1}\cup b_k$ and the annulus $B_m$ has boundary $b_m\cup b_1$,
both contained in $X^+$. Thus $b_{k-1}$ and $b_k$ are contained in
different components of $\Phi_1^+$, and so are $b_{1}$ and $b_m$. In
particular the conclusion of the lemma follows immediately if $k=1$,
i.e. if the permutation $\s$ is trivial. So we may assume that $\s$ is
non-trivial, $b_1$ is contained $T_1$, and $b_m$ in $T_2$. We now only
need to show that $b_k$ is in $T_1$. Since $k\ne 1$, the bigon face
$R_{k-1}$ is mapped into $X^+$ (since $k$ is odd) and $e_k^*$ connects
$b_k$ to $b_{\s(k)}$. If $\s(k)=1$, then we are done. If $\s(k)\ne 1$,
then $R_{\s(k)-1}^*$ is in $X^+$ (since $\s(k)$ is odd) and
$e^*_{\s(k)}$ connects $b_{\s(k)}$ to $b_{\s^2(k)}$ (recall that the
indices here are defined mod (m)).  Repeat in this way for finitely
many times until $\s^n(k)=1$ for some positive integer $n$ (actually
$n=\frac{m}{d}-1$ is the number of elements in the orbit minus one).
\end{proof}

Let $\bar\G$ denote the reduced graph of $\G$, obtained from $\G$ by
amalgamating parallel edges into a single edge. Then $\bar\G$ is a
graph with no $1$-edge or $2$-edge disk faces. If an edge $\bar e$ of
$\bar\G$ represents $n$ parallel edges of $\G$, we say the edge $\bar
e$ has weight $n$.

\begin{lemma}\label{Corner}
  Let $\bar e$ be an edge of $\bar \G$ with weight $m$. Then no face
  $f$ of $\bar \G$ with $\p f$ containing $\bar e$ is an triangle
  face, i.e.  a 3-edge disk face.
\end{lemma}

\begin{proof}
  Suppose otherwise that there is a triangle face f in $\bar\G$
whose boundary contains $\bar e$. Note that $f$ is also a face in
the graph $\G$. Let $v$ and $v'$ be the two vertices connected by
$\bar e$, and $e_1,...,e_m$ be the family of parallel edges of
$\G$ represented by $\bar e$.

First we consider the case that $v$ and $v'$ have different
orientations. We may assume that the labels of the endpoints of
the edges $e_i$ are given as in Figure 2.  Now consider the two
``corners'' of the face $f$ at the vertices $v$ and $v'$, i.e. the
intersection arcs of the boundary of $f$ with the boundary of the
fat vertices. From Figure 2, we see that the two corners have
labels $k, k-1$ and $m, 1$ respectively. If we follow $\p f$ in
the clockwise direction, the four labels appear in the order $k,
k-1, m, 1$.

Recall the setting in Lemma \ref{homology}. The punctured solid torus
$U$ carried a generator of the first homology of the punctured lens
space $\hat X^+$ and the disk $D_{1/2}$ was a meridian disk of $U$.
Note that $U$ contains all the 2-handles $H_i$ for $i$ even and that
the disk $D_{1/2}$ intersects each $H_i$, $i$ even, in a single
meridian disk of $H_i$.

If a disk face $f'$ of $\G$ has $n$ corners and is on the $X^+$ side,
then it is not hard to see that $(\p f')^*$ is contained in $U$ and
intersects $D_{1/2}$ transversely in $n$ points (all the intersections
occur precisely one each within the corner arcs of $\p f'$). In our
current situation, the triangle face $f$ is indeed on the $X^+$ side
since $R_{m-1}$ is on the $X^-$ side.  Further the algebraic
intersection number of $\p f$ with $D_{1/2}$ is $1$ or $-1$ because of
the label orders on $\p f$ together with Lemma \ref{orbit}.  Now we
see that the existence of such a triangle face $f$ implies that the
first homology of $\hat X^+$ is trivial, contradicting to the fact
that $\hat X^+$ is a punctured non-trivial lens space.  This completes
the proof for the case when the vertices $v$ and $v'$ have different
orientation.

Now we consider the case when $v$ and $v'$ have the same orientation.
Then by Lemma \ref{NoSameOri}, we have $m=4$. By Lemma \ref{Scycle}
(2), Lemma \ref{No+Scycle} and Lemma \ref{NoEScycle}, we may assume
that four edges form two $S$-cycles and the labels on the tails and
heads of the four edges are as shown in Figure 3. From the figure, we
see directly that $b_1\cup b_2$ is contained in one the tight
components of $\Phi_1^+$ and $b_3\cup b_4$ is contained in the other.
And the labels $2,3,4,1$ appeared consecutively on $\p f$ (in
clockwise direction). This would imply, as in the previous case, that
the first homology of $\hat X^+$ is trivial, giving the same
contradiction. \end{proof}

\begin{figure}
\centerline{\includegraphics{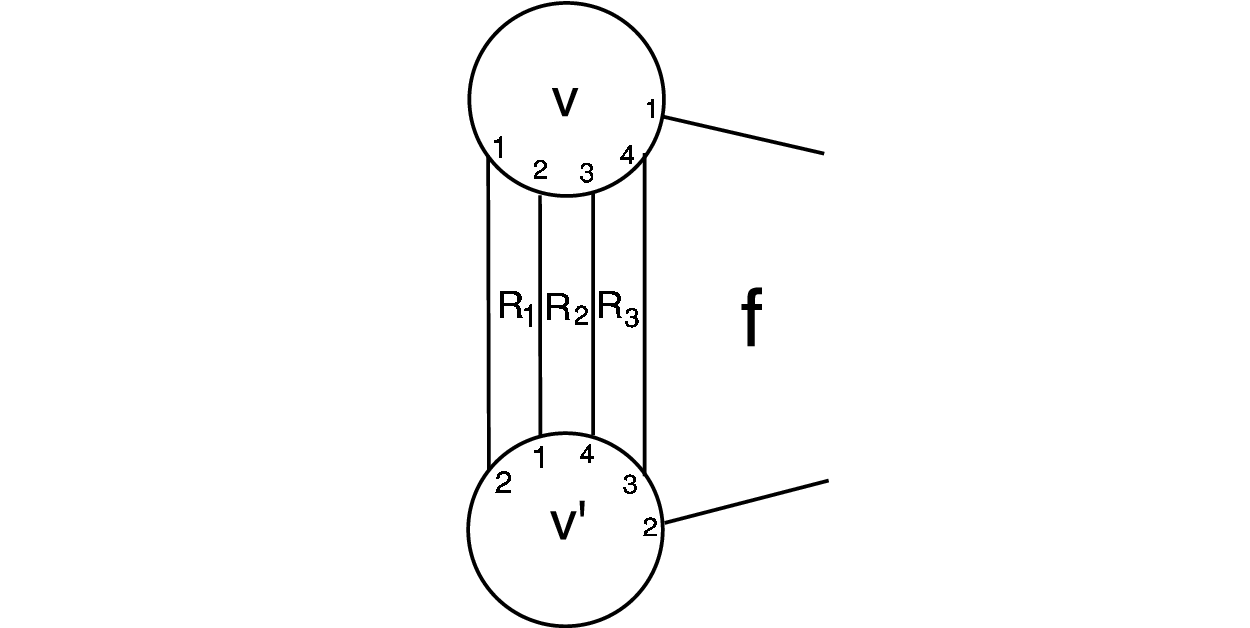}}
\caption{A pair of
vertices of the same orientation connected by $m=4$ parallel edges
which form two $S$-cycles}
\end{figure}

Now we can finish the proof of Proposition \ref{TLS}. Suppose
otherwise that $\D(\alpha,\beta)\geq 5$.  An Euler characteristic
calculation shows that either the reduced graph $\bar \G$ on the torus
$T$ has a vertex of valence less than $6$ or every vertex of $\bar\G$
has valence $6$ and every face of $\bar\G$ is a triangle face.  By
Lemma \ref{Length}, every edge of $\bar\G$ has weight at most $m$. So
we may assume that $\bar\G$ has no vertex of valence less than $5$.

Consideration of Lemma \ref{Length} yields $\D=5$ and the fact that
every edge of $\bar\G$ incident a valence $5$ vertex has weight
exactly equal to $m$.  So it follows from Lemma \ref{Corner} that the
graph $\bar\G$ in the torus $T$ has the following properties:
\begin{itemize}
\item no disk face has $1$ or $2$ edges,
\item every vertex has valence at least $5$,
\item no triangle face is incident to a valence $5$
vertex.
\end{itemize}
For a vertex $v$ and face $f$ of $\bar \Gamma$, we write $v \in
\partial f$ to signify that $v$ is incident to $f$.  Consider
$$\chi_f = \chi(f) + \sum_{v \in \partial f} (\frac{1}{\hbox{valency}(v)} - \frac12).$$
By construction, if $\partial f$ has three edges then $\hbox{valency}(v) \geq 6$ for each $v \in \partial f$. Hence $\chi_f \leq 1 + 3(-\frac13) = 0$ with equality if and only if $f$ is a triangle face and each of its vertices has valency $6$. On the other hand, if $\partial f$ has at least four edges, then $\chi_f \leq 1 + 4(\frac15 - \frac12) = -\frac15 < 0$. Thus since $0 = \chi(T) = \sum_f \chi_f$, each face of $\bar \Gamma$ is a triangle face and vertex has valency $6$.

The proof of the following
lemma is similar to that of Lemma \ref{Corner}.

\begin{lemma}\label{Corner2}
  The graph $\bar \G$ cannot have an edge with weight larger than
  $m-2$.
\end{lemma}

\begin{proof}
 Suppose otherwise that $\bar e$ is an edge with weight at least
$m-1$. Since every face of $\bar\G$ is a triangle face, Lemma
\ref{Corner} shows that the weight of $\bar e$ is exactly $m-1$. Let
$v$ and $v'$ be the two vertices connected by $\bar e$ and let
$e_1,...,e_{m-1}$ be the family of parallel edges of $\G$ represented
by $\bar e$, oriented such that their tails are at $v$ and heads at
$v'$.

If $v$ and $v'$ have the same orientation, then we have $m=4$ (Lemma
\ref{NoSameOri}). Since $v$ has valence $6$ while there are
$4\D(\alpha,\beta)\geq 20$ endpoints of edges of $\G$ incident to $v$,
some edge of $\bar \G$ incident to $v$ will have weight $m=4$,
contrary to the conclusion of Lemma \ref{Corner}.

Suppose then that $v$ and $v'$ have different orientations. We may
assume that the labels of the endpoints of the edges $e_i$ are given
as in Figure 4. By the parity rule, for each of the edges
$e_1,...,e_{m-1}$, the two labels at its endpoints are congruent (mod
2). In particular the label $k$ in Figure 4 is an odd number.  Denote
by $R_1,..., R_{m-2}$ the $m-2$ bigon faces defined by the $m-1$
edges, where $R_j$ contains the edges $e_j$ and $e_{j+1}$.  By our
convention $R_1$ lies on the $X^-$ side since $R_1^*$ intersects the
annulus $B_1$ which lies on the $X^-$ side (cf.  Figure 4).  So the
triangle face $f$ of $\bar\G$ which contains the edge $e_1$ (shown in
Figure 4) lies on the $X^+$ side. It also follows that for each
$i=1,..., m-2$, $R_i$ lies on the $X^-$-side if $i$ is odd or on the
$X^+$ side if $i$ is even.

By Corollary \ref{Length4}, $\Phi_1^+$ is a pair of tight components
$T_1$ and $T_2$ and contains all components of $\p F$. We want to show
that $b_{1}\cup b_{k}$ is contained in one component of $\Phi_1^+$ and
$b_m\cup b_{k-1}$ is contained in the other. This is obviously true if
$k=1$ since $B_m$ is contained in $\Sigma_1^+$. So suppose that $k>1$.
As in the proof of Lemma \ref{Corner}, by considering the orbit of the
label $1$ under the permutation of odd integers $\{1,3,5, \ldots ,
m-1\}$ given by the $m-1$ edges, we see that there is a sequence of
odd labels $k_1=k, k_2,..., k_n\in \{3, 5, ..., m-1\}$ and edges
$e_{i_1},...,e_{i_n}\in \{e_3, e_5, ...,e_{m-1}\}$ such that for
$1\leq j<n$, the edge $e_{i_j}$ has tail label $k_j$ and head label
$k_{j+1}$, and the edge $e_{i_n}$ has tail label $k_n$ and head label
$1$.  Since $R_{i_j-1}$ lies on the $X^+$ side, we see that all
$e_{i_j}^*$, $j=1,...,n$, are contained $\Phi_1^+$.  Since these $n$
edges $e^*_{i_1}$, ..., $e_{i_n}^*$ connect $b_{k}=b_{k_1}$,
$b_{k_2}$,...  $b_{k_n}$ and $b_1$, we see that $b_{1}\cup b_{k}$ is
contained in one component of $\Phi_1^+$, say $T_1$.  It follows that
$b_m\cup b_{k-1}$ is contained in $T_2$, the other component of
$\Phi_1^+$, since the annuli $B_m$ and $B_{k-1}$ are contained in
$\Sigma_1^+$.

From Figure 4, we see that the two corners of $f$ at $v$ and $v'$ have
labels $m, 1$ and $k,k-1$ respectively in clockwise direction. Now
combining with Lemma \ref{homology}, we see that the first homology of
$\hat X^+$ is trivial, which is a contradiction.  \end{proof}

\begin{figure}
\centerline{\includegraphics{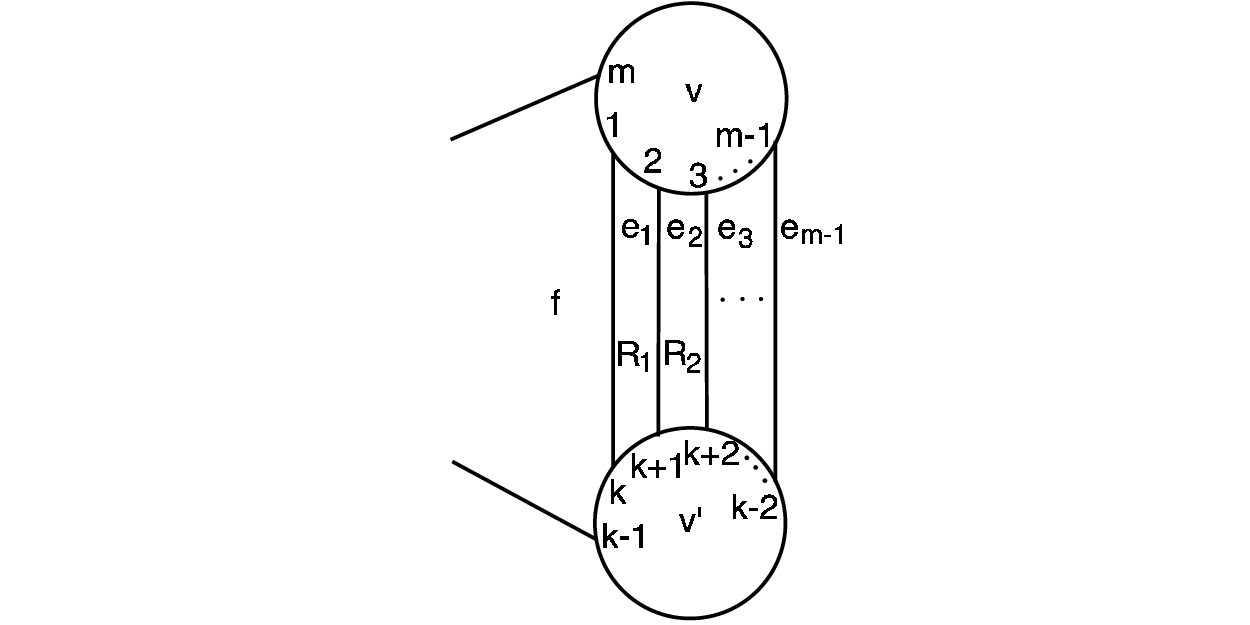}}
\caption{A pair
of vertices of different orientations connected by $m-1$ parallel
edges}
\end{figure}

We call an edge of $\bar\G$  positive (respectively negative) if
it connects two vertices of the same orientation (respectively
different orientations). We call the endpoint of an edge at a
vertex positive or  negative if the edge is positive or negative.
We define the weight of an endpoint of an edge to be the weight of
the  edge. The sum of the weights of the endpoints at any vertex
is $\D(\alpha,\beta) m$.

\begin{lemma}\label{NoTwoPositive}
  Let $v$ be a vertex of $\bar\G$.  Then among the six endpoints at
  $v$, at most one is positive.
\end{lemma}

\begin{proof}
 If there are two positive endpoints at $v$, then their weight
sum is at most $m+4$ by Lemmas \ref{NoEScycle} and \ref{Scycle}.
So the rest four endpoints have total weight at least $4m-4$. So
at least one endpoint has weight $m-1$. This gives a contradiction
with Lemma \ref{Corner2}. \end{proof}

\begin{lemma}\label{TwoPositive}
  There is a vertex of $\bar\G$ with at least two positive endpoints.
\end{lemma}

\begin{proof}
 The previous lemma implies that the graph $\bar \G$ has no loops.
Pick any vertex $v_0$ of $\bar\G$ and let $p_1,...,p_6$ be the six
endpoints at $v$ in clockwise order. We may assume that $p_5,...,p_6$
are all negative endpoints. Let $\bar e_i$ be the edge of $\bar\G$
with endpoint $p_i$ and observe that they are distinct edges since
there are no loops.  Let $v_i$ be the other vertex that $\bar e_i$ is
incident to. Then the $v_i$ have the same orientation for $i=2,...,6$.
Now $v_5\ne v_4$ since there are no loops, and there is an edge of
$\bar \G$ connecting them. Similarly there is an edge connecting $v_3$
and $v_4$. Note that $v_3\ne v_5$ since otherwise there is either a
non-triangle face of $\bar\G$ or $v_4$ has valence less than $6$. Thus
$v_4$ has at least two positive endpoints. \end{proof}

\begin{proof}[Proof of Proposition \ref{TLS}]. 
The contradiction between Lemmas \ref{NoTwoPositive} and
\ref{TwoPositive} completes the proof that $\Delta(\alpha, \beta) \leq
4$. If we have equality, then Proposition 8.4 of [BCSZ] and Theorem
\ref{redfin} imply that $M(\alpha)$ is Seifert fibred with base
orbifold of the form $S^2(r,s,t)$ where $(r,s,t)$ is a hyperbolic
triple and $\hbox{lcm}(r,s,t)$ divides $4$. Thus Proposition \ref{TLS}
holds.  \end{proof}

\vspace{8mm}
\noindent
{\large {\bf References}}

{\small 
\noindent[BB] L. Ben Abdelghani and S. Boyer,
{\it A calculation of the Culler-Shalen seminorms associated to small
  Seifert Dehn fillings}, Proc. Lond. Math. Soc.  {\bf 83} (2001)
235-256.

\noindent[B] S. Boyer,
{\it On the local structure of $SL_2(\mathbb C)$-character varieties
  at reducible characters}, Top. Appl. {\bf 121} (2002), 383-413.

\noindent[BCSZ] S. Boyer, M. Culler, P. Shalen and X. Zhang,
{\it Characteristics subsurfaces and Dehn filling}, Trans. Amer. Math.
Soc., to appear.

\noindent[BGZ] S. Boyer, C. McA. Gordon, and X. Zhang, {\it Dehn
  fillings of large hyperbolic 3-manifolds}, J. Diff. Geom., {\bf 58}
(2001) 263-308.

\noindent[BN] S. Boyer and A. Nicas, {\it Varieties of group
  representations and Casson's invariant for rational homology
  3-spheres}, Trans. Amer. Math. Soc.  {\bf 322} (1990), 507-522.

\noindent [BZ1] S. Boyer and X. Zhang, {\it Finite Dehn surgery on
  knots}, J. Amer. Math. Soc. {\bf 9} (1996) 1005-1050.

\noindent [BZ2] ---, {\it On Culler-Shalen Seminorms and Dehn filling},
  Ann. Math. {\bf 148} (1998) 737-801.

\noindent [BZ3] ---, {\it Virtual Haken 3-manifolds and Dehn filling}, Topology {\bf 39} (2000), 103-114. 

\noindent [BZ4] --- S. Boyer and X. Zhang,
{\it On simple points of the character varieties of $3$-manifolds},
{\bf Knots in Hellas '98}, Proceedings of the International Conference
on Knot Theory and its Ramifications, World Scientific Publishing Co.
Pte. Ltd., 2000.

\noindent [BZ5] ---, {\it A proof of the finite filling conjecture},
J. Diff. Geom., {\bf 59} (2001) 87-176.

\noindent
[CJ] A. Casson and D. Jungreis, {\it Convergence groups and Seifert
  fibred $3$-manifolds}, Invent. Math. {\bf 118} (1994), 441-456.

\noindent [CGLS] M. Culler, C. Gordon, J. Luecke and P. Shalen,
{\it Dehn surgery on knots}, Ann. of Math. {\bf 125} (1987) 237-300.

\noindent [CS] M. Culler and P. Shalen, {\it Varieties of group
  representations and splittings of 3-manifolds}, Ann. of Math. {\bf
  117} (1983) 109-146.

\noindent [Dn] N. Dunfield,
{\it Cyclic surgery, degrees of maps of character curves, and volume
  rigidity for hyperbolic manifolds}, Inv. Math. {\bf 136} (1999),
623-657.

\noindent [Ga1]
D. Gabai, {\it Foliations and the topology of $3$-manifolds II}, J.
Diff. Geom. {\bf 26} (1987) 461-478.

\noindent [Ga2] ---, {\it Homotopy hyperbolic 3-manifolds are virtually
  hyperbolic}, J. Amer. Math. Soc.  {\bf 7} (1994), 193-198.

\noindent [Ga3]   ---,
{\it Convergence groups are Fuchsian groups}, Ann. Math. {\bf 136}
(1992), 447-510.

\noindent [GMT] D. Gabai, R. Meyerhoff and N. Thurston,
{\it Homotopy hyperbolic 3-manifolds are hyperbolic}, preprint. 

\noindent [Gld]
W. Goldman, {\it The symplectic nature of fundamental
groups of surfaces}, Advances in Math. {\bf 54} (1984), 200-225.

\noindent [Go]{}
C. McA. Gordon, {\it Dehn filling: A survey}, Proceedings of the
Mini Semester in Knot Theory, Banach Center, Warsaw, Poland, 1995.

\noindent [GLu] C. McA. Gordon and J. Luecke, {\it Reducible manifolds
and Dehn surgery}, Topology {\bf 35} (1996), 385--409.

\noindent[H] J. Hempel, 3-Manifolds, Annals of Mathematics Studies 86,
Princeton Univ. Press 1976.

\noindent[Ja] W. Jaco, {\it lectures on three-manifold topology},
CBMS CBMS Regional Conf. Ser. Math. {\bf 43} (1980).

\noindent[Lee] S. Lee,
{\it Reducing and toroidal Dehn fillings on $3$-manifolds bounded by
  two tori}, Math. Res. Lett. {\bf 13} (2006), 287-306.

\noindent[NR] W. Neumann and F. Raymond,
{\it Seifert manifolds, plumbing, $\mu$-invariant and orientation
  reversing maps}, in Algebraic and geometric topology, Proc. Sympos.,
Univ.  California, Santa Barbara, Calif., 1977, pp. 163--196, Lecture
Notes in Math., 664, Springer, Berlin, 1978.

\noindent 
[Ne] P. E. Newstead, {\it Introduction to moduli problems and orbit
  spaces}, Tata Institute Lecture Notes, Narosa Publishing House, New
Delhi, 1978.

\noindent[Oh] S. Oh, {\it Reducible and toroidal $3$-manifolds obtained
  by Dehn fillings}, Topology Appl. {\bf 75} (1997), 93--104.

\noindent[Sc] P. Scott, {\it A new proof of the annulus and torus theorems},
Amer. J. Math. {\bf 102} (1980), 241--277.

\noindent[Ser] J.-P. Serre,
{\it Repr\'esentations lin\'eaires des groupes finis}, 3e ed.,
Hermann, Paris, 1978.

\noindent [Wa] F. Waldhausen, {\it On irreducible $3$-manifolds which
are sufficiently large}, Ann. Math. {\bf 87} (1968), 56-88.

\noindent [Wu1]{} -----, {\it Incompressibility of surfaces in surgered
  $3$-manifolds}, Topology {\bf 31} (1992), 271-279.

\noindent[Wu2] Y. Wu, {\it Dehn fillings producing reducible manifolds
and toroidal manifolds}, Topology {\bf 37} (1998), 95-108.

\noindent [Wu3] Y. Wu, {\it Standard graphs in lens spaces},
Pacific J. Math. (to appear). Preprint available at
http://arXiv.org/abs/math/0011007.
}

\end{document}